\documentclass[10pt]{article}
\usepackage{amssymb,amsmath,amsthm}
\usepackage{latexsym}

\input amssym.def
\input amssym

\newtheorem{theorem}{Theorem}[section]
\newtheorem{propo}[theorem]{Proposition}
\newtheorem{theo}[theorem]{Theorem}
\newtheorem{defi}[theorem]{Definition}
\newtheorem{corol}[theorem]{Corollary}
\newtheorem{remark}{Remark}

\renewcommand{\theequation}{\thesection.\arabic{equation}}

\def\sect#1{\section{#1}\setcounter{equation}{0}\setcounter{remark}{0}}
\def\ssect#1{\subsection{#1}}

\newcommand{\bc}{\begin{center}}
\newcommand{\ec}{\end{center}}
\newcommand{\ba}{\begin{array}}
\newcommand{\ea}{\end{array}}
\newcommand{\beq}{\begin{equation}}
\newcommand{\eeq}{\end{equation}}
\newcommand{\beqa}{\begin{eqnarray}}
\newcommand{\eeqa}{\end{eqnarray}}
\def\bi{\begin{itemize}}
\def\ei{\end{itemize}}
\newcommand{\no}{\nonumber}
\newcommand{\nn}{\nonumber \\ }
\renewcommand{\d}{\partial}
\def\frc#1#2{\frac{#1}{#2}}
\def\b#1{\bar{#1}}
\def\t#1{\tilde{#1}}
\def\h#1{\hat{#1}}
\def\lt#1{\left#1}
\def\rt#1{\right#1}
\newcommand{\<}{\langle}
\renewcommand{\>}{\rangle}

\newcommand{\Z}{{\Bbb Z}}
\newcommand{\C}{{\Bbb C}}
\newcommand{\N}{{\Bbb N}}

\newcommand{\D}{{\cal D}}
\def\a#1{{\goth{#1}}}
\newcommand{\End}{{\rm End\,}}

\newcommand{\Res}{{\rm Res}}
\newcommand{\Tr}{{\rm Tr}}
\newcommand{\halmos}{\rule{1ex}{1.4ex}}
\newcommand{\eproof}{\hspace*{\fill}\mbox{$\halmos$}}

\def\for#1{\ \ \ \ \mbox{#1 }\ \ {}}
\def\com#1{\quad(#1)}
\newcommand{\om}{\omega_p}

\renewcommand{\:}{\mbox{\tiny ${\bullet\atop\bullet}$}}
\newcommand{\nordbullet}{\mbox{\tiny ${\bullet\atop\bullet}$}}
\renewcommand{\aa}{\mbox{\scriptsize ${\ast\atop\ast}$}}
\newcommand{\pp}{\mbox{\scriptsize ${+ \atop +}$}}
\newcommand{\xx}{\mbox{\scriptsize ${\times\atop\times}$}}

\def\dpf#1#2{\delta\left(\frac{#1}{#2}\rt)}

\title{Twisted vertex operators and Bernoulli polynomials}

\author{B. Doyon,  J. Lepowsky, A. Milas}

\date{}

\begin{document}

    \maketitle

\begin{abstract}

Using general principles in the theory of vertex operator algebras
and their twisted modules, we obtain a bosonic, twisted construction
of a certain central extension of a Lie algebra of differential
operators on the circle, for an arbitrary twisting automorphism. The
construction involves the Bernoulli polynomials in a fundamental
way. We develop new identities and principles in the theory of
vertex operator algebras and their twisted modules, and explain the
construction by applying general results, including an identity that
we call {\em modified weak associativity}, to the Heisenberg vertex
operator algebra. This paper gives proofs and further explanations
of results announced earlier. It is a generalization to twisted
vertex operators of work announced by the second author some time
ago, and includes as a special case the proof of the main results of
that work.

\end{abstract}
{\em Keywords:}{ \ Vertex operator algebras; twisted modules;
Bernoulli polynomials} \\
\noindent Mathematics Subject Classification 2000: Primary: 17B69,
Secondary: 11B68, 17B66
\tableofcontents

\renewcommand{\theequation}{\thesection.\arabic{equation}}
\renewcommand{\thetheorem}{\thesection.\arabic{theorem}}
\setcounter{equation}{0} \setcounter{theorem}{0}
\setcounter{section}{0}

\sect{Introduction}

\subsection{Twisted vertex operators and twisted modules}

The theory of vertex operator algebras and their twisted modules
is a powerful theory with applications in many areas of
mathematics and physics. In mathematics, vertex operators, twisted
and untwisted, are involved in a fundamental way in the
construction of modules for many infinite-dimensional algebras,
including affine Lie algebras and central extensions of algebras
of differential operators, for instance, the Virasoro algebra. In
particular, twisted vertex operators were discovered and used in
\cite{LW}, where a formal--differential--operator construction of
the affine Lie algebra $A_{1}^{(1)}$ was found. Moreover, various
types of vertex operators, including twisted vertex operators, and
relations between them, were used in \cite{FLM3} to construct the
``moonshine module'' for the Monster group, the largest sporadic
finite simple group, in such a way that the Monster was realized
as the automorphism group of this ``algebra of vertex operators.''
In physics, vertex operator algebras form the structure at the
foundation of conformal field theory \cite{BPZ}.  What came to be
understood as twisted modules for vertex operator algebras are the
main building blocks for the so-called {\em orbifold models} in
conformal field theory (see \cite{DHVW1}, \cite{DHVW2}, \cite{DFMS}
and \cite{DVVV}).  In particular, the moonshine module \cite{FLM1},
\cite{FLM3} was the first example of an orbifold construction -- it
appeared before the concept of orbifold compactification was developed
in string theory.

The precise notion of {\em vertex operator algebra} arose {}from
formalizing the natural algebraic structure underlying (untwisted)
vertex operators, which are objects based on a lattice (cf.
\cite{Bo}, \cite{FLM3}), at least when the lattice is even. For
the basic theory, based on formal calculus and the Jacobi
identity, of vertex operator algebras, including lattice vertex
operator algebras, the reader can consult \cite{FLM3} and
\cite{LL}. Twisted vertex operators, parametrized by the vectors
in a space having a structure of vertex operator algebra together
with a lattice isometry of finite order, led to the notion of {\em
twisted module} for a vertex operator algebra. A systematic
treatment of certain ``twisted formal calculus'' underlying
various vertex operator constructions associated to an even
lattice was first carried out in \cite{L1}. A more complete study
of twisted vertex operators was then carried out in \cite{FLM2},
where in particular the formal operator $e^{\Delta_x}$ was
introduced as an essential part of the definition and as the main
addition compared to the construction of untwisted vertex
operators (see also \cite{FLM3}, \cite{DL2}). The definition of
``twisted module'' for vertex operator algebras was introduced in
\cite{FFR} (see also \cite{D}); this notion summarizes the key
properties that had been discovered for twisted vertex operators,
in particular, the twisted Jacobi identity. A large part of the
theory of twisted modules was further developed in \cite{Li2}.
More recently, a ``coordinate-free'' approach was developed in
\cite{FrS} along the lines of \cite{FrB} (cf. \cite{H2}). It is
important to mention here the work \cite{BDM} where, by using the
theory of twisted modules, a conceptual explanation of the
so-called {\em cyclic orbifold} construction in conformal field
theory (cf. \cite{BHS}) was given.  In addition to \cite{BHS}, we
would like to mention the program developed in \cite{dBHO}, \cite{HO},
\cite{GHHO} and \cite{HH} including studies of twisted vertex operators
associated with twisted abelian current algebras.

The present paper contributes to the general theory of twisted
modules and to the deeper study of twisted vertex operators in
important special settings as well; see Section
\ref{sectpresentwork} below.  New general principles in vertex
operator algebra theory are introduced and developed, and then are
applied to twisted modules and twisted vertex operators associated
with lattice vertex operator algebras, to shed light on the
relations between certain values of the Riemann zeta function and
certain algebras of differential operators on the circle. The
results in this paper were announced some time ago in \cite{L4}
and \cite{L5} in the untwisted case, and recently in \cite{DLMi}
in the more general, and more subtle, twisted case.  In
particular, the present paper is the first giving the full
treatment of even the untwisted case of these ideas, methods,
results and proofs.  It also happens that the new ideas introduced
here (and announced in \cite{DLMi}) enable us to further elucidate
even the earlier untwisted case of \cite{L4} and \cite{L5}. Our
main results are the following: theorems on the new and useful
concept of ``resolving factors'' (Theorems \ref{theovwa} and
\ref{theomwa} in the untwisted case and Theorems \ref{theovwat},
\ref{theomwat} and \ref{theomwaht} in the more general twisted
setting; see Remarks \ref{remares} and \ref{remaresh}); the new
general commutator formula Theorem \ref{commiter}; and the main
result in the special context incorporating differential operators
on the circle, Theorem \ref{main1}.

Throughout this paper, we extensively use the framework and
technology of formal calculus (cf. \cite{FLM3}), and in fact it is
necessary to use formal calculus both to express and to prove our
results.  From the point of view of a physicist, this
formal-calculus framework for vertex operator algebra theory is
revealing itself as the most appropriate way of both rigorously
and effectively studying the chiral sector in conformal field
theory.  Our results and proofs, in both the general and special
settings in this paper, in fact include new applications of formal
calculus; this serves to emphasize its efficiency as a tool in
uncovering surprising structures in conformal field theory and in
addition, exploits the intuition that one gets from studying
conformal field theory as a physicist.  Specifically, for
instance, our twisted construction in the special context
mentioned above is highly nontrivial, even though it involves only
a small part of the general identities that we establish; what we
call the ``correction'' terms (involving Bernoulli polynomials),
crucial ingredients of our constructions (as we explain below),
are in general very hard to compute, and our formal-calculus
methods provide an extremely efficient way of computing them.

We have tried to be careful to make this paper reasonably
self-contained and accessible even to non-experts in vertex
operator algebra theory.  In the literature, there are a number of
useful variants of certain basic definitions and concepts, and
because of the subtle nature of the results and proofs presented
here, we have found it appropriate to record the precise
background and the exact versions of a number of basic notions
that we need.

\subsection{Twisted vertex operators and Bernoulli polynomials}

In this subsection, we shall present some background and
motivation for this work and in particular, we shall discuss how
new ideas in vertex operator algebra theory enable one to shed
light on properties of zeta values and Bernoulli polynomials. Some
of the following motivation is extracted from our earlier research
announcements \cite{L4}, \cite{L5} and \cite{DLMi} of the present
work. We have decided to recall this motivation here in order to
make the present paper more self-contained, since it is the first
paper to supply all the proofs of the announced results.

In order to understand the basic ideas behind our work, consider
the famous classical ``formula''
\begin{equation}\label{1+}
1 + 2 + 3 + \cdots = - \frac{1}{12},
\end{equation}
which has the rigorous meaning
\begin{equation}\label{zeta-1}
\zeta(-1) = - \frac{1}{12}.
\end{equation}
Here $\zeta$ is the Riemann zeta function
\begin{equation}\label{zeta}
\zeta(s) = \sum_{n > 0} n^{-s}
\end{equation}
(analytically continued), and (\ref{1+}) is classically
generalized by the formal equality
\begin{equation}\label{zeta-s}
\sum_{n>0}n^s = \zeta(-s)
\end{equation}
for $s = 0,1,2,\dots$.  It is well known that the classical number
theory underlying this analytic continuation is related to the
issue of regularizing certain infinities in quantum field theory,
in particular, in conformal field theory. One of the goals of the
present work is to develop some general principles of vertex
operator algebra theory that elucidate the passage {}from the
unrigorous but suggestive formula (\ref{1+}) to formula
(\ref{zeta-1}), and the generalization (\ref{zeta-s}), as
announced in \cite{L4} and \cite{L5}. This passage was noticed by
S. Bloch \cite{Bl} to have simplifying effects on certain infinite
dimensional Lie algebras of formal differential operators on the
circle, and his work was the main motivation behind the ideas
initially announced in \cite{L4} and \cite{L5}. Developing further
the theory of twisted vertex operators, as announced in
\cite{DLMi}, we will also generalize this to analogous formal
relations involving the Hurwitz zeta function
\beq\label{hurdef}
    \zeta(s,v) = \sum_{n \ge 0} (n+v)^{-s},
\eeq
that is, to formal unrigorous relations of the type
\beq\label{hurwitzformal}
    \sum_{n\ge0} (n+v)^s = \zeta(-s,v)
\eeq
for $v$ certain rational numbers and for $s=0,1,2,\ldots$. In terms of
the Lie algebras studied by Bloch, this turns out to be related to
representations of these Lie algebras on certain twisted spaces that
we discuss in this paper.  We will hence ``explain'' and considerably
generalize this work of Bloch's, and in the process, discover new
structures in the theory of vertex operator algebras. In fact, in
\cite{DLMi} new developments were announced that are of interest for
this ``explanation'' of Bloch's work and more generally in the theory
of vertex operator algebras both in the twisted and untwisted
settings. In the following we will review this work of Bloch's and
then explain the meaning of our results in this context.

Before reviewing this work, though, let us recall, from \cite{L4}, a
variant of Euler's heuristic interpretation of (\ref{1+}) and its
generalization (\ref{zeta-s}). Consider the following formal
expansion (in powers of $x$), which defines the Bernoulli numbers
$B_k$ for $k\ge0$:
\begin{equation}\label{fromBer}
\frac{1}{1 - e^x} = - \sum_{k \geq 0} \frac{B_{k}}{k!} x^{k-1}.
\end{equation}
Expand the left--hand side {\it unrigorously} as the formal
geometric series
\begin{equation}\label{unrig}
1 + e^x + e^{2x} + \cdots = 1 + \sum_{k \geq 0} \frac{1^k}{k!} x^k
+ \sum_{k \geq 0} \frac{2^k}{k!} x^k + \cdots.
\end{equation}
For $k>1$, the coefficient of $x^{k-1}$ in this formal expression
is
\begin{equation}
\frac{1}{(k-1)!} (1^{k-1} + 2^{k-1} + \cdots),
\end{equation}
which formally resembles $\frac{1}{(k-1)!} \zeta (-k+1)$.  Also,
the coefficient of $x^0$ in (\ref{unrig}) is formally
\beq\label{casek=1}
1 + \frac{1}{0!} (1^0 + 2^0 + \cdots),
\eeq
which we formally view as $1 + \zeta (0)$ (and not as $\zeta
(0)$). Thus, formally equating the coefficients of $x^l$ for $l
\geq 0$ in (\ref{fromBer}) tells us what (\ref{zeta-s}) means and
``explains'' the well--known relation between Bernoulli numbers
and the zeta function
\begin{equation}\label{zetaB}
\zeta(-k+1) = - \frac{B_{k}}{k}
\end{equation}
for $k>1$, and the fact that
\beq\label{zeta0B1}
\zeta(0) = - B_1 - 1\ \ \left(= - \frac{1}{2}\right).
\eeq
Note, incidentally, that the coefficient of $x^{-1}$ in
(\ref{fromBer}), which is $-B_0 = -1$, does not appear in
(\ref{unrig}) -- (\ref{zeta0B1}); our unrigorous equating of
coefficients of powers of $x$ in the two series (\ref{unrig}) and
(\ref{fromBer}) only applies to nonnegative powers.

Similar considerations apply to the formal relation
(\ref{hurwitzformal}): Let $v$ be a complex number and consider
the formal expansion (again in powers of $x$) which defines the
Bernoulli polynomials $B_k(v)$ for $k\ge0$:
\begin{equation}\label{fromBerHur}
\frac{e^{vx}}{1 - e^x} = - \sum_{k \geq 0} \frac{B_{k}(v)}{k!}
x^{k-1}.
\end{equation}
Note that
\beq
    B_k(0) = B_k.
\eeq
Expanding the geometric series on the left--hand side unrigorously
in powers of $e^x$ and multiplying by $e^{vx}$ formally gives
\begin{equation}\label{unrig-Hurwitz}
e^{vx} + e^{(v+1)x} + e^{(v+2)x} + \cdots = \sum_{k\geq0}
\frc{v^k}{k!} x^k + \sum_{k \geq 0} \frac{(v+1)^k}{k!} x^k +
\sum_{k \geq 0} \frac{(v+2)^k}{k!} x^k + \cdots.
\end{equation}
Now, for $k\geq1$, the coefficient of $x^{k-1}$ in this formal
expression is
\begin{equation}
\frac{1}{(k-1)!} (v^{k-1} + (v+1)^{k-1} + \cdots)
\end{equation}
(which agrees with (\ref{casek=1}) when $k=1$ and $v=0$), and this
is formally $\frac{1}{(k-1)!} \zeta (-k+1,v)$ (so that $\zeta(0,0)
= 1+\zeta(0)$, for example). Thus, formally equating the
coefficients of $x^l$ for $l \geq 0$ in (\ref{fromBerHur}) gives
an interpretation of (\ref{hurwitzformal}) and ``explains'' the
well--known relation between Bernoulli numbers and the Hurwitz
zeta function
\begin{equation}\label{hurzetaB}
\zeta(-k+1,v) = - \frac{B_{k}(v)}{k}
\end{equation}
for $k\geq1$. Note that setting $v=1$ in (\ref{hurdef}) gives the
Riemann zeta function:
\beq
    \zeta(s,1) = \zeta(s).
\eeq
Comparing (\ref{zetaB}) with (\ref{hurzetaB}) at $v=1$ then gives
the well-known relations
\beq
    B_k(1) = B_k \ \ (= B_k(0))
\eeq
for $k>1$ and
\beq\label{B11}
    B_1(1) = B_1 + 1 \ \ (= B_1(0) + 1).
\eeq
Note that as above, considering the pole in (\ref{fromBerHur})
gives $B_0(v) = 1$ as a polynomial, but this pole does not arise
in (\ref{unrig-Hurwitz}) -- (\ref{B11}).

The key point in these formal unrigorous arguments is the
interplay between the formal geometric series expansion (in powers
of $e^x$) and the expansion in powers of $x$.

Bloch \cite{Bl} found very interesting phenomena relating the values
$\zeta(-n)$, $n = 1, 3, 5, \dots$, of the zeta function at negative
odd integers to the commutators of certain formal differential
operators on the circle, and used (\ref{zeta-s}) in order to interpret
them.  We now sketch a part of his work and interpretation. Consider
the Lie algebra
\begin{equation}
{\frak {d}} = \mbox{\rm Der}\ \Bbb C [t,t^{-1}]
\end{equation}
of formal vector fields on the circle, with basis $\{t^n D | n \in
\Bbb Z\}$, where
\begin{equation}\label{D}
D = D_t = t \frac{d}{dt},
\end{equation}
and recall the Virasoro algebra ${\frak {v}}$, the central
extension
\begin{equation}\label{Vir}
0 \rightarrow \Bbb C c \rightarrow {\frak {v}} \rightarrow {\frak
{d}} \rightarrow 0,
\end{equation}
where ${\frak {v}}$ has basis $\{L(n) | n \in \Bbb Z \}$ together
with a central element $c$; the bracket relations are given by
\begin{equation}\label{Virbrackets}
[L(m),L(n)] = (m-n)L(m+n) + \frac{1}{12} (m^3 - m) \delta_{m+n,0}
c
\end{equation}
and $L(n)$ maps to $-t^n D$ in (\ref{Vir}).  This Lie algebra has
the following well-known realization: Start with the Heisenberg
Lie algebra with basis consisting of the symbols $h(n)$ for $n \in
\Bbb Z$, $n \neq 0$ and a central element $1$, with the bracket
relations
\begin{equation}\label{Heisbrackets}
[h(m),h(n)] = m \delta_{m+n,0} 1.
\end{equation}
For convenience we adjoin an additional central basis element
$h(0)$, so that the relations (\ref{Heisbrackets}) hold for all
$m,n \in \Bbb Z$.  This Lie algebra acts irreducibly on the
polynomial algebra
\begin{equation}\label{S}
S = \Bbb C [h(-1),h(-2),h(-3),\dots]
\end{equation}
as follows: For $n<0$, $h(n)$ acts by multiplication; for $n>0$,
$h(n)$ acts as $n \frac{\partial}{\partial h(-n)}$; $h(0)$ acts as
0; and $1$ acts as the identity operator.  Then ${\frak {v}}$ acts
on $S$ via:
\begin{equation}\label{Lnnormal}
L(n) \mapsto \frac{1}{2} \:\sum_{j \in \Bbb Z} h(j)h(n-j)\:
\end{equation}
for all $n \in \Bbb Z$.  Here the operators (\ref{Lnnormal}) are
in ``normal-ordered form,'' denoted by the colons, that is, the
$h(n)$ for $n>0$ act to the right of the $h(n)$ for $n<0$. It is
well known that the operators (\ref{Lnnormal}) indeed satisfy the
bracket relations (\ref{Virbrackets}).  (This exercise and the
related constructions are presented in \cite{FLM2}, for example,
where the standard generalization of this construction of ${\frak
{v}}$ using a Heisenberg algebra based on a finite-dimensional
space of operators $h(n)$ for each $n$ is also carried out.)

As is well known, removing the normal ordering in the definition
of $L(0)$ introduces an infinity which formally equals
$\frac{1}{2} \zeta(-1)$, since the unrigorous expression
\begin{equation}\label{Lbar(0)}
{\bar L}(0) = \frac{1}{2} \sum_{j \in \Bbb Z} h(-j)h(j)
\end{equation}
formally equals (by (\ref{Heisbrackets}))
\begin{equation}\label{Lbar(0)2}
L(0) + \frac{1}{2} (1 + 2 + 3 + \cdots),
\end{equation}
which itself formally equals, according to our considerations
above,
\begin{equation}\label{Lbar(0)3}
L(0) + \frac{1}{2} \zeta(-1) = L(0) - \frac{1}{24}.
\end{equation}
Rigorizing ${\bar L}(0)$ by defining it as
\begin{equation}\label{Lbar(0)rig}
{\bar L}(0) = L(0) + \frac{1}{2} \zeta(-1),
\end{equation}
we set
\begin{equation}\label{Lbar(n)}
{\bar L}(n) = L(n) \ \ \mbox{for} \ n \neq 0,
\end{equation}
to get a new basis of ${\frak {v}}$.  (We are identifying the
elements of ${\frak {v}}$ with operators on $S$.)  The brackets
become:
\begin{equation}\label{newVirbrackets}
[{\bar L}(m),{\bar L}(n)] = (m-n){\bar L}(m+n) + \frac{1}{12} m^3
\delta_{m+n,0};
\end{equation}
that is, $m^3 - m$ in (\ref{Virbrackets}) has become the pure
monomial $m^3$.

In \cite{Bl}, Bloch considered the larger Lie algebra of formal
differential operators, spanned by
\begin{equation}
\{ t^n D^m | n \in \Bbb Z, \ m \geq 0 \}
\end{equation}
or more precisely, we restrict to $m > 0$ and further, to the Lie
subalgebra ${\cal D}^+$, containing $\frak {d}$, spanned by the
differential operators of the form $D^r (t^nD) D^r$ for $r \geq 0,
\ n \in \Bbb Z$.  Then we can construct a central extension of
${\cal D}^+$ using generalizations of the operators
(\ref{Lnnormal}):
\begin{equation}\label{Lrnnormal}
L^{(r)}(n) = \:\frac{1}{2} \sum_{j \in \Bbb Z} j^r h(j) (n-j)^r
h(n-j)\:
\end{equation}
for $n \in \Bbb Z$.  These operators provide \cite{Bl} a central
extension of ${\cal D}^+$ such that
\begin{equation}
L^{(r)}(n) \mapsto (-1)^{r+1} D^r (t^nD) D^r.
\end{equation}

A central point of \cite{Bl} is that the formal removal of the
normal-ordering procedure in the definition (\ref{Lrnnormal}) of
$L^{(r)}(0)$ adds the infinity `` $\sum_{n>0} n^{2r+1}$'' $=
(-1)^r \frac{1}{2} \zeta(-2r-1)$ (generalizing
(\ref{Lbar(0)})--(\ref{Lbar(0)3})), and if we correspondingly
define
\begin{equation}\label{lbarr(0)}
{\bar L}^{(r)} (0) = L^{(r)} (0) + (-1)^r \frac{1}{2} \zeta(-2r-1)
\end{equation}
and ${\bar L}^{(r)} (n) = L^{(r)} (n)$ for $n \neq 0$
(generalizing (\ref{Lbar(0)rig}) and (\ref{Lbar(n)})), the
commutators simplify in a remarkable way: The complicated
polynomial in the scalar term of $[{\bar L}^{(r)} (m),{\bar
L}^{(s)} (-m)]$ reduces to a pure monomial in $m$, by analogy
with, and generalizing, the passage {}from $m^3 - m$ to $m^3$ in
(\ref{newVirbrackets}); see \cite{Bl} (and below) for the formulas
and further results.

We now recast these results in the more natural language of formal
calculus. This will make more apparent the generalization and
interpretation stemming from the theory of vertex operator
algebras.

Using a formal variable $x$, we form the generating functions
\begin{equation}\label{h(x)}
h(x) = \sum_{n \in \Bbb Z} h(n)x^{-n}
\end{equation}
and
\begin{equation}
L^{(r)}(x) = \sum_{n \in \Bbb Z} L^{(r)}(n)x^{-n},
\end{equation}
and using $D_x$ to denote the operator $x \frac{d}{dx}$, we
observe that
\begin{equation}\label{Lrx}
L^{(r)}(x) = {\frac{1}{2}}\nordbullet(D_x^r h(x))^2\nordbullet,
\end{equation}
where the colons, as always, denote normal ordering.

Then we introduce suitable generating functions over the number of
{\it derivatives}, and we use the formal multiplicative analogue
\begin{equation}\label{infinitdil}
e^{y D_x} f(x) = f(e^y x)
\end{equation}
of the formal Taylor theorem
\begin{equation}\label{Taylor}
e^{y \frac{d}{dx}} f(x) = f(x+y),
\end{equation}
where $f(x)$ is an arbitrary formal series of the form $\sum_n a_n
x^n$; $n$ is allowed to range over something very general, like
$\Bbb Z$ or even $\Bbb C$, say; and the $a_n$ lie in a fixed
vector space (cf. \cite{FLM2}, Proposition 8.3.1).  Although
$\nordbullet(D_x^r h(x))^2\nordbullet$ (recall (\ref{Lrx})) is
hard to put into a ``good'' generating function over $r$, we make
the problem easier by making it more general: Consider
independently many derivatives on each of the two factors $h(x)$
in $\nordbullet h(x)^2\nordbullet$, use two new independent formal
variables $y_1$ and $y_2$, and form the generating function
\begin{equation}\label{Ly1y2}
L^{(y_1,y_2)}(x) = {\frac{1}{2}} \nordbullet(e^{y_1 D_x}
h(x))(e^{y_2 D_x} h(x))\nordbullet = {\frac{1}{2}} \nordbullet
h(e^{y_1} x)h(e^{y_2} x)\nordbullet
\end{equation}
(where we use (\ref{infinitdil})), so that $L^{(r)}(x)$ is a
``diagonal piece'' of this generating function.  Using formal
vertex operator calculus techniques, we can calculate
\begin{equation}\label{bracketofquadratics}
[\nordbullet h(e^{y_1} x_1)h(e^{y_2} x_1)\nordbullet,\nordbullet
h(e^{y_3} x_2)h(e^{y_4} x_2)\nordbullet].
\end{equation}

Now, formally removing the normal ordering gives the formal
expression $h(e^{y_1} x)h(e^{y_2} x)$, which is not rigorous, as
we see by (for example) trying to compute the constant term in the
variables $y_1$ and $y_2$ in this expression; the failure of this
expression to be defined in fact corresponds exactly to the
occurrence of formal sums like $\sum_{n>0} n^r$ with $r>0$, as we
have been discussing.  However, we have
\begin{equation}\label{hx1hx2}
h(x_1)h(x_2) = \nordbullet h(x_1)h(x_2)\nordbullet + x_2
\frac{\partial}{\partial x_2} \frac{1}{1 - x_2 / x_1}
\end{equation}
and it follows that
\begin{equation}\label{hex1hex2}
h(e^{y_1} x_1)h(e^{y_2} x_2) = \nordbullet h(e^{y_1} x_1)h(e^{y_2}
x_2)\nordbullet + x_2 \frac{\partial}{\partial x_2} \frac{1}{1 -
e^{y_2} x_2 / e^{y_1} x_1};
\end{equation}
note that $x_2\frac{\partial}{\partial x_2}$ can be replaced by $-
\frac{\partial}{\partial y_1}$ in the last expression.  The
expression $\frac{1}{1 - e^{y_2} x_2 / e^{y_1} x_1}$ came {}from,
and is, a geometric series expansion (recall (\ref{hx1hx2})).

If we try to set $x_1 = x_2 \ (= x)$ in (\ref{hex1hex2}), the
result is unrigorous on the left--hand side, {\it but the result
has rigorous meaning on the right--hand side}, because the
normal-ordered product $\nordbullet h(e^{y_1} x)h(e^{y_2}
x)\nordbullet$ is certainly well defined, and the expression $-
\frac{\partial}{\partial y_1} \frac{1}{1 - e^{-y_1 + y_2}}$ can be
interpreted rigorously as in (\ref{fromBer}); more precisely, we
take $\frac{1}{1 - e^{-y_1 + y_2}}$ to mean the formal (Laurent)
series in $y_1$ and $y_2$ of the shape
\begin{equation}\label{defof1/1-e}
\frac{1}{1 - e^{-y_1 + y_2}} = (y_1 - y_2)^{-1}F(y_1,y_2),
\end{equation}
where $(y_1 - y_2)^{-1}$ is understood as the expansion $y_1^{-1}
\sum_{j\geq0} y_2^j y_1^{-j}$ in nonnegative powers of $y_2$ and
negative powers of $y_1$, and $F(y_1,y_2)$ is an (obvious) formal
power series in (nonnegative powers of) $y_1$ and $y_2$.  This
motivates us to define a new ``normal-ordering'' procedure
\begin{equation}\label{hexhex}
\pp h(e^{y_1} x)h(e^{y_2} x)\pp = \nordbullet h(e^{y_1}
x)h(e^{y_2} x)\nordbullet - \frac{\partial}{\partial y_1}
\frac{1}{1 - e^{-y_1 + y_2}},
\end{equation}
with the last part of the right--hand side being understood as we
just indicated. This gives us a natural ``explanation'' of the
zeta--function--modified operators defined in (\ref{lbarr(0)}): We
use (\ref{hexhex}) to define the following analogues of the
operators (\ref{Ly1y2}):
\begin{equation}\label{Lbary1y2}
{\bar L}^{(y_1,y_2)}(x) = {\frac{1}{2}} \pp h(e^{y_1} x)h(e^{y_2}
x)\pp,
\end{equation}
and the operator ${\bar L}^{(r)} (n)$ (with $r\geq0,\,n\in
\mathbb{Z}$) is exactly $(r!)^2$ times the coefficient of $y_1^r
y_2^r x_0^{-n}$ in (\ref{Lbary1y2}); the significant case is the
case $n = 0$.
\begin{remark}
{\it This ``rigorization'' of the undefined formal expression
$h(e^{y_1} x)h(e^{y_2} x)$ corresponds exactly to Euler's heuristic
interpretation of (\ref{zetaB}) discussed above} (as in
\cite{L4}). Note that the formal series on the right-hand side of
(\ref{hexhex}) involves terms with negative powers of $y_1$, and
therefore, so does the expression (\ref{Lbary1y2}). However, these
terms singular in $y_1$ are not involved in the definition of the
operators $\bar L^{(r)}(n)$, which, as we just mentioned, are related
to the nonnegative powers of $y_1$ in $\bar L^{(y_1,y_2)}(x)$. This
omission of the singularities in $y_1$ corresponds to the fact that
the Bernoulli number $B_0$ did not appear in our unrigorous
considerations above. In calculating commutators of normal-ordered
bilinear expressions (\ref{hexhex}), one obtains, as in Proposition
\ref{propLbarbracket} below, for instance, similar normal-ordered
bilinear expressions, in such a way that these singularities exactly
cancel out, as they must.  This cancellation process (see proof of
Proposition \ref{propLbarbracket}) is quite subtle. Actually,
Proposition \ref{propLbarbracket} deals with the general {\it twisted}
case, and as such involves Bernoulli {\it polynomials} and is related
to our discussion of the Hurwitz zeta function above (recall
(\ref{fromBerHur}) -- (\ref{hurzetaB})), not just to our discussion of
the Riemann zeta function.
\end{remark}

With the new normal ordering (\ref{hexhex}) replacing the old one,
remarkable cancellation occurs in the commutator
(\ref{bracketofquadratics}). This gives rise to the simple form of
the central term in the bracket relations involving the new basis
elements ${\bar L}^{(r)}(n)$, noted by Bloch.

The removal of the normal-ordering operation, interpreted
rigorously as above by means of the zeta regularization, has a
simple but considerable generalization in the context of the
theory of vertex operator algebras. It is recalled below (see
(\ref{usualnormord})) that the normal-ordering operation in vertex
operator algebra theory has the form
\beq
    \: Y(u,x_2+x_0) Y(v,x_2) \: = \mbox{regular part in $x_0$ of }
Y(Y(u,x_0)v,x_2).
\eeq
Here $u$ lies in a vector space, say $V$, that forms a vertex
operator algebra, and $Y(\cdot,x)$ is the vertex operator map (see
our short review in Section \ref{sec:VOA}). In fact, formal series
like (\ref{h(x)}) are examples of what we call {\it homogeneous}
vertex operators $X(u,x)$ (see (\ref{homo})), for which the
normal-ordering operation takes the form
\beq\label{intro-XX}
    \: X(u,e^y x_2) X(v,x_2) \: = \mbox{regular part in $x_0$ of }
    X(Y[u,y]v,x_2)
\eeq
where on both sides one replaces the formal variable $y$ by
the formal series $\log\lt(1+\frc{x_0}{x_2}\rt)$, and where the
map $Y[\cdot,y]$ (see \cite{Z1}, \cite{Z2}; see (\ref{Zhu}) below)
defines on $V$ a vertex operator algebra (in ``cylindrical
coordinates'') isomorphic to that defined by the map $Y(\cdot,x)$.
Now, the generalization of the zeta-regularized removal of the
normal-ordering, as in (\ref{hexhex}), to {\it arbitrary} vertex
operators takes the extremely simple form (see
(\ref{ppX-untwisted}))
\beq\label{intro-ppXXpp}
    \pp X(u,e^y x_2) X(v,x_2) \pp = X(Y[u,y]v,x_2),
\eeq
where this time, $y$ does not need to be replaced by any formal series
and no regular part needs to be extracted. This simple but extremely
general operation is the generalization to vertex operator algebras of
the heuristic argument of Euler justifying formal unrigorous relations
like (\ref{1+}) and (\ref{zeta-s}). Note that the expression
(\ref{intro-ppXXpp}) generically contains singularities in $y$, as we
have remarked above in our simple example, and as we have also remarked,
this phenomenon corresponds to the presence of the pole (in $x$) in
(\ref{fromBer}) and in (\ref{fromBerHur}).

In this paper we introduce and develop an extension to twisted vertex
operators of the operation (\ref{intro-ppXXpp}) (see (\ref{ppX})), an
extension that generalizes the heuristic argument leading to
(\ref{hurwitzformal}). Also, we establish a new relation in the theory
of vertex operator algebras and twisted modules, which
we call ``modified weak associativity.''  (See our announcement
\cite{DLMi}.) This new relation essentially expresses
(\ref{intro-ppXXpp}) as a formal limit applied to an ordinary product,
without any normal ordering, of homogeneous vertex operators
multiplied by an appropriate ``resolving factor'' (see Theorems
\ref{theomwa}, \ref{theomwat} and \ref{theomwaht} below).  This
relation gives a simple way of calculating brackets of
$\pp\cdot\pp$-normal-ordered bilinear expressions like those in
(\ref{intro-ppXXpp}) (see Theorem \ref{commiter}), which we use in
turn to elucidate the simplification of the central term observed by
Bloch (see Theorem \ref{main1}).

Now that we have put our results in context, let us outline a few
important and interesting consequences and potential applications.

Our general commutator formula (Theorem \ref{commiter}) can be
used to obtain nontrivial identities involving Bernoulli
polynomials and Bernoulli numbers. For instance, formula (3.12) in
the proof of the commutator formula Proposition 3.4 in \cite{Bl},
a formula relating the coefficients of powers of $m$ in
$\sum_{k=1}^{m-1} (m-k)^n k^n$ to Bernoulli numbers, becomes a
consequence of our commutator formula, after we plug in
appropriate vectors and equate coefficients of certain powers on
both sides of the formula. This classical fact is just a sample; we
anticipate interesting results on properties of values of
Bernoulli polynomials by comparing different twistings in
different realizations of Lie algebras of differential operators
on the circle.

It is interesting that our work also sheds some light on the
representation theory of Lie algebras of differential operators on the
circle. It is known that irreducible highest weight modules of Lie
algebras of differential operators on the circle are classified in
terms of values of certain Bernoulli polynomials (see for instance
\cite{KR}). Our Corollary \ref{kacradul} (a consequence of our
construction) illustrates how bosonic twisted constructions of
untwisted modules offer a new understanding of this phenomenon.

In addition, \cite{L4} and \cite{L5} were in retrospect a starting
point for a study of certain multi-point trace functions carried out
by the third author in \cite{M3}. These multi-point trace functions
have several interesting properties (such as quasi-modularity) and
have been computed in a special case by Bloch and Okounkov \cite{BO}. A
fascinating fact is that these functions also yield stationary
Gromov-Witten invariants of an elliptic curve (after Okounkov and
Pandharipande \cite{OP}).  Thus it is of great interest to understand
these multi-point trace functions in the general, twisted setting. Our
present work, together with \cite{DLM}, serves as a natural framework
for such a study.

\subsection{The Virasoro algebra and cylindrical coordinates}

In our work two notions will be of particular importance: the notion
of twisted module for a vertex operator algebra, and, in a less
predominant but still important way, the algebra isomorphism
(\cite{Z1}, \cite{Z2}) corresponding geometrically to a change to
cylindrical coordinates. The second notion appeared briefly (but
crucially) above in (\ref{intro-XX}) and (\ref{intro-ppXXpp}), through
the vertex operator map denoted $Y[\cdot,x]$. In order to motivate
this appearance, let us briefly explain how the change to cylindrical
coordinates gives a well-known simple geometrical interpretation of
the passage from the usual bracket relations (\ref{Virbrackets}) for
the Virasoro algebra to (\ref{newVirbrackets}), discussed above in
terms of the zeta regularization.

Consider the Virasoro algebra ${\frak {v}}$ (\ref{Vir}) with basis
$\{L(n) | n \in \Bbb Z,\; c \}$, $c$ central and with the bracket
relations (\ref{Virbrackets}). Consider
\beq \label{somefor}
    T(z)=\sum_{n \in \mathbb{Z}} L(n)
    z^{-n-2},
\eeq
where $z$ is a complex variable and the generators $L(n)$ act on a
certain module for the Virasoro algebra. Let us assume that the
formula (\ref{somefor}) makes sense for $z \neq 0$ (this can be
made rigorous if we use matrix coefficients \cite{FLM2}). We are
interested in how $T(z)$ (the ``stress--energy tensor'')
transforms under a change of coordinates $z \mapsto z'=f(z)$,
where $f(z)$ is a holomorphic function. That is, we would like to
compute $T'(z')$ as an expansion in $z$ with coefficients
expressed in terms of the generators $L(n)$, $n \in \mathbb{Z}$.
This can be achieved by using the well--known transformation
formula for $T(z)$ (cf. \cite{G}):
$$
    \left(\frac{dz'}{dz}\right)^2T'(z')=T(z)-\frac{c}{12} \{z',z
    \}
$$
with
$$
    \{z',z \}=\left(\frac{f'''(z)}{f'(z)}-
    \frac{3}{2}\left( \frac{f''(z)}{f'(z)} \right)^2 \right).
$$
Now, consider an infinite cylinder constructed by
identifying $z'$ and $z'+2 k \pi i$ for every $z' \in \mathbb{C}$ and
$k \in \mathbb{Z}$. Then the map
\[
   z \mapsto z'={\rm ln}(z)
\]
defines a holomorphic
map from the punctured plane $\mathbb{C} \ \backslash \ \{0\}$ onto
the cylinder. The transformation of the
stress--energy tensor under $z'={\rm ln}(z)$ is given by
\beq
\label{Tcylcoord}
    T'(z') = z^2T(z) - \frc{c}{24},
\eeq
and if we write $T'(z') = \sum_{n\in \Z} \b{L}(n) z^{-n}$ we have
\[
   \b{L}(n) = L(n) - \frc{c}{24}\delta_{n,0}.
\]
Notice how the modes $\b{L}(n)$ are defined: there is no shift of
$-2$ in the powers of $z$. This is closely related to the concept
of homogeneous vertex operator alluded to above, which will play
an important role in our twisted construction. The new modes
$\b{L}(n)$ do not satisfy the standard bracket relations
(\ref{Virbrackets}) of the Virasoro algebra; instead, the central
term in the commutator is a pure monomial:
\beq
    [\b{L}(m),\b{L}(n)]=(m-n)\b{L}(m+n)+\frac{m^3}{12} \delta_{m+n,0}c.
\eeq
This well-known fact is the same phenomenon as the one described above
in the context of Bloch's results; see (\ref{newVirbrackets}). This
cylindrical coordinate transformation brings important modular
invariance properties to objects constructed out of the new modes
$\b{L}(n)$; in particular, the graded dimension
$\Tr|_{M(1)}\lt(q^{\b{L}(0)}\rt)$, closely related to the partition
function in statistical systems, is $1/\eta(q)$ where $\eta(q)$ is
Dedekind's eta--function, which has important modular properties when
viewed as a function of the variable $\tau$ with $q=e^{2\pi i
\tau}$. The transformation to cylindrical coordinates will play an
important role in our considerations. {}From the vertex operator
algebra point of view, this transformation is an isomorphism between
two vertex operator algebras. This isomorphism is well understood
\cite{Z1}, \cite{Z2}, \cite{H1}, \cite{H2}, and in \cite{Z1} and
\cite{Z2} is used in the course of explaining modular invariance
phenomena.  In this paper we will obtain additional results that will
put the previous description of its effects on the Virasoro algebra
into the general framework of vertex operator algebra theory.

\subsection{The present work} \label{sectpresentwork}

As we have been saying, work of Bloch's \cite{Bl} revealed a
connection between values of the (analytically continued) Riemann zeta
function at the negative integers and a central extension
$\hat{\mathcal{D}}^+$ of a certain classical Lie algebra of
differential operators on the circle, a subalgebra of the central
extension of the Lie algebra of formal differential operators on the
circle of all nonnegative orders (denoted $W_{1+\infty}$ in the
physics literature).  In addition to the work \cite{L4}, \cite{L5}
discussed above (and in this paper), in \cite{M1}--\cite{M3} the
general theory of vertex operator algebras was used in order to extend
all these results to the context of central extensions of classical
Lie superalgebras of differential operators on the circle, in
connection with values of the Hurwitz zeta functions \cite{M2}.

This paper is a continuation of these works and also an extension in
several directions, giving proofs and further explanation of results
announced in \cite{L4}, \cite{L5} and, for the generalization to the
twisted setting, \cite{DLMi}. The present goal is, first, to develop
new concepts and identities in the general theory of vertex operator
algebras and their twisted modules (this is carried out in Sections
\ref{sec:VOA}, \ref{sec:homoVOA} and \ref{sec:normord}), in particular
to obtain a new general Jacobi identity for twisted operators (Theorem
\ref{theojacobiXt}) and a general commutator formula for related
iterates of such operators (Theorem \ref{commiter}). These new
identities will then be specialized to the Heisenberg vertex operator
algebra (cf. \cite{LL}) and used to obtain proofs of results announced
in \cite{L4}, \cite{L5} and, mainly, to obtain a twisted construction,
announced in \cite{DLMi}, of the algebra $\hat{\mathcal{D}}^+$ studied
in \cite{Bl} (cf. Proposition \ref{propLbarbracket} and Theorem
\ref{main1}), combining and extending methods {}from \cite{L4},
\cite{L5}, \cite{M1}--\cite{M3}, \cite{FLM2}, \cite{FLM3} and
\cite{DL2}.

In these earlier papers, we used vertex operator techniques to analyze
untwisted actions of the Lie algebra $\hat{\mathcal{D}}^+$ on a module
for a Heisenberg Lie algebra of a certain standard type, based on a
finite-dimensional vector space equipped with a nondegenerate
symmetric bilinear form. Now consider an arbitrary isometry $\nu$ of
period say $p$, that is, with $\nu^p={1}$. We prove that the
corresponding $\nu$--twisted modules carry an action of the Lie
algebra $\hat{\mathcal{D}}^+$ in terms of twisted vertex operators,
parametrized by certain quadratic vectors in the untwisted module. In
particular, we extend a result {}from \cite{FLM2}, \cite{FLM3},
\cite{DL2} where actions of the Virasoro algebra were constructed
using twisted vertex operators.  In addition, we explicitly compute
certain ``correction'' terms for the generators of the ``Cartan
subalgebra'' of $\hat{\mathcal{D}}^+$ that naturally appear in any
twisted construction. These correction terms are expressed in terms of
special values of certain Bernoulli polynomials. They can in principle
be generated, in the theory of vertex operator algebras, by the formal
operator $e^{\Delta_x}$ involved in the construction of a twisted
action for a certain type of vertex operator algebra.  We generate
these correction terms in an easier way, using a new {\em modified
weak associativity} (see Theorem \ref{theomwaht}) that is a
consequence of the twisted Jacobi identity.

In \cite{KR} Kac and Radul established a relationship between the Lie
algebra of differential operators on the circle and the Lie algebra
$\widehat{\frak{gl}}(\infty)$. We believe that their work can be
modified for purposes of classification of quasi--finite highest
weight representations of $\hat{\mathcal{D}}^+$ as well (for related
constructions, generalizations and a relationship with dual pairs see
\cite{AFOQ}, \cite{AFMO}, \cite{FKRW}, \cite{KWY}). Our new methods
and motivation for studying Lie algebras of differential operators,
based on vertex operator algebras, are of a different scope, so we do
not pursue their direction (however, see Corollary \ref{kacradul}).

Various research directions are suggested by the present work.
These include understanding multi--point correlation functions
\cite{M3} \cite{BO} (certain graded $q$--series) built {}from the
twisted vertex operators considered in this paper, as well as
investigating relationships between $\mathcal{W}$--algebras (in
the sense of \cite{FKRW}) and the present work. Also, the relation
between the present work and \cite{BDM} is interesting. We shall
investigate these directions in future publications.

\sect{The Lie algebra $\h{\D}^+$ and its untwisted construction}
\label{sec:untwistedconstruction}

Let $\mathcal{D}$ be the Lie algebra of formal differential
operators on $\mathbb{C}^\times$ spanned by $t^n D^r$, where
$D=t\,\frac{d}{dt}$ and $n \in \mathbb{Z}$, $r \in \N$ (the
nonnegative integers). This Lie algebra has an essentially unique
one-dimensional central extension $\hat{\mathcal{D}}=\mathbb{C}c
\oplus \mathcal{D}$ (denoted in the physics literature by
$\mathcal{W}_{1+\infty}$).

The representation theory of the highest weight modules of
$\hat{\mathcal{D}}$ was initiated in \cite{KR}, where, among other
things, the complete classification problem of the so-called {\em
quasi-finite} representations\footnote{These are representations
with finite-dimensional homogeneous subspaces.} was settled.  The
detailed study of the representation theory of certain subalgebras
of $\hat{\mathcal{D}}$ having properties related to those of
certain infinite--rank ``classical'' Lie algebras was initiated in
\cite{KWY} along the lines of \cite{KR}. In \cite{Bl} and
\cite{M2}, related Lie algebras (and superalgebras) are considered
{}from different viewpoints. We will follow these lines and
concentrate on the Lie subalgebra $\hat{\mathcal{D}}^+$ described
in \cite{Bl} and recalled below.

View the elements $t^n D^r \;(n \in \mathbb{Z},\; r \in \N)$ as
generators of the central extension $\hat{\mathcal{D}}$. They can
be taken to satisfy the following commutation relations (cf.
\cite{KR}):
\beqa \label{dcom}
    && [t^m f(D), t^n g(D)] \\ &&
    \qquad = t^{m+n}(f(D+n)g(D)-g(D+m)f(D))+\Psi(t^m f(D),t^ng(D))
    c,\no
\eeqa
where $f$ and $g$ are polynomials and $\Psi$ is the
$2$--cocycle (cf. \cite{KR}) determined by
$$
    \Psi(t^m f(D), t^n g(D))=-\Psi(t^n g(D), t^m f(D))=
    \delta_{m+n,0} \sum_{i=1}^m f(-i)g(m-i), \; m > 0.
$$
We consider the Lie subalgebra $\D^+$ of $\mathcal{D}$ generated
by the formal differential operators
\beq
\label{Lnrdef}
    L_{n}^{(r)}=(-1)^{r+1} D^r (t^n D) D^r ,
\eeq
where $n\in\Z,\;r \in\N$ \cite{Bl}. The subalgebra $\D^+$ has
an essentially unique central extension (cf. \cite{N}) and this
extension may be obtained by restriction of the $2$--cocycle
$\Psi$ to $\mathcal{D}^+$. Let $\hat{\mathcal{D}}^+=\mathbb{C}c
\oplus \mathcal{D}^+$ be the nontrivial central extension defined
via the slightly normalized $2$--cocycle $-\frac{1}{2}\Psi$, and
view the elements $L_{n}^{(r)}$ as elements of
$\hat{\mathcal{D}}^+$. This normalization gives, in particular,
the usual Virasoro algebra bracket relations
\begin{equation}\label{virasoro}
    [L_m^{(0)},L_n^{(0)}]=(m-n)L_{m+n}^{(0)}+\frac{m^3-m}{12}
    \delta_{m+n,0}\,c.
\end{equation}

In \cite{Bl} Bloch discovered that the Lie algebra
$\hat{\mathcal{D}}^+$ can be defined in terms of generators that
lead to a simplification of the central term in the Lie bracket
relations. Consider generators $\b{L}_n^{(r)} \;( n\in\Z,\;
r\in\N)$ for the algebra $\hat{\mathcal{D}}^+$, defined by
\beq\label{bLnrbloch}
    \bar{L}_n^{(r)}=L_n^{(r)}+
    \frac{(-1)^r}{2}\zeta(-1-2r)\delta_{n,0}c.
\eeq
Bloch found the following commutation relations for these
generators:
\beq \label{bl2coc}
    [\bar{L}_{m}^{(r)},\bar{L}_{n}^{(s)}]= \sum_{i={\rm min} (r,s)}^{r+s}
    a_i^{(r,s)}(m,n) \bar{L}_{m+n}^{(i)}+\frc{(r+s+1)!^2}{2(2(r+s)+3)!}
    m^{2(r+s)+3}  \delta_{m+n,0} \, c.
\eeq
The structure constants $a_i^{(r,s)}(m,n)$ are consequences of the
bracket relation (\ref{dcom}) and can be defined as follows:
Consider the following symmetric polynomial in two variables
$x_1,x_2$, with parameters $n\in\Z,\;r,s\in\N$:
$$
    f^{(r,s)}(n;x_1,x_2) = (-1)^{s+1} \lt( (x_2+n)^{r+s+1} x_1^r
    x_2^s + (x_1+n)^{r+s+1} x_1^s x_2^r \rt).
$$
In general, any symmetric polynomial in two variables $x_1,x_2$
can be written in a unique way as a polynomial in $x_1+x_2$ and
$x_1x_2$. The polynomial above is obviously symmetric, so can be written
in the following form:
$$
    f^{(r,s)}(n;x_1,x_2) = \sum_{i={\rm min}(r,s)}^{r+s}
    f_i^{(r,s)}(n;x_1+x_2) (x_1x_2)^i
$$
where $f_i^{(r,s)}(n,x)$ are homogeneous polynomials in $n$ and $x$
of total degree $2(r+s-i)+1$.
Then the coefficients $a_i^{(r,s)}(m,n)$ are defined by
$$
    a_i^{(r,s)}(m,n) = f_i^{(r,s)}(n;-m-n).
$$
Notice that, oddly enough, the central term in the commutator
(\ref{bl2coc}) is a pure monomial in $m$, in contrast to the
central term in (\ref{virasoro}) and in other bracket relations
that can be found {}from (\ref{dcom}). As was announced in
\cite{L4}, \cite{L5} and reviewed in the Introduction above, in
order to conceptualize this simplification (especially the
appearance of zeta-values) one can construct certain
infinite-dimensional projective representations of $\D^+$ using
vertex operators. But first, let us explain Bloch's construction
and heuristic conceptualization \cite{Bl}.

As in \cite{FLM3}, consider the (infinite-dimensional) Lie algebra
$\h{\a{h}}$, the affinization of an abelian Lie algebra $\a{h}$ of
dimension $d$ (over $\mathbb{C}$) with nondegenerate symmetric
bilinear form $\<\cdot,\cdot\>$:
\[
    \h{\a{h}} = \coprod_{n\in\Z} \a{h} \otimes t^n \oplus \C C,
\]
with the commutation relations
\begin{eqnarray}
    [\alpha\otimes t^m,\beta\otimes t^n] &=& \<\alpha,\beta\>
    m\delta_{m+n,0} C \com{\alpha,\beta \in\a{h},\; m,n\in\Z} \no\\ {}
    [C,\h{\a{h}}] &=& 0. \no
\end{eqnarray}
Set
\[
    \h{\a{h}}^+ = \coprod_{n>0} \a{h}\otimes t^n ,\qquad
    \h{\a{h}}^- = \coprod_{n<0} \a{h}\otimes t^n.
\]
The subalgebra
\[
    \h{\a{h}}^+ \oplus \h{\a{h}}^- \oplus \C C
\]
is a Heisenberg Lie algebra. Form the induced (level--one)
$\h{\a{h}}$--module
\[
    S = \mathcal{U}(\hat{\a{h}}) \otimes_{\mathcal{U}\lt(
    \hat{\a{h}}^+ \oplus \a{h} \oplus \C C\rt)} \C
        \simeq S(\hat{\a{h}}^-) \for{(linearly),}
\]
where $\hat{\a{h}}^+\oplus \a{h}$ acts trivially on $\C$ and $C$ acts
as 1; $\mathcal{U}(\cdot)$ denotes universal enveloping algebra and
$S(\cdot)$ denotes the symmetric algebra.  Then $S$ is irreducible
under the Heisenberg algebra $\h{\a{h}}^+ \oplus \h{\a{h}}^- \oplus \C
C$.  We will use the notation $\alpha(n)\;(\alpha\in\a{h},\, n\in\Z)$
for the action of $\alpha\otimes t^n \in \hat{\a{h}}$ on $S$. Then the
correspondence
\beq\label{Lnr}
    L_n^{(r)} \mapsto  \frc12 \sum_{q=1}^d \sum_{j\in\Z} j^r (n-j)^r
    \: \alpha_q(j) \alpha_q(n-j) \: \com{n\in\Z}\,,\;
    c\mapsto d,
\eeq
where $\{\alpha_q\}$ is an orthonormal basis of $\a{h}$, and
where $\:\cdot\:$ is the usual normal ordering, which brings
$\alpha(n)$ with $n>0$ to the right, gives a representation of
$\hat{\mathcal{D}}^+$. Let us denote the operator on the
right--hand side of (\ref{Lnr}) by $L^{(r)}(n)$. In particular,
the operators $L^{(0)}(m) \;(m\in\Z)$ give a well-known
representation of the Virasoro algebra with central charge
$c\mapsto d$,
$$
    [L^{(0)}(m),L^{(0)}(n)] = (m-n)L^{(0)}(m+n) + d \,\frc{m^3-m}{12}\,
    \delta_{m+n,0},
$$
and the construction (\ref{Lnr}) for these operators is the
standard realization of the Virasoro algebra on a module for a
Heisenberg Lie algebra (cf. \cite{FLM3}).

As we explained in the Introduction, the appearance of zeta-values in
(\ref{bLnrbloch}) can be conceptualized by noting (see \cite{Bl}) that
if we remove the normal ordering in (\ref{Lnr}) and use the relation
$[\alpha_q(m),\alpha_q(-m)]=m$ to rewrite $\alpha_q(m) \alpha_q(-m)$,
with $m \geq 0$, as $\alpha_q(-m) \alpha_q(m)+m$, then the resulting
expression contains an infinite formal divergent series of the form
$$1^{2r+1}+2^{2r+1}+3^{2r+1}+ \cdots .$$
The procedures discussed and analyzed in the Introduction lead to the
operators (\ref{bLnrbloch}), and these operators satisfy the bracket
relations (\ref{bl2coc}).  In the rest of this paper, we proceed to
supply the details and proofs of the results announced in \cite{L4},
\cite{L5} (for a different proof of certain of these results, see
\cite{M2}) and \cite{DLMi}.  We remark that our results in the theory
of twisted modules for vertex operator algebras are much more general
than their consequences in the representation theory of the algebra
$\hat{\mathcal{D}}^+$; they can be applied to other examples of vertex
operator algebras, and we expect them to have other interesting
consequences.

\sect{Twisted modules for vertex operator algebras and
commutativity and associativity relations} \label{sec:VOA}

In this section, we recall the definition of vertex operator
algebra, (untwisted) module and
twisted module. We derive important commutativity and associativity
relations, some of which are known relations for vertex operator
algebras and their untwisted modules, and some of which are new
identities even when specialized to untwisted modules.  For the basic
theory of vertex operator algebras and modules, we will use the
viewpoint of \cite{LL}.

In the theory of vertex operator algebras, formal calculus plays a
fundamental role.  Here we recall some basic elements of formal
calculus (cf. \cite{LL}).  Formal calculus is the calculus of formal
doubly--infinite series of formal variables, denoted below by $x$,
$y$, and by $x_1,\,x_2,\ldots$, $y_1,y_2,\ldots$.  The central object
of formal calculus is the formal delta--function
\[
    \delta(x) = \sum_{n\in\Z} x^n
\]
which has the property
\[
    \delta\lt(\frc{x_1}{x_2}\rt) f(x_1) = \delta\lt(\frc{x_1}{x_2}\rt)
    f(x_2)
\]
for any formal series $f(x_1)$. The formal delta--function enjoys
many other properties, two of which are:
\beq\label{2delta}
    x_2^{-1}\delta\lt(\frc{x_1-x_0}{x_2}\rt) = x_1^{-1}
    \delta\lt(\frc{x_2+x_0}{x_1}\rt)
\eeq
and
\beq\label{3delta}
    x_0^{-1}\delta\lt(\frc{x_1-x_2}{x_0}\rt) + x_0^{-1}
    \delta\lt(\frc{x_2-x_1}{-x_0}\rt) =
    x_2^{-1}\delta\lt(\frc{x_1-x_0}{x_2}\rt).
\eeq
In these equations, binomial expressions of the type
$(x_1-x_2)^n,\,n\in\Z$ appear. Their meaning as formal series in
$x_1$ and $x_2$, as well as the meaning of powers of more
complicated formal series, is summarized in the ``binomial
expansion convention'' -- the notational device according to which
binomial expressions are understood to be expanded in nonnegative
integral powers of the second variable. When more elements of
formal calculus are needed below, we shall recall them.

\ssect{Vertex operator algebras and untwisted modules}

We recall {}from \cite{FLM3} the definition of the notion of vertex
operator algebra, a variant of Borcherds' notion \cite{Bo} of vertex
algebra:
\begin{defi}\label{VOA}
A {\bf vertex operator algebra} $(V,Y,{\bf 1},\omega)$, or $V$ for
short, is a $\mathbb{Z}$--graded vector space
\[
    V=\coprod_{n \in \mathbb{Z}} V_{(n)};
    \ \mbox{\rm for}\ v\in V_{(n)},\;\mbox{\rm wt}\ v = n,
\]
such that
\beqa
     &&  V_{(n)} = 0 \;\; \mbox{ for }\; n \mbox{ sufficiently
    negative,} \nn &&  \mbox{\rm dim }V_{(n)}<\infty\;\;\mbox{ for }\;
    n \in {\mathbb Z} , \no
\eeqa
equipped with a linear map $Y(\cdot,x)$:
\begin{eqnarray}
    Y(\cdot,x)\,: \ V&\to&(\mbox{\rm End}\; V)[[x, x^{-1}]]\nonumber \\
    v&\mapsto& Y(v, x)={\displaystyle \sum_{n\in{\mathbb Z}}}v_{n} x^{-n-1}
    \,,\;\;v_{n}\in \mbox{\rm End} \;V,
\end{eqnarray}
where $Y(v,x)$ is called the {\em vertex operator} associated with
$v$, and two particular vectors, ${\bf 1},\,\omega\in V$, called
respectively the {\em vacuum vector} and the {\em conformal
vector}, with the following properties:\\
{\em truncation condition:} For every $v,w \in V$
\beq
    v_n w=0
\eeq
for $n\in\Bbb Z$ sufficiently large;\\
{\em vacuum property:}
\beq
    Y({\bf 1},x) = 1_V \for{($1_V$ is the identity on $V$);}
\eeq
{\em creation property:}
\beq
    Y(v,x){\bf 1} \in V[[x]] \for{and} \lim_{x\to0} Y(v,x){\bf 1} = v \;;
\eeq
{\em Virasoro algebra conditions:} Let
\begin{equation}
    L(n)=\omega _{n+1}\;\; \mbox{\rm for} \;n\in{\mathbb Z},
        \;\;{\rm i.e.},\;\;
        Y(\omega, x)=\sum_{n\in{\mathbb Z}}L(n)x^{-n-2} \;.
\end{equation}
Then
\[
    [L(m),L(n)]=
    (m-n)L(m+n)+c_V\frc{m^{3}-m}{12}\,\delta_{n+m,0}\,1_V
\]
for $m, n \in {\mathbb Z}$, where $c_V\in {\mathbb C}$ is the
central charge (also called ``rank'' of $V$),
\[
    L(0) v=({\rm wt}\ v) v
\]
for every homogeneous element $v$, and we have the {\em
$L(-1)$--derivative property:}
\beq\label{Lm1der}
    Y(L(-1)u,x)=\frac{d}{dx}Y(u,x) \;;
\eeq
{\em Jacobi identity:}
\beqa\label{jacobi}
    && x_0^{-1}\dpf{x_1-x_2}{x_0} Y(u,x_1) Y(v,x_2)
        - x_0^{-1}\dpf{x_2-x_1}{-x_0} Y(v,x_2) Y(u,x_1) \no\\
    &&  \qquad = x_2^{-1} \delta\lt(\frc{x_1-x_0}{x_2}\rt)
            Y(Y(u,x_0)v,x_2)\;.
\eeqa
\end{defi}

An important property of vertex operators is skew--symmetry, which
is an easy consequence of the Jacobi identity (cf. \cite{FHL}):
\beq \label{skewsymm}
    Y(u,x)v = e^{x L(-1)} Y(v,-x)u.
\eeq
Another easy consequence of the Jacobi identity is the
$L(-1)$--bracket formula:
\beq\label{Lm1bracket}
    [L(-1),Y(u,x)] = Y(L(-1)u,x).
\eeq

Fix a vertex operator algebra $(V,Y,{\bf 1},\omega)$, with central
charge $c_V$.
\begin{defi}\label{VOAmodule}
A (${\Bbb Q}$--graded) {\bf module} $W$ for the vertex operator
algebra $V$ (or $V$--{\bf module}) is a ${\Bbb Q}$--graded vector
space,
$$
    W=\coprod_{n\in {\Bbb Q}}W_{(n)}; \for{for} v\in
    W_{(n)},\;\mbox{\rm wt}\ v = n,
$$
such that \beqa &&  W_{(n)} = 0 \for{for} n \mbox{ sufficiently
negative,} \nn &&  \dim W_{(n)}<\infty\for{for} n \in {\Bbb Q} ,
\no \eeqa equipped with a linear map
\begin{eqnarray}
    Y_W(\cdot,x)\,: \ V&\to&(\mbox{\rm End}\; W)[[x, x^{-1}]]\nonumber \\
    v&\mapsto& Y_W(v, x)={\displaystyle \sum_{n\in{\Bbb Z}}}v_{n}^W x^{-n-1}
    \,,\;\; v_{n}^W\in \mbox{\rm End} \; W,
\label{vo}
\end{eqnarray}
where $Y_W(v,x)$ is still called the {\em vertex operator}
associated with $v$, such that the following conditions hold:\\
{\em truncation condition:} For every $v \in V$ and $w \in W$
\beq\label{truncation}
    v_n^W w=0
\eeq
for $n\in\Bbb Z$ sufficiently large;\\
{\em vacuum property:}
\begin{equation} \label{vacuum}
    Y_W({\bf 1}, x)=1_{W};
\end{equation}
{\em Virasoro algebra conditions:} Let
\[
    L_W(n)=\omega _{n+1}^W\;\; \mbox{\rm for} \;n\in{\mathbb Z},
    \;\;{\rm i.e.},\;\;
    Y_W(\omega,x)=\sum_{n \in \mathbb{Z}} L_W (n)x^{-n-2}.
\]
We have
\beqa
    [L_W(m),L_W(n)]&=&(m-n)L_W(m+n)+
        c_V\frac{m^3-m}{12}\,\delta_{m+n,0}\, 1_W,\nn
    L_W(0) v&=&({\rm wt}\ v) v \no
\eeqa
for every homogeneous element $v\in W$, and
\beq\label{Lm1derm}
        Y_W(L(-1)u,x)=\frac{d}{dx}Y_W(u,x) \;;
\eeq
{\em Jacobi identity:}
\beqa\label{jacobim}
    &&  x_0^{-1}\dpf{x_1-x_2}{x_0} Y_W(u,x_1) Y_W(v,x_2)
    - x_0^{-1}\dpf{x_2-x_1}{-x_0} Y_W(v,x_2) Y_W(u,x_1) \no\\
    &&  \qquad = x_2^{-1} \delta\lt(\frc{x_1-x_0}{x_2}\rt)
        Y_W(Y(u,x_0)v,x_2).
\eeqa
\end{defi}

We now recall the main commutativity and associativity properties of
vertex operators in the context of modules (\cite{FLM3}, \cite{FHL},
\cite{DL1}, \cite{Li1}; cf. \cite{LL}), and then we will derive new
identities somewhat analogous to these.  All these identities will be
generalized to twisted modules below. Note that taking the module to
be the vertex operator algebra $V$ itself, the relations below
specialize to commutativity and associativity properties in vertex
operator algebras.

{}From the Jacobi identity (\ref{jacobim}), one can derive the
weak commutativity and weak associativity relations, respectively:
\beqa
    (x_1-x_2)^{k(u,v)}Y_W(u,x_1) Y_W(v,x_2) &=& (x_1-x_2)^{k(u,v)}
    Y_W(v,x_2) Y_W(u,x_1) \no \\
    \label{wcom} \\
    (x_0+x_2)^{l(u,w)} Y_W(u,x_0+x_2) Y_W(v,x_2)w &=& (x_0+x_2)^{l(u,w)}
    Y_W(Y(u,x_0)v,x_2)w,
    \no \\ \label{wass}
\eeqa
where $u,\,v\,\in V$ and $w\in W$, valid for large enough integers
$k(u,v)$ and $l(u,w)$, their minimum value depending respectively on
$u,\,v$ and on $u,\,w$. For definiteness, we will pick the integers
$k(u,v)$ and $l(u,w)$ to be the smallest integers for which the
relations above are valid. These relations imply the main ``formal''
commutativity and associativity properties of vertex operators, which,
along with the fact that these properties are equivalent to the
Jacobi identity, can be formulated as follows (see \cite{LL}):

\begin{theo} \label{theostructY}
Let $W$ be a vector space (not assumed to be graded) equipped with
a linear map $Y_W(\cdot,x)$ (\ref{vo}) such that the truncation
condition (\ref{truncation}) and the Jacobi identity
(\ref{jacobim}) hold. Then for $u,\,v\,\in V$ and $w\,\in W$,
there exist $k(u,v) \in\N$ and $l(u,w) \in \N$ and a (nonunique)
element $F(u,v,w;x_0,x_1,x_2)$ of $W((x_0,x_1,x_2))$ such that
\beqa
    &&  x_0^{k(u,v)} F(u,v,w;x_0,x_1,x_2) \in W[[x_0]]((x_1,x_2)), \no\\
    \label{propFuvw} &&  x_1^{l(u,w)} F(u,v,w;x_0,x_1,x_2) \in
    W[[x_1]]((x_0,x_2))
\eeqa
and
\beqa
    Y_W(u,x_1) Y_W(v,x_2)w &=& F(u,v,w;x_1-x_2,x_1,x_2), \no\\
    Y_W(v,x_2) Y_W(u,x_1)w &=& F(u,v,w;-x_2+x_1,x_1,x_2), \no\\
    \label{structY}
    Y_W(Y(u,x_0)v,x_2)w &=& F(u,v,w;x_0,x_2+x_0,x_2)
\eeqa
(where we are using the binomial expansion convention).
Conversely, let $W$ be a vector space equipped with a linear map
$Y_W(\cdot,x)$ (\ref{vo}) such that the truncation condition
(\ref{truncation}) and the statement above hold, except that
$k(u,v)$ ($\in\N$) and $l(u,w)$ ($\in \N$) may depend on all three
of $u,\,v$ and $w$. Then the Jacobi identity (\ref{jacobim}) holds.
\end{theo}

It is important to note that since $k(u,v)$ can be (and typically is)
greater than 0, the formal series $F(u,v,w;x_1-x_2,x_1,x_2)$ and
$F(u,v,w;-x_2+x_1,x_1,x_2)$ are not in general equal. Along with
(\ref{propFuvw}), the first two equations of (\ref{structY}) represent
formal commutativity, while the first and last equations of
(\ref{structY}) represent formal associativity, as formulated in
\cite{LL} (see also \cite{FLM3} and \cite{FHL}). We will not prove
this theorem; instead we will prove its twisted generalization below.

{}From the equations in Theorem \ref{theostructY}, we can derive
a number of relations similar to weak commutativity and weak
associativity but involving formal limit procedures (the meaning
of such formal limit procedures is recalled below).
Even though only one of these will be of use
in the following sections, we state here for completeness of the
discussion the two relations that are not ``easy'' consequences of
weak commutativity and weak associativity.

The first relation can be expressed as follows:

\begin{theo}\label{theovwa}
With $W$ as in Theorem \ref{theostructY},
\beq\label{varwass}
    \lim_{x_0 \to -x_2+x_1} \lt( (x_0+x_2)^{l(u,w)}
    Y_W(Y(u,x_0)v,x_2)w\rt) =
    x_1^{l(u,w)} Y_W(v,x_2) Y_W(u,x_1)w
\eeq
for $u,\,v\,\in V$.
\end{theo}

The meaning of the formal limit
\begin{equation}\label{formallimit}
    \lim_{x_0 \to -x_2+x_1} \lt((x_0+x_2)^{l(u,w)}
    Y_W(Y(u,x_0)v,x_2)w\rt)
\end{equation}
is that one replaces each power of the formal variable $x_0$ in
the formal series $(x_0+x_2)^{l(u,w)} Y_W(Y(u,x_0)v,x_2)w$ by the
corresponding power of the formal series $-x_2+x_1$ (defined using
the binomial expansion convention). Notice again that the order of
$-x_2$ and $x_1$ is important in $-x_2+x_1$, according to the
binomial expansion convention.

{\em Proof of Theorem \ref{theovwa}:} Apply the limit
$\lim_{x_0 \to -x_2+x_1}$ to the expression
\[
    (x_0+x_2)^{l(u,w)} Y_W(Y(u,x_0)v,x_2)w
\]
written as in the right--hand side of the third equation of
(\ref{structY}). This limit is well defined; indeed, the only
possible problems are the negative powers of $x_2+x_0$ in
$F(u,v,w,x_0,x_2+x_0,x_2)$, but they are cancelled out by the
factor $(x_0+x_2)^{l(u,w)}$. The resulting expression is read off
the second relation of (\ref{structY}). \eproof

\begin{remark}\label{rematrivialrelation}
It is instructive to consider the following relation, deceptively
similar to (\ref{varwass}), but that is in fact an immediate
consequence of weak associativity (\ref{wass}):
\beq \label{varwass2}
    \lim_{x_0 \to x_1-x_2} \lt( (x_0+x_2)^{l(u,w)}
    Y_W(Y(u,x_0)v,x_2)w\rt) =
    x_1^{l(u,w)} Y_W(u,x_1) Y_W(v,x_2)w.
\eeq
More precisely, it can be obtained by noticing that the
replacement of $x_0$ by $x_1-x_2$ independently in each factor in
the expression as written on the left--hand side of (\ref{wass})
is well defined. We emphasize that, by contrast, the relation
(\ref{varwass}) {\em cannot} be obtained in such a manner. Indeed,
although the formal limit procedure $\lim_{x_0 \to -x_2+x_1}$ is
of course well defined on the series on both sides of
(\ref{wass}), one cannot replace $x_0$ by $-x_2+x_1$ either in the
factor $Y_W(u,x_0+x_2)Y_W(v,x_2)w$ on the left--hand side or in
the factor $Y_W(Y(u,x_0)v,x_2)w$ on the right--hand side of
(\ref{wass}).
\end{remark}

The second nontrivial relation, which we call {\em modified weak
associativity}, will be important when
generalized to twisted modules. It is stated as:

\begin{theo}\label{theomwa}
With $W$ as in Theorem \ref{theostructY},
\beq\label{modwass}
    \lim_{x_1 \to x_2+x_0} \lt( (x_1-x_2)^{k(u,v)} Y_W(u,x_1)
    Y_W(v,x_2) \rt) =
    x_0^{k(u,v)} Y_W(Y(u,x_0)v,x_2)
\eeq
for $u,\,v\,\in V$.
\end{theo}

\proof Apply the limit $\lim_{x_1 \to x_2+x_0}$ to
the expression
\[
    (x_1-x_2)^{k(u,v)} Y_W(u,x_1) Y_W(v,x_2)
\]
written as in the right--hand side of the first equation of
(\ref{structY}). This limit is well defined, since negative
powers of $x_1-x_2$ in $F(u,v,w,x_1-x_2,x_1,x_2)$ are cancelled
out by the factor $(x_1-x_2)^{k(u,v)}$. The resulting expression
is read off the third relation of (\ref{structY}). \eproof

\begin{remark}\label{remares}
Equation (\ref{modwass}) can be written in the following form:
\beq
    \lim_{x_1 \to x_2+x_0} \lt( \lt(\frc{x_1-x_2}{x_0}\rt)^{k(u,v)}
    Y_W(u,x_1) Y_W(v,x_2) \rt) =
    Y_W(Y(u,x_0)v,x_2).
\eeq
The factor $\lt(\frc{x_1-x_2}{x_0}\rt)^{k(u,v)}$ appearing in
front of the product of two vertex operators on the left--hand
side is crucial in giving a well--defined limit, but when the
limit is applied to this factor without the product of vertex
operators, the result is simply $1$.  We will call such a factor a
``resolving factor''. Its power will be apparent, in particular,
in the proof of the commutator formula (\ref{gencomform}); as will
be explained, it allows one to evaluate nontrivial limits of sums
of terms with cancelling ``singularities'' in a straightforward
fashion, evaluating the limit of each term independently.
\end{remark}

\ssect{Twisted modules for vertex operator algebras}

The notion of twisted module for a vertex operator algebra was
formalized in \cite{FFR} and \cite{D} (see also the geometric
formulation in \cite{FrS}; see also \cite{DLM}),
summarizing the basic properties of the
actions of twisted vertex operators discovered in \cite{FLM2},
\cite{FLM3} and \cite{L2}; the main nontrivial axiom in this notion is
the twisted Jacobi identity of \cite{FLM3} (and \cite{L2});
cf. \cite{FLM2}.

A critical ingredient in formal calculus needed in the theory of
twisted modules is the appearance of fractional powers of formal
variables, like $x^{1/p},\, p\in \Z_+$ (the positive integers).
For the purpose of formal calculus, the object $x^{1/p}$ is to be
treated as a new formal variable whose $p$--th power is $x$.  The
binomial expansion convention is applied as stated at the
beginning of Section \ref{sec:VOA} to binomials of the type
$(x_1+x_2)^{1/p}$. {}From a geometrical point of view, these rules
correspond to choosing a branch in the ``orbifold structure''
described (locally) by the twisted vertex operator algebra module.

We now fix a positive integer $p$ and a primitive $p$--th root of unity
\beq\label{primitiveroot}
    \om \in \mathbb{C}.
\eeq
We record here two important properties of the formal delta--function
involving fractional powers of formal variables:
\beq\label{deltaid1}
    \delta(x) = \frc1p\sum_{r=0}^{p-1}\delta(\omega^{r}_p
    x^{1/p})
\eeq
and
\beq\label{deltaid2}
    x_2^{-1} \delta\lt(\om^r\lt(\frc{x_1-x_0}{x_2}\rt)^{1/p}\rt) =
    x_1^{-1}
    \delta\lt(\om^{-r}\lt(\frc{x_2+x_0}{x_1}\rt)^{1/p}\rt).
\eeq

Recall the vertex operator algebra $(V, Y, {\bf 1}, \omega)$ with
central charge $c_V$ of the previous subsection.  Fix an automorphism
$\nu$ of period $p$ of the vertex operator algebra $V$, that is, a
linear automorphism of the vector space $V$ preserving $\omega$ and
${\bf 1}$ such that
\beq
    \nu Y(v,x)\nu^{-1} = Y(\nu v,x) \ \mbox{for}\ v\in V,
\eeq
and
\beq
    \nu^p=1_V.
\eeq

\begin{defi}\label{VOAmodulet}
A (${\Bbb Q}$-graded) $\nu$-{\bf twisted} $V$-{\bf module} $M$ is
a ${\Bbb Q}$-graded vector space,
$$
    M=\coprod_{n\in {\Bbb Q}}M_{(n)}; \ \mbox{\rm for}\ v\in
    M_{(n)},\;\mbox{\rm wt}\ v = n,
$$
such that
\beqa
    &&  M_{(n)} = 0 \;\; \mbox{ for }\; n \mbox{ sufficiently negative,} \nn
    &&  \mbox{\rm dim }M_{(n)}<\infty\;\;\mbox{ for }\; n \in {\Bbb Q} , \no
\eeqa
equipped with a linear map
\begin{eqnarray}
    Y_M(\cdot,x)\,: \ V&\to&(\mbox{\rm End}\; M)[[x^{1/p},
    x^{-1/p}]]\nonumber \\
    v&\mapsto& Y_M(v, x)={\displaystyle \sum_{n\in{\frc1p\Bbb
    Z}}}v_{n}^\nu x^{-n-1}
    \,,\;\; v_{n}^\nu\in \mbox{\rm End} \; M,
\label{tvo}
\end{eqnarray}
where $Y_M(v,x)$ is called the {\em twisted vertex operator}
associated with $v$, such that the following conditions hold:\\
{\em truncation condition:} For every $v \in V$ and $w \in M$
\beq\label{truncationt}
    v_n^\nu w=0
\eeq
for $n\in\frc1p\Bbb Z$ sufficiently large;\\
{\em vacuum property:}
\begin{equation} \label{vacuumt}
    Y_M({\bf 1}, x)=1_{M};
\end{equation}
{\em Virasoro algebra conditions:} Let
\[
    L_M(n)=\omega _{n+1}^\nu\;\; \mbox{\rm for} \;n\in{\mathbb Z},
    \;\;{\rm i.e.},\;\;
    Y_M(\omega,x)=\sum_{n \in \mathbb{Z}} L_M (n)x^{-n-2}.
\]
We have
\beqa
    [L_M(m),L_M(n)]&=&(m-n)L_M(m+n)+
        c_V\frac{m^3-m}{12}\,\delta_{m+n,0}\, 1_M,\nn
    L_M(0) v&=&({\rm wt}\ v) v \no
\eeqa
for every homogeneous element $v$, and
\beq\label{Lm1dermt}
    Y_M(L(-1)u,x)=\frac{d}{dx}Y_M(u,x)\;;
\eeq
{\em Jacobi identity:}
\beqa\label{jacobitm}
    && x_0^{-1}\dpf{x_1-x_2}{x_0} Y_M(u,x_1) Y_M(v,x_2)
    - x_0^{-1}\dpf{x_2-x_1}{-x_0} Y_M(v,x_2) Y_M(u,x_1) \no\\
    &&  \qquad = \frc1p x_2^{-1} \sum_{r=0}^{p-1}
    \delta\lt(\om^r\lt(\frc{x_1-x_0}{x_2}\rt)^{1/p}\rt)
    Y_M(Y(\nu^r u,x_0)v,x_2).
\eeqa
\end{defi}

Note that when restricted to the fixed--point subalgebra $\{v \in
V\,|\, \nu v = v\}$, a twisted module becomes a {\em true} module: the
twisted Jacobi identity (\ref{jacobitm}) reduces to the untwisted one
(\ref{jacobim}), by (\ref{deltaid1}). This will enable us to construct
natural representations of the algebra $\h{\D}^+$ on suitable twisted
modules.

We derive below various commutativity and associativity properties of
twisted vertex operators. In order to express some of these
properties, we need one more element of formal calculus: a certain
projection operator. Consider the operator $P_{[[x_0,x_0^{-1}]]}$
acting on the space $\mathbb{C}\{x_0\}$ of formal series with any
complex powers of $x_0$, which projects to the formal series with
integral powers of $x_0$:
\beq\label{projection}
    P_{[[x_0,x_0^{-1}]]} \;:\;
    \mathbb{C}\{x_0\} \to \mathbb{C}[[x_0,x_0^{-1}]] \;.
\eeq
We will extend the meaning of this notation in the obvious way to
projections acting on formal series with coefficients lying in vector
spaces other than $\C$, vector spaces which might themselves be spaces
of formal series in other formal variables. Notice that when this
projection operator acts on a formal series in $x_0$ with powers that
are in $\frc1p \Z$, for instance on $f(x_0) \in \C
[[x_0^{1/p},x_0^{-1/p}]]$, it can be described by an explicit formula:
\[
    P_{[[x_0,x_0^{-1}]]} f(x_0) = \frc1p \sum_{r=0}^{p-1}
    \lt( \lim_{x^{1/p}\to \om^r x_0^{1/p}} f(x) \rt).
\]
(See Remark \ref{remaformallimitfrac} below for the meaning of formal
limit procedures involving fractional powers of formal variables.)  We
will also extend this projection notation to different kinds of formal
series in obvious ways.  For instance,
\[
    P_{x_0^{q/p} [[x_0,x_0^{-1}]]} \;:\;
    \mathbb{C}\{x_0\} \to \mathbb{C}\,x_0^{q/p}[[x_0,x_0^{-1}]] \;.
\]
Again, of course, we will extend the meaning of this notation to
formal series with coefficients in vector spaces other than $\C$.

The twisted Jacobi identity (\ref{jacobitm}) implies twisted
versions of weak commutativity and weak associativity, proved
below ($u,\,v\,\in V,\;w\in M$):
\beqa \label{wcomt}
    &&      (x_2-x_1)^{k} Y_M(v,x_2)Y_M(u,x_1) =
    (x_2-x_1)^{k} Y_M(u,x_1) Y_M(v,x_2)
    \nn && \\
    \label{wasst}
    && P_{[[x_0,x_0^{-1}]]} \lt( (x_0+x_2)^{l}
    Y_M(u,x_0+x_2)Y_M(v,x_2)w \rt) \no\\
    &&      \qquad = (x_2+x_0)^{l} \frc1p
    \sum_{r=0}^{p-1}
    \omega^{-lrp}_p Y_M(Y(\nu^ru,x_0)v,x_2)w.
\eeqa
These relations are valid for all large enough $k \in \N$ and $l
\in \frc1p \N$, their minimum value depending respectively on
$u,\,v$ and on $u,\,w$. For definiteness, we will denote these minimum
values by $k(u,v)$ and $l(u,w)$, respectively (they depend also on the
module $M$; in particular, they differ from the integer numbers
$k(u,v)$ and $l(u,w)$ used in the previous subsection in
connection with the module $W$). As in the untwisted case, these
relations imply the main ``formal'' commutativity and
associativity properties of twisted vertex operators \cite{Li2},
which, along with with the fact that these properties are
equivalent to the Jacobi identity, can be formulated as follows:

\begin{theo} \label{theostructYt}
Let $M$ be a vector space (not assumed to be graded) equipped with
a linear map $Y_M(\cdot,x)$ (\ref{tvo}) such that the truncation
condition (\ref{truncationt}) and the Jacobi identity
(\ref{jacobitm}) hold. Then for $u,\,v\,\in V$ and $w\,\in M$,
there exist $k(u,v) \in\N$ and $l(u,w) \in \frc1p\N$ and a
(nonunique) element $F(u,v,w;x_0,x_1,x_2)$ of
$M((x_0,x_1^{1/p},x_2^{1/p}))$ such that
\beqa
    &&  x_0^{k(u,v)} F(u,v,w;x_0,x_1,x_2)
    \in M[[x_0]]((x_1^{1/p},x_2^{1/p})), \no\\
    \label{propFuvwt}
    &&  x_1^{l(u,w)} F(u,v,w;x_0,x_1,x_2) \in
    M[[x_1^{1/p}]]((x_0,x_2^{1/p}))
\eeqa
and
\beqa
    Y_M(u,x_1) Y_M(v,x_2)w &=& F(u,v,w;x_1-x_2,x_1,x_2), \no\\
    Y_M(v,x_2) Y_M(u,x_1)w &=& F(u,v,w;-x_2+x_1,x_1,x_2), \no\\
    \label{structYt}
    Y_M(Y(\nu^{-s}u,x_0)v,x_2)w &=&
    \lim_{x_1^{1/p} \to \om^s (x_2+x_0)^{1/p}}
     F(u,v,w;x_0,x_1,x_2)
\eeqa
for $s\in\Z$ (where we are using the binomial expansion
convention). Conversely, let $M$ be a vector space equipped with a
linear map $Y_M(\cdot,x)$ (\ref{tvo}) such that the truncation
condition (\ref{truncation}) and the statement above hold, except
that $k(u,v)$ ($\in\N$) and $l(u,w)$ ($\in \frc1p\N$) may depend
on all three of $u,\,v$ and $w$. Then the Jacobi identity
(\ref{jacobitm}) holds.
\end{theo}

\begin{remark}\label{remaformallimitfrac}
Formal limit procedures involving fractional powers of formal
variables like $x_1^{1/p}$ have the same meaning as in
(\ref{formallimit}), but with $x_1^{1/p}$ being treated as a
formal variable by itself. For instance, the formal limit
procedure
\[
    \displaystyle{\lim_{x_1^{1/p} \to \om^s
    (x_2+x_0)^{1/p}}F(u,v,w;x_0,x_1,x_2)}
\]
above means that one replaces each integral power of the formal
variable $x_1^{1/p}$ in the formal series $F(u,v,w;x_0,x_1,x_2)$
by the corresponding power of the formal series $\om^s
(x_2+x_0)^{1/p}$ (defined using the binomial expansion
convention).
\end{remark}

{\em Proof of Theorem \ref{theostructYt}:} By applying $\Res_{x_0}
x_0^{k(u,v)}$ to the twisted Jacobi identity (\ref{jacobitm}) for
$k(u,v) \in\N$ minimal such that
\beq
    x_0^{k(u,v)} Y_M(\nu^r u,x_0)v \in M[[x_0]] \for{for all} r\in\Z
\eeq
(which exists by the truncation condition), we find
\beq\label{wcomt1}
    (x_1-x_2)^{k(u,v)} [Y_M(u,x_1),Y_M(v,x_2)] = 0 \com{u,\,v\,\in V}.
\eeq
This relation is the expression of weak commutativity
(\ref{wcomt}), which is the same as for untwisted modules
(\ref{wcom}). It immediately implies, using the truncation
condition and applying the equation to a vector $w\in M$, that
there exists an element of $M((x_0,x_1^{1/p},x_2^{1/p}))$
depending on three vectors $u,\,v\,\in V$ and $w\,\in M$, which we
denote $F(u,v,w;x_0,x_1,x_2)$, satisfying the first condition of
(\ref{propFuvwt}) and such that the product of two twisted vertex
operators can be written
\beqa
    Y_M(u,x_1) Y_M(v,x_2)w &=& F(u,v,w;x_1-x_2,x_1,x_2) \no\\
\label{YYF}
    \quad Y_M(v,x_2) Y_M(u,x_1)w &=& F(u,v,w;-x_2+x_1,x_1,x_2).
\eeqa
Many series $F(u,v,w;x_0,x_1,x_2)$ can be chosen to satisfy these
relations. In such a series $F(u,v,w;x_0,x_1,x_2)$, consider the
coefficient $c_n(x_0,x_2) \in M((x_0,x_2^{1/p}))$ of the monomial
$x_1^n$ for $n\in\frc1p\Z$. Suppose that
$(x_0+x_2)^{-j}c_n(x_0,x_2) \in M((x_0,x_2^{1/p}))$ for some
positive integer $j$. Then one can define a new series
$\t{F}(u,v,w;x_0,x_1,x_2)$, which still satisfies the first
condition of (\ref{propFuvwt}) and equations (\ref{YYF}), where
the term $x_1^n c_n(x_0,x_2)$ is replaced by $x_1^{n+j}
(x_0+x_2)^{-j} c_n(x_0,x_2)$. In view of such transformations, it
is always possible to choose $F(u,v,w;x_0,x_1,x_2)$ in the form
\[
    F(u,v,w;x_0,x_1,x_2) = \sum_{n\in\frc1p\Z,\;n\ge n_0}
        x_1^n c_n(x_0,x_2) \for{for some} n_0\in\frc1p\Z,
\]
such that for each $n$, the maximal
integer $j_n$ giving $(x_0+x_2)^{-j_n}c_n(x_0,x_2) \in
M((x_0,x_2^{1/p}))$ is $j_n=0$. Let us choose such a formal series
$F(u,v,w;x_0,x_1,x_2)$.

Now we pick $l\in\Z$ and $q\in\N,\,0\le q\le p-1$ such that $l+q/p
= l(u,w) - 1 + 1/p$, where $l(u,w)\in\frc1p\N$ is the minimal rational number
satisfying
\beq\label{condluw}
    x_1^{l(u,w)} Y_M(u,x_1)w \in M[[x_1^{1/p}]]
\eeq
(which exists by the truncation condition). Using the formal
delta--function relation (\ref{deltaid2}) and applying
$\Res_{x_1}x_1^{l+q/p}$ to the twisted Jacobi identity
(\ref{jacobitm}), itself applied to a vector $w$, we find
\beqa
    &&  P_{[[x_0,x_0^{-1}]]}
    \lt( (x_0+x_2)^{l+q/p} Y_M(u,x_0+x_2)Y_M(v,x_2)w \rt)\no\\
    \label{wasst1}
    &&  \qquad = (x_2+x_0)^{l+q/p} \frc1p \sum_{r=0}^{p-1} \omega^{-qr}_p
    Y_M(Y(\nu^ru,x_0)v,x_2)w.
\eeqa
(Here $P_{[[x_0,x_0^{-1}]]}$ is the projection defined by
(\ref{projection}).) {}From (\ref{YYF}) and (\ref{wasst1}), we
obtain
\beqa \label{ppp1} &&
    \lim_{x_1 \to x_0+x_2} P_{[[x_1,x_1^{-1}]]}\,\lt( x_1^{l+q/p}\,
        F(u,v,w;x_0,x_1,x_2)  \rt)
    \\&&   \qquad = (x_2+x_0)^{l+q/p}
    \frc1p \sum_{r=0}^{p-1} \om^{-qr}Y_M(Y(\nu^r u,x_0) v,x_2)w \;.
    \no
\eeqa
The right--hand side of (\ref{ppp1}) contains only finitely many
negative powers of $x_0$. In view of the comment after
(\ref{YYF}), this implies that the expression inside the limit on
the left--hand side does not contain negative powers of $x_1$.
That is, the part of the series $F(u,v,w;x_0,x_1,x_2)$ for which
the powers of $x_1$ have a fractional part equal to $-q/p$, is of
the form $x_1^{-l-q/p}M[[x_1]]((x_0,x_2^{1/p}))$.

Now, the argument is left unchanged, and formula (\ref{ppp1})
stays valid, if we choose $l\in\Z$ and $q\in\N,\,0\le q\le p-1$
such that $l+q/p=l(u,w)-1+t/p$ for all $t\in \Z_+$ (recall that
$l(u,w)\in \frc1p\Z$ is the minimal rational number satisfying
(\ref{condluw})). In particular, repeating the argument for the
cases $1\le t\le p$ gives, {}from Equation (\ref{wasst1}), the
twisted version of weak associativity (\ref{wasst}), and gives,
{}from the discussion after (\ref{ppp1}), the second condition of
(\ref{propFuvwt}). Other values of $t\in\Z_+$ do not give new
information.

Then, for any $l+q/p=l(u,w)-1+t/p$ with $t\in \Z_+$, we can change
the limit to $\lim_{x_1 \to x_2+x_0}$ on the left--hand side of
(\ref{ppp1}), and we obtain:
\beqa
    &&  \lim_{x_1 \to x_2+x_0} \lt(
    P_{x_1^{-q/p} [[x_1,x_1^{-1}]]}
    F(u,v,w;x_0,x_1,x_2) \rt) \nn
    && \qquad = \frc1p \sum_{r=0}^{p-1}
    \om^{-qr}Y_M(Y(\nu^r u,x_0) v,x_2)w \;.
\eeqa
Choosing values of $t$ large enough such that $l\ge l(u,w)$, one can
see that the formula above is valid for all $q\in\Z$. Hence we can
apply the summation $\sum_{q=0}^{p-1} \om^{-qs}$ for any $s\in\Z$ to
this formula, obtaining
\beq
    \lim_{x_1^{1/p} \to
    \om^s (x_2+x_0)^{1/p}}F(u,v,w;x_0,x_1,x_2)=
    Y_M(Y(\nu^{-s}u,x_0)v,x_2)w \com{s\in\Z}.
\eeq
Next, we prove the converse. Assume (\ref{propFuvwt}),
(\ref{structYt}) and the truncation condition.  In fact, the
truncation condition is essential for the statements (\ref{structYt})
to be valid. Using, in the last term of (\ref{3delta}), the relation
(\ref{deltaid1}), and applying the resulting identity to
$F(u,v,w;x_0,x_1,x_2)$, we obtain (\ref{jacobitm}).  \eproof

\begin{remark}
Note that this proof illustrates the phenomenon, which arises
again and again throughout the theory of vertex operator algebras,
that formal calculus inherently involves just as much ``analysis''
as ``algebra'': in many relations there are integers that can be
left unspecified, except for their minimum values, and the proof
involves taking these integers ``large enough''. Recall that
essentially the same issues arose for example in the use of formal
calculus for the proof of the Jacobi identity for (twisted) vertex
operators in \cite{FLM3} (see Chapters 8 and 9). This is certainly
not surprising, since we are using the Jacobi identity (for all
twisting automorphisms) in order to prove, in a different
approach, properties of (twisted) vertex operators.
\end{remark}

Along with (\ref{propFuvwt}), the first two equations of
(\ref{structYt}) represent what we call {\em formal commutativity}
for twisted vertex operators, while the first and last equations
of (\ref{structYt}) represent {\em formal associativity} for
twisted vertex operators. When specialized to the untwisted case
$p=1$ ($\nu=1_V$), these two relations lead respectively to the
usual formal commutativity and formal associativity for vertex
operators, as described in (\ref{structY}).

As in the case of ordinary vertex operators, one can write other
relations involving formal limit procedures. Among them, two
cannot be directly obtained {}from weak commutativity and weak
associativity. One of these, the relation generalizing
(\ref{varwass}), is stated as follows:

\begin{theo}\label{theovwat}
With $M$ as in Theorem \ref{theostructYt},
\beqa &&
    \lim_{x_0 \to -x_2+x_1} \lt( (x_2+x_0)^{l} \frc1p
    \sum_{r=0}^{p-1}
    \omega^{-lrp}_p Y_M(Y(\nu^ru,x_0)v,x_2)w \rt) \no\\
    \label{varwasst}
    &&  \qquad =  P_{[[x_1,x_1^{-1}]]} \lt(
        x_1^{l} Y_M(u,x_2)Y_M(v,x_1)w \rt),
\eeqa
for all $l\in\frc1p\Z,\; l\ge l(u,w)$.
\end{theo}

\proof This is proved along the lines of the proof of Theorem
\ref{theovwa}, with some additions due to the fractional
powers. One uses the third equation of (\ref{structYt}) in order to
rewrite the left--hand side of (\ref{varwasst}) as
\[
    \lim_{x_0 \to -x_2+x_1} \lt( (x_2+x_0)^{l} \frc1p \sum_{r=0}^{p-1}
    \omega^{-lrp}_p \lim_{x_3^{1/p} \to
    \om^{-r}(x_2+x_0)^{1/p}}
    F(u,v,w;x_0,x_3,x_2) \rt).
\]
The sum over $r$ keeps only the terms in which $x_2+x_0$ is raised
to a power which has a fractional part equal to the negative of
the fractional part of $l$. Multiplying by $(x_2+x_0)^l$, for any
$l\in\frc1p\Z,\; l\ge l(u,w)$, brings
the remaining series to a series with finitely many negative
powers of $x_2$ (as well as $x_0$), to which it is possible to
apply the limit $\lim_{x_0 \to -x_2+x_1}$. This limit of course
brings only integer powers of $x_1$, and the right--hand side of
(\ref{varwasst}) can be obtained {}from the second equation of
(\ref{structYt}). \eproof

\begin{remark}
A relation similar to the last one, but that is a direct
consequence of weak associativity (\ref{wasst}), is
\beqa
    && \lim_{x_0 \to x_1-x_2} \lt( (x_2+x_0)^{l} \frc1p
        \sum_{r=0}^{p-1} \omega^{-lrp}_p Y_M(Y(\nu^ru,x_0)v,x_2)w
        \rt) \no\\
    \label{varwass2t}
    && \qquad= P_{[[x_1,x_1^{-1}]]} \lt( x_1^{l}
    Y_M(u,x_1)Y_M(v,x_2)w \rt)
\eeqa
for all $l\in\frc1p\Z,\; l\ge l(u,w)$. This generalizes
(\ref{varwass2}) (see the comments in Remark
\ref{rematrivialrelation}). It can be obtained by applying the
formal limit involved in the left--hand side to both sides of
(\ref{wasst}).
\end{remark}

The most important relation for our purposes, generalizing
(\ref{modwass}) and which we call {\em modified weak
associativity} for twisted vertex operators, is given by the
following theorem:

\begin{theo}\label{theomwat}
With $M$ as in Theorem \ref{theostructYt},
\beqa
    &&  \lim_{x_1^{1/p} \to \om^{s}(x_2+x_0)^{1/p}}
    \lt( (x_1-x_2)^{k(u,v)} Y_M(u,x_1) Y_M(v,x_2) \rt) \no\\
    \label{modwasst} && \qquad\qquad= x_0^{k(u,v)}Y_M(Y(\nu^{-s}u,x_0)v,x_2)
\eeqa
for $u,\,v\,\in V$ and $s\in \Z$.
\end{theo}

\proof The proof is a straightforward generalization of the proof
of Theorem \ref{theomwa}. \eproof

\begin{remark}
The specialization of Theorems \ref{theostructYt} and
\ref{theomwat} to the untwisted case $p=1$ and $M=V$ gives,
respectively, Theorems \ref{theostructY} and \ref{theomwa}.
\end{remark}

Finally, we derive a simple relation that specifies the structure
of the formal series $Y_M(u,x)$ (see \cite{DL2}).

\begin{theo}\label{theotransnu}
With $M$ as in Theorem \ref{theostructYt},
\beq \label{transnu}
    \lim_{x_1^{1/p} \to \om^s x^{1/p}} Y_M(\nu^{s}u,x_1) = Y_M(u,x)
\eeq
for $u\in V$ and $s\in \Z$.
\end{theo}

\proof In the Jacobi identity (\ref{jacobitm}), replace $u$ by
$\nu^s u$ and $x_1^{1/p}$ by $\om^s x^{1/p}$. The right--hand side
becomes
\[
    \frc1p x_2^{-1} \sum_{r=0}^{p-1}
    \delta\lt(\om^{r+s}\lt(\frc{x-x_0}{x_2}\rt)^{1/p}\rt)
        Y_M(Y(\nu^{r+s} u,x_0)v,x_2),
\]
which is independent of $s$, as is apparent if we make the
shift in the summation variable $r\mapsto r-s$. Hence the
left--hand side is also independent of $s$. Choosing $v={\bf 1}$
and using the vacuum property (\ref{vacuumt}), this gives
\beqa
    &&    \lt( x_0^{-1}\dpf{x-x_2}{x_0}
    - x_0^{-1}\dpf{x_2-x}{-x_0}
    \rt) \lim_{x_1^{1/p} \to \om^s x^{1/p}} Y_M(\nu^s
    u,x_1) \\
    &&\qquad\qquad\qquad    =
    \lt( x_0^{-1}\dpf{x-x_2}{x_0}
    - x_0^{-1}\dpf{x_2-x}{-x_0}
    \rt) Y_M(u,x) \no
\eeqa
which, upon using (\ref{3delta}) and taking $\Res_{x_2}$, gives
(\ref{transnu}). \eproof

{}From Theorem \ref{theotransnu}, we directly have the following
corollary:

\begin{corol}
With $M$ as in Theorem \ref{theostructYt},
\[
    Y_M(u,x) = \sum_{n\in\Z+q/p} u_nx^{-n-1} \for{for} u\in V,\ \nu u = \om^q
    u,\ q\in\Z.
\]
\end{corol}

\sect{Homogeneous twisted vertex operators and associated
relations} \label{sec:homoVOA}

\ssect{Homogeneous untwisted vertex operators}

Our construction of $\h{\D}^+$ {}from a particular vertex operator
algebra, and our explanation of its connection with the Riemann
zeta function and with Bernoulli polynomials (in the twisted
case), involve homogeneous vertex operators. It is natural to
consider homogeneous vertex operators acting on a $V$-module $W$
as defined in Definition \ref{VOAmodule}. The homogeneous vertex operator
associated to $v\in V$ will be denoted $X_W(v,x)$. It is defined
mainly by the property
\beq
    {\rm wt} \lt( \Res_{x} x^{-1} X_W(v,x) \rt) = 0.
\eeq
This can be achieved through:
\beq\label{homo}
    X_W(v,x) = Y_W(x^{L(0)}v,x)
\eeq
for $v \in V$ (see \cite{FLM3}). Homogeneous vertex operators
clearly inherit the vacuum property (\ref{vacuum}) of ordinary
vertex operators, and their main properties can be expressed
through the equivalent of the Jacobi identity (\ref{jacobim}), a
rewriting of Theorem 4.2 of \cite{L4} (see also \cite{L5},
\cite{M2}):

\begin{theo} \label{theojacobiX}
For a $V$-module $W$ as in Definition \ref{VOAmodule} and for
$u,\,v\,\in V$, we have
\beqa
    && x_0^{-1}\delta\lt(\frc{x_1-x_2}{x_0}\rt) X_W(u,x_1)X_W(v,x_2) -
    x_0^{-1}\delta\lt(\frc{x_2-x_1}{-x_0}\rt) X_W(v,x_2)X_W(u,x_1) \no\\
    \label{jacobiX}
    && \  = x_1^{-1} \delta\lt( e^{y} \frc{x_2}{x_1}\rt)
        X_W(Y[u,y] v,x_2)
\eeqa
where
$$
    y = \log\lt(1+\frc{x_0}{x_2}\rt)
$$
and
\beq\label{Zhu}
    Y[u,y] = Y(e^{yL(0)}u,e^y-1).
\eeq
\end{theo}

\proof Consider the Jacobi identity (\ref{jacobim}) with the
replacement $u\mapsto x_1^{L(0)}u$ and $v\mapsto x_2^{L(0)}v$. The
left--hand side of (\ref{jacobim}) is then directly the left--hand
side of (\ref{jacobiX}). After using the identity (\ref{2delta}),
the right--hand side becomes
\[
    x_1^{-1} \delta\lt( \frc{x_2+x_0}{x_1}
    \rt)Y_W(Y(x_1^{L(0)}u,x_0)x_2^{L(0)}v,x_2).
\]
Using the property
\beq\label{propL0Y}
    x_2^{-L(0)} Y(u,x_0) x_2^{L(0)} =
    Y\lt(x_2^{-L(0)}u,\frc{x_0}{x_2}\rt),
\eeq
the right--hand side is
\[
    x_1^{-1} \delta\lt( \frc{x_2+x_0}{x_1}
    \rt)X_W\lt(Y\lt(\lt(\frc{x_1}{x_2}\rt)^{L(0)}u,\frc{x_0}{x_2}\rt)v,x_2\rt).
\]
Making the replacement $x_1 \mapsto x_2+x_0$ inside the $Y$ vertex
operator, allowed {}from the properties of the formal
delta--function, one obtains the right--hand side of
(\ref{jacobiX}). \eproof

The expressions $\log(1+x)$ and $e^y$ in (\ref{jacobiX}), where
$x$ and $y$ are formal variables, are defined by their series
expansions in nonnegative powers of $x$ and $y$, respectively. The
use of the formal variable $y$ is natural here in particular
because of the appearance of the vertex operator $Y[u,y]$
(\ref{Zhu}) defined and studied in \cite{Z1}, \cite{Z2}. These
operators give a new vertex operator algebra isomorphic to $V$,
and the isomorphism corresponds geometrically to a
change--of--coordinates transformation expressed formally by
$y\mapsto e^y-1$ (\cite{Z1}, \cite{Z2}; see \cite{H1}, \cite{H2}
for the generalization to arbitrary coordinate changes).

\begin{remark}\label{remacylcoord}
The Jacobi identity (\ref{jacobiX}) thus suggests that there is a
close relationship between homogeneous vertex operators and the
change to ``cylindrical coordinates'' mentioned in our
Introduction. More precisely, see the right--hand side of the
Jacobi identity (\ref{jacobiX}) as a generating function in $x_2$
of endomorphisms of $W$ associated to some elements of $V$. These
elements are computed using the vertex operator algebra structure
given by the ``cylindrical coordinates'' vertex operator map
$Y[\cdot,y]$. Hence, the endomorphisms thus generated generalize
in some sense the modes ${\bar L}_n$ introduced after formula
(\ref{Tcylcoord}). Indeed, similar iterates of vertex operators
will be used below in Section \ref{sec:main} to generate a basis
for the algebra $\h\D^+$, along with its realizations on certain
twisted spaces, with properties similar to those for the basis
${\bar L}_n$.
\end{remark}

{}From the Jacobi identity (\ref{jacobiX}), one can obtain the
following commutator formula (\cite{L4}, \cite{L5}):

\begin{corol}
With $X_W(u,x)$ defined by (\ref{homo}) and $u,\,v\,\in V$, we have
\beq\label{comform}
    [X_W(u,x_1),X_W(v,x_2)] = \Res_y \delta\lt( e^y\frc{x_2}{x_1}\rt)
        X_W(Y[ u,y] v,x_2).
\eeq
\end{corol}

\proof Consider the following general fact concerning formal
series:
$$
    \Res_x h(x) = \Res_y \lt( h(F(y))\,\frc{d}{dy} F(y)\rt) \for{for}
    h(x)\in A((x)),\; F(y)\in yA[[y]]
$$
where $A$ is a commutative associative algebra (or more generally,
a module for it) and where the coefficient of $y^1$ in $F(y)$ is
invertible. Apply $\Res_{x_0}$ on the  Jacobi identity
(\ref{jacobiX}). On the left--hand side this gives the commutator,
and on the right--hand side use the general fact above with $F(y)
= x_2(e^y-1)$. The commutator formula follows. \eproof

An important property that will be used below is what we call the
{\em $L[-1]$--derivative property} for homogeneous vertex
operators:

\begin{theo}\label{theoLm1X}
With $X_W(u,x)$ defined by (\ref{homo}), $u\in V$ and $L[-1] =
L(-1) + L(0)$, we have
\beq\label{Lm1X}
    X_W(e^{L[-1]y}u,x) = X_W(u,e^{y}x).
\eeq
\end{theo}

\proof Using the $L(-1)$--derivative property (\ref{Lm1derm}) and
the identity $L(-1)x^{L(0)} = x^{L(0)-1} L(-1)$, which is a
consequence of the Virasoro algebra commutation relations, we have
in general
\[
    x\frc{d}{dx} X_W(u,x) = X_W((L(0)+L(-1))u,x),
\]
which proves Theorem \ref{theoLm1X}. \eproof

As is suggested by the notation $L[-1]$, the combination
$L(-1)+L(0)$ is the mode of the conformal vector representing the
Virasoro element $L_{-1}$ in the vertex operator algebra $V$
endowed with the vertex operator map defined by (\ref{Zhu})
\cite{Z1}, \cite{Z2}.

We will also use below some properties of Zhu's vertex operators
$Y[u,y]$ defined by (\ref{Zhu}): their skew--symmetry property,
the equivalent of (\ref{skewsymm}), and their $L[-1]$--derivative
and $L[-1]$--bracket properties, equivalent to, respectively,
(\ref{Lm1der}) and (\ref{Lm1bracket}). These properties are
consequences of the work of Zhu \cite{Z1}, \cite{Z2} (see also the
work of Huang \cite{H1}, \cite{H2}), that is, of the fact that the
transformation (\ref{Zhu}) is a vertex operator algebra
isomorphism. In fact, more simply, they are consequences of the
Jacobi identity for the operators $Y[u,y]$ (which was proven from
their definition and from the Jacobi identity for the operators
$Y(u,x)$ using formal variable techniques in \cite{L3}). However,
for completeness, we give here simple direct proofs. We have:

\begin{theo}\label{theoZhu}
With $Y[u,y]$ defined by (\ref{Zhu}), $u,\,v\,\in V$ and $L[-1] =
L(-1) + L(0)$, we have
\beq\label{skewsymmZhu}
    Y[u,y]v = e^{L[-1]y}Y[v,-y]u,
\eeq
\beq\label{Lm1derZhu}
    Y[L[-1]u,y] = \frc{d}{dy} Y[u,y]
\eeq
and
\beq\label{Lm1bracketZhu}
    [L[-1],Y[u,y]] = Y[L[-1]u,y].
\eeq
\end{theo}

\proof First, note that the left--hand side of (\ref{jacobiX}) is
invariant if we replace $(x_0,x_1,x_2)$ by $(-x_0,x_2,x_1)$ and
$(u,v)$ by $(v,u)$.
Require this invariance on the right--hand side of
(\ref{jacobiX}). Under this transformation the delta--function
$x_1^{-1}\delta\lt( e^{y} \frc{x_2}{x_1}\rt) =
x_1^{-1}\delta\lt(\frc{x_2+x_0}{x_1}\rt)$ is invariant, and using
this delta--function, $y$ transforms as $y\mapsto-y$ and $x_2$
transforms as $x_2 \mapsto e^yx_2$. Hence we have
\beq\label{XYskewsymm}
    X_W(Y[u,y] v,x_2) = X_W(Y[v,-y] u,e^yx_2).
\eeq
Then
\[
    X_W(Y[u,-y] v,e^yx_2) = X_W(e^{y (L(-1) + L(0))}Y[v,-y]u,x_2),
\]
and the injectivity of the vertex operator map $X_W(\cdot,x_2)$
along with (\ref{Lm1X}) gives (\ref{skewsymmZhu}).

Second, using the $L(-1)$--derivative property (\ref{Lm1der}) and
again the identity $L(-1)x^{L(0)} = x^{L(0)-1} L(-1)$, it is a
simple matter to obtain (\ref{Lm1derZhu}).

Third, {}from the Jacobi identity (\ref{jacobi}), one can obtain
the bracket relation
\[
[L(0),Y(u,x)] = Y((xL(-1) + L(0))u,x).
\]
{}From this and {}from the $L(-1)$--bracket relation
(\ref{Lm1bracket}), equation (\ref{Lm1bracketZhu}) follows.
\eproof

Homogeneous vertex operators satisfy other important properties,
similar to the formal commutativity and associativity properties
of vertex operators. We will not state them here, rather we will
state and prove their twisted generalization below, which can be
easily specialized to the untwisted case.

\ssect{Homogeneous twisted vertex operators}

{}From now on we fix a ${\Bbb Q}$-graded $\nu$-twisted $V$-module
$M$, as in Definition \ref{VOAmodulet}. We define homogeneous
twisted vertex operators by a simple twisted generalization of
(\ref{homo}):
\beq\label{thomo}
    X_M(u,x) = Y_M(x^{L(0)}u,x),
\eeq
as in \cite{FLM3}. We now state and prove the following twisted
generalizations of the results above: a Jacobi identity, a
commutator formula similar to (\ref{comform}), and formal
commutativity and associativity properties, including modified
weak associativity, for these homogeneous twisted vertex
operators.

Recalling the definition (\ref{Zhu}), we have the following
twisted generalization of the Jacobi identity (\ref{jacobiX}) for
homogeneous vertex operators:

\begin{theo}\label{theojacobiXt}
For a $\nu$-twisted $V$-module $M$ as in Definition
\ref{VOAmodulet}, for $u,\,v\,\in V$ and with $Y[u,y]$ defined by
(\ref{Zhu}), we have
\beqa
    && x_0^{-1}\delta\lt(\frc{x_1-x_2}{x_0}\rt) X_M(u,x_1)X_M(v,x_2) -
    x_0^{-1}\delta\lt(\frc{x_2-x_1}{-x_0}\rt) X_M(v,x_2)X_M(u,x_1) \no\\
    \label{jacobiXt}
    && \qquad = \frc1p x_1^{-1} \sum_{r=0}^{p-1}
    \delta\lt( \om^{-r} \lt(e^y \frc{x_2}{x_1}\rt)^{1/p}\rt)
        X_M(Y[\nu^r u,y] v,x_2)
\eeqa
where
\beq
    y = \log\lt(1+\frc{x_0}{x_2}\rt).
\eeq
\end{theo}

\proof Consider the twisted Jacobi identity (\ref{jacobitm}) with
the replacement $u\mapsto x_1^{L(0)}u$ and $v\mapsto x_2^{L(0)}v$.
The left--hand side of (\ref{jacobitm}) is then directly the
left--hand side of (\ref{jacobiXt}). After using the identity
(\ref{deltaid2}), the right--hand side becomes
\[
    \frc1p x_1^{-1} \sum_{r=0}^{p-1} \delta\lt( \om^{-r}
    \lt(\frc{x_2+x_0}{x_1}\rt)^{1/p}\rt)
        Y_M(Y(\nu^r x_1^{L(0)}u,x_0) x_2^{L(0)} v,x_2).
\]
Using the property (\ref{propL0Y}) and the fact that $\nu L(0)
\nu^{-1} = L(0)$, the right--hand side is
\[
    \frc1p x_1^{-1} \sum_{r=0}^{p-1} \delta\lt( \om^{-r}
        \lt(\frc{x_2+x_0}{x_1}\rt)^{1/p}\rt) X_M\lt(Y\lt(
        \lt(\frc{x_1}{x_2}\rt)^{L(0)} \nu^ru,\frc{x_0}{x_2}\rt)
        v,x_2\rt).
\]
Doing the replacement $x_1 \mapsto x_2+x_0$ inside the $Y$ vertex
operator, allowed {}from the properties of the formal
delta--function, one obtains the right--hand side of
(\ref{jacobiXt}). \eproof

This Jacobi identity leads to the following commutator formula:

\begin{corol}
With $X_M(u,x)$ defined by (\ref{thomo}) and $u,\,v\,\in V$, we have
\beq \label{comformt}
    [X_M(u,x_1),X_M(v,x_2)] = \Res_y \frc1p \sum_{r=0}^{p-1}
        \delta\lt( \om^{-r} \lt(e^y\frc{x_2}{x_1}\rt)^{1/p} \rt)
        X_M(Y[\nu^r u,y] v,x_2).
\eeq
\end{corol}

\proof The proof is a straightforward generalization of the proof
of (\ref{comform}). \eproof

The modified weak associativity relation (\ref{modwasst}) yields
{\em modified weak associativity} for homogeneous twisted vertex
operators:

\begin{theo}\label{theomwaht}
With $X_M(u,x)$ defined by (\ref{thomo}), $u,\,v\,\in V$ and
$s\in\Z$, we have
\beqa \label{modwassht} &&
    \lim_{x_1^{1/p} \to \om^s (e^y x_2)^{1/p}}
    \lt(\lt(\frc{x_1}{x_2}-1\rt)^{k(u,v)} X_M(u,x_1)X_M(v,x_2)\rt) \nn
    &&  \qquad\qquad = \lt(e^y-1\rt)^{k(u,v)}
        X_M(Y[\nu^{-s}u,y]v,x_2).
\eeqa
\end{theo}

\proof Apply $\lim_{x_0 \to (e^y-1)x_2}$ to both sides of the
modified weak associativity relation (\ref{modwasst}). Since
$(e^y-1) \in y \C[[y]]$, this limit is applicable to any series in
$x_0$ with finitely many negative powers of $x_0$, which is
obviously the case of the left--hand side and independently of
both factors $x_0^{k(u,v)}$ and $Y_M(Y(\nu^{-s}u,x_0)v,x_2)$ on
the right--hand side of (\ref{modwasst}). On the left--hand side,
use the formula
\[
    \lim_{x_0 \to (e^y-1)x_2} (x_2 + x_0)^{1/p} = (e^y x_2)^{1/p}.
\]
Then make the replacement $u\mapsto x_1^{L(0)}u$ and $v\mapsto
x_2^{L(0)}v$ (and keep unchanged the integer number $k(u,v)$), and
recall techniques used in the proof of Theorem \ref{theojacobiXt}
to obtain homogeneous vertex operators {}from ordinary vertex
operators. Equation (\ref{modwassht}) is obtained by multiplying
through the result by the factor $x_2^{-k(u,v)}$. \eproof

\begin{remark}\label{remaresh}
Remark \ref{remares} concerning resolving factors generalizes to
the homogeneous twisted case upon rewriting (\ref{modwassht}) in
the form
\beq\label{eqremaresh}
    \lim_{x_1^{1/p} \to \om^s (e^y x_2)^{1/p}}
    \lt(\lt(\frc{x_1/x_2-1}{e^y-1}\rt)^{k(u,v)} X_M(u,x_1)X_M(v,x_2)\rt)
    = X_M(Y[\nu^{-s}u,y]v,x_2).
\eeq
The factor $\lt(\frc{x_1/x_2-1}{e^y-1}\rt)^{k(u,v)}$ appearing in
front of the product of two vertex operators on the left--hand
side is crucial in having a well defined limit, but when the limit
is applied to this factor without the product of vertex operators,
the result is simply 1.
\end{remark}

\begin{remark}\label{remanormord}
Another important concept that becomes apparent {}from the
rewriting (\ref{eqremaresh}) of Equation (\ref{modwassht}) is the
concept of generalized normal ordering. In the next section,
normal orderings in vertex operator algebra will be recalled and
generalized. A normal ordering is essentially an operation that
regularizes, in some well--defined way, the product of two (or
more) vertex operators with the same formal variable (that is,
``at the same point''). The usual one, denoted $\:\cdot\:$, has
the effect, in the Heisenberg algebra, of putting modes
$\alpha(n)$ with positive $n$ to the right of modes with negative
$n$. This normal ordering also has a general definition in the
theory of vertex operator algebras. On the other hand, the
resolving factor on the left--hand side of Equation
(\ref{eqremaresh}) allows one, in particular, to make the
replacement of $x_1$ by $x_2$ in the formal series consisting of
this resolving factor multiplied by a product of two homogeneous
twisted vertex operators with formal variables $x_1$ and $x_2$.
This is in the spirit of a normal ordering operation, and such a
limit involving a resolving factor will be identified as a normal
ordering operation of a generalized type in the next section. In
particular, this is a reinterpretation of the normal ordering
denoted $\pp\cdot\pp$ in \cite{L4}, \cite{L5}.
\end{remark}

We also have the homogeneous counterpart of the relations
(\ref{structYt}):

\begin{theo} \label{theostructXt}
For $u,\,v\,\in V$ and $w\in M$, there exists $k(u,v) \in\N$ and
$l'(u,w) \in \frc1p\N$ and an element $G(u,v,w;x_0,x_1,x_2)$ of
$M((x_0,x_1^{1/p},x_2^{1/p}))$ such that
\beqa
    && x_0^{k(u,v)} G(u,v,w;x_0,x_1,x_2) \in
    M[[x_0]]((x_1^{1/p},x_2^{1/p})), \no\\
    \label{propGuvw} &&  x_1^{l'(u,w)} G(u,v,w;x_0,x_1,x_2) \in
    M[[x_1^{1/p}]]((x_0,x_2^{1/p})),
\eeqa
and
\beqa
    X_M(u,x_1) X_M(v,x_2)w &=& G(u,v,w;x_1-x_2,x_1,x_2), \no\\
    \label{structX}
    X_M(v,x_2) X_M(u,x_1)w &=& G(u,v,w;-x_2+x_1,x_1,x_2), \\
    X_M(Y[\nu^{-s}u,y]v,x_2)w &=&
    \lim_{ x_1^{1/p} \to \om^s (e^y x_2)^{1/p}}
    G(u,v,w;(e^{y}-1)x_2,x_1,x_2), \no
\eeqa
for $s\in\Z$. Here $k(u,v)$ can be taken to be the same as in
Theorem \ref{theostructYt}, and $l'(u,w)$ can be taken to be
$l(u,w) - {\rm wt}\; u$ if $u$ is homogeneous (that is, if $L(0) u
= ({\rm wt}\; u) u$).
\end{theo}

\proof Writing the first two equations of (\ref{structYt}) with
$u$ and $v$ replaced, respectively, by $x_1^{L(0)}u$ and
$x_2^{L(0)}v$, and using the definition (\ref{thomo}), gives the
first two equations of (\ref{structX}), with $G(u,v,w;x_0,x_1,x_2)
= F(x_1^{L(0)}u,x_2^{L(0)}v,w;x_0,x_1,x_2)$. This last
identification in particular gives both conditions in
(\ref{propGuvw}), along with the fact that we can take $k(u,v)$ as
being the same as in Theorem \ref{theostructYt} and $l'(u,w)$ as
being equal to $l(u,w) - {\rm wt}\; u$ if $u$ is homogeneous. The
homogeneous counterpart of the last equation of (\ref{structYt})
can be obtained by using modified weak associativity for
homogeneous twisted vertex operators (\ref{modwassht}). Multiply
through the first equation of (\ref{structX}) by the factor
$\lt(x_1/x_2-1\rt)^{k(u,v)}$, and apply to both sides of the
resulting equation the limit $\lim_{x_1^{1/p} \to \om^s(e^y
x_2)^{1/p}}$. On the right--hand side, the limit can be taken
directly, giving
\[
    \lt(e^y-1\rt)^{k(u,v)} \lim_{ x_1^{1/p} \to \om^s (e^y x_2)^{1/p}}
    G(u,v,w;(e^y-1)x_2,x_1,x_2).
\]
On the left--hand side, use (\ref{modwassht}), giving
\[
    \lt(e^y-1\rt)^{k(u,v)} X_M(Y[\nu^{-s}u,y]v,x_2).
\]
Every factor in both resulting expressions can be multiplied by
$\lt(e^y-1\rt)^{-k(u,v)} \in \C ((y))$ independently. One can then
cancel out the factors $\lt(e^y-1\rt)^{k(u,v)}$ and obtain the
last equation of (\ref{structX}). \eproof

Again, along with (\ref{propGuvw}), the first two equations of
(\ref{structX}) represent what we call {\em formal commutativity}
for homogeneous twisted vertex operators, while the first and last
equations of (\ref{structX}) represent {\em formal associativity}
for homogeneous twisted vertex operators.

Finally, we mention two straightforward consequences of results
established above: the $L[-1]$--derivative property of homogeneous
twisted vertex operators and the transformation properties of
twisted vertex operators under the automorphism $\nu$. First,
using the $L(-1)$--derivative property (\ref{Lm1dermt}), we have,
as in the untwisted case (\ref{Lm1X}):
\beq\label{Lm1Xt}
    X_M(e^{L[-1]y}u,x) = X_M(u,e^{y}x).
\eeq
In particular, with the skew--symmetry of Zhu's vertex operators
(\ref{skewsymmZhu}), this gives:
\beq\label{XYskewsymmt}
    X_M(Y[u,y]v,x) = X_M(Y[v,-y]u,e^{y}x).
\eeq
Second, since weights of elements of $V$ are integers,
Equation (\ref{transnu}) can directly be written in terms of
homogeneous twisted vertex operators:
\beq \label{transnuX}
    \lim_{x_1^{1/p} \to \om^s x^{1/p}} X_M(\nu^{s}u,x_1) =
    X_M(u,x).
\eeq

\begin{remark}
Taking $p=1$ in both Theorem \ref{theomwaht} and Theorem
\ref{theostructXt} gives, respectively, modified weak
associativity, and formal commutativity and formal associativity
for homogeneous untwisted vertex operators.
\end{remark}

\sect{Normal orderings in vertex operator algebras and
generalizations}\label{sec:normord}

This section is not essential for establishing our main results,
but it extends important concepts in the theory of vertex operator
algebras and provides interesting interpretations of some of our
results. The concept of normal ordering appears naturally in the
construction of vertex operator algebras. A new type of normal
ordering denoted $\pp\cdot\pp$ was introduced in \cite{L4},
\cite{L5}, which, as we mentioned in Remark \ref{remanormord}, can
be reinterpreted as a limit process. This reinterpretation gives a
direct link between a natural definition of this normal ordering
and general vertex operator algebra principles. As was stated in
\cite{L4}, \cite{L5}, this new normal ordering is the most natural
one to use in order to define operators $\b{L}^{(k)}_n$ generating
the algebra $\h\D^+$ and having the particular properties
mentioned in Section \ref{sec:untwistedconstruction}. We will see
in the next section that a proper generalization of this normal
ordering is also most naturally related to the definition of the
action of these particular generators on a twisted space. We now
generalize and unify various normal-ordering operations in the
general framework of vertex operator algebras.

\ssect{Untwisted vertex operators}

In order to make sense of the product of two vertex operators at
the same ``point'' (that is, with the same formal variable),
normal-ordering operations have to be introduced. We now recall
the standard definition of normal ordering in vertex operator
algebras (cf. \cite{FLM3}, \cite{LL}). Split the vertex operator
into two parts:
$$
    Y(u,x)=Y^+(u,x)+Y^-(u,x),
$$
where
$$
    Y^+(u,x)=\sum_{n <0} u_n x^{-n-1},
$$
and
$$
    Y^-(u,x)=\sum_{n \geq 0} u_n x^{-n-1}.
$$
The $\: \ \:$ normal ordering is defined by:
$$
    \: Y(u,x) Y(v,x) \: = Y^+(u,x)Y(v,x)+Y(v,x)Y^-(u,x),
$$
and more generally:
$$
    \: Y(u,x_1) Y(v,x_2) \: = Y^+(u,x_1)Y(v,x_2)+Y(v,x_2)Y^-(u,x_1).
$$
It is not hard to see directly {}from the Jacobi identity that
$$
    Y(u_{-1} v,x)=\: Y(u,x) Y(v,x) \:.
$$
{}From the formulas (cf. \cite{FHL})
\beqa
    Y(e^{x_0L(-1)}u,x_2) &=&e^{x_0 \frac{d}{dx_2}}Y(u,x_2),\nn
    e^{x_0L(-1)} Y(u,x_2)e^{-x_0L(-1)}&=& Y(u,x_2+x_0)\no
\eeqa
and
$$
    Y^+(u,x_0)=e^{x_0L(-1)}u_{-1}e^{-x_0L(-1)},
$$
it follows that
\beqa \label{usualnormord}
    && Y(Y^+(u,x_0)v,x_2)= \nn &&
    = Y(e^{x_0L(-1)}u_{-1}e^{-x_0L(-1)}v,x_2)=e^{x_0 \frac{d}{dx_2}}
    Y(u_{-1}e^{-x_0L(-1)}v,x_2)\nn &&= e^{x_0
    \frac{d}{dx_2}}\left(Y^+(u,x_2)Y(e^{-x_0L(-1)}
    v,x_2)+Y(e^{-x_0L(-1)} v,x_2)Y^-(u,x_2) \right) \nn &&=
    Y^+(u,x_2+x_0)Y(v,x_2)+Y(v,x_2)Y^-(u,x_2+x_0) \nn &&= \:
    Y(u,x_2+x_0)Y(v,x_2) \: \com{u,\,v\,\in V}.
\eeqa
In other words, the $\: \ \:$ normal ordering gives the regular
part in $x_0$ (the part with nonnegative powers of $x_0$) of the
iterate $Y(Y(u,x_0)v,x_2)$. It is convenient to have a
normal-ordering notion that is {\em equal} to $Y(Y(u,x_0)v,x_2)$:
define $\xx \cdot \xx$ by
\beq \label{newxx}
    \xx Y(u,x_2+x_0)Y(v,x_2)\xx=Y(Y(u,x_0)v,x_2).
\eeq
If we recall the modified weak associativity relation
(\ref{modwass}), this is equivalent to
\beq\label{xxYYxxlim}
    \xx Y(u,x_2+x_0)Y(v,x_2)\xx =\lim_{x_1 \to x_2+x_0}
    \lt(\lt(\frc{x_1-x_2}{x_0}\rt)^{k}Y(u,x_1)Y(v,x_2)\rt).
\eeq
Note that on the right--hand side, we have a limit applied to an
ordinary (that is, not normal-ordered) product of vertex
operators. Below, we will find similar expressions for homogeneous
twisted vertex operators, recovering and generalizing the
definition of the $\pp\cdot\pp$ normal ordering introduced in
\cite{L4}, \cite{L5}.

\ssect{Twisted vertex operators}

We want to generalize the normal-ordering operations $\: \ \:$ and
$\xx \ \xx$ to twisted vertex operators. As we already mentioned,
more general {\em untwisted} vertex operators can be obtained by
taking $\: \ \:$ normal-ordered product of several ``generating''
vertex operators. For example,
\[
    Y(u_{-1}v,x)=\: Y(u,x)Y(u,x) \:.
\]
On the other hand, for twisted vertex operators there is no simple
formula of this sort (cf. \cite{FLM3}, \cite{DL2}) and it is a
nontrivial matter, as we explained in the introduction, to
construct general twisted operators.

For our purposes, it is natural to define the normal-ordering
operation $\:\cdot \:$ in the following way, simply generalizing a
similar formula for untwisted vertex operators (and motivated by
modified weak associativity (\ref{modwasst})):
\beq \label{xxx}
    \lim_{x_1^{1/p} \to \om^s (x_2+x_0)^{1/p}}   \:
    Y_M(u,x_1)Y_M(v,x_2) \: = Y_M(Y^+(\nu^{-s}u,x_0)v,x_2),
\eeq
for $s\in\Z$ and $u,\,v\,\in V$. Also, generalizing the $\xx \ \xx$
normal ordering, we define:
\beqa \label{xxxx}
    &&\lim_{x_1^{1/p} \to \om^s (x_2+x_0)^{1/p}} \xx Y_M(u,x_1)
    Y_M(v,x_2)\xx  \nn && \qquad \qquad = \lim_{x_1^{1/p} \to
    \om^s (x_2+x_0)^{1/p}}
    \lt(\lt(\frc{x_1-x_2}{x_0}\rt)^{k(u,v)}
    Y_M(u,x_1)Y_M(v,x_2)\rt) \no\\
    &&\qquad \qquad = Y_M(Y(\nu^{-s}u,x_0)v,x_2)
\eeqa
for $s\in\Z$ and $u,\,v\,\in V$.

The definitions (\ref{xxx}), (\ref{xxxx}) can be generalized to
more general changes of variables in the following way: Let $f(y)
\in \mathbb{C}[[y]]$ be a formal power series such that
$f(y)=1+a_1y+a_2y^2+\cdots$. Then
\beqa \label{xxxxy}
    &&\lim_{x_1^{1/p} \to \om^s f(y)^{1/p} x_2^{1/p}} \xx
    Y_M(u,x_1) Y_M(v,x_2)\xx \nn &&\qquad \qquad = \lim_{x_1^{1/p} \to
    \om^s f(y)^{1/p} x_2^{1/p}}
    \lt(\lt(\frc{x_1/x_2-1}{f(y)-1}\rt)^{k(u,v)}
    Y_M(u,x_1)Y_M(v,x_2)\rt) \no\\
    &&\qquad \qquad = Y_M(Y(\nu^{-s}u,x_2(f(y)-1))v,x_2).
\eeqa
The last formal expression is well defined because
$(f(y)-1)^k,\,k\in\Z$ is a well--defined element inside
$\mathbb{C}((y))$. {}From this, one can take the regular part in
$y$ in order to define new normal-ordering operations $\:\cdot
\:_f$ for every $f(y)$ as specified:
\beqa \label{xxxy}
    && \lim_{x_1^{1/p} \to \om^s f(y)^{1/p} x_2^{1/p}} \:
    Y_M(u,x_1)Y_M(v,x_2) \:_f \nn && \qquad \qquad = {\rm Reg}_y \
    Y_M(Y(\nu^{-s}u,x_2(f(y)-1))v,x_2),
\eeqa
where ${\rm Reg}_y$ stands for the operator taking the regular
part of a formal series in $y$. Those normal-ordering operations
are well defined when one sets $y=0$, so that they truly
regularize the product of two vertex operators (or their
derivatives) at the same point, as does the usual $\: \cdot \:$.
The special case where $y=x_0/x_2$ and $f(y)=1+y$ is the case of
this usual normal ordering.

\ssect{Homogeneous untwisted vertex operators}

The homogeneous counterpart of the normal ordering defined by
(\ref{xxYYxxlim}) will be denoted $\pp\cdot\pp$, and can naturally
be defined by
\beqa
    \pp X_W(u,e^y x_2)X_W(v,x_2)\pp
    &=& \lim_{x_1 \to e^y
        x_2}\lt(\lt(\frc{x_1/x_2-1}{e^y-1}\rt)^{k(u,v)}
        X_W(u,x_1)X_W(v,x_2)\rt) \nn
\label{ppX-untwisted}
    &=& X_W(Y[u,y]v,x_2),
\eeqa
for $u,\,v\,\in V$. This corresponds essentially to the formulas
(2.25), (2.26) of
\cite{L5}. There is also a homogeneous counterpart for the usual
normal ordering $\:\cdot\:$ in a vertex operator algebra: one just
takes the regular part in $y$ of the previous expression. We
denote it $\aa\cdot\aa$, and define it by
\beq\label{aaXXaa}
    \aa X_W(u,e^y x_2)X_W(v,x_2) \aa = X_W(Y^+[\nu^{-s}u,y]v,x_2),
\eeq
where $Y^+[u,y]$ is the regular
part in $y$ of $Y[u,y]$.

\ssect{Homogeneous twisted vertex operators}

The previous construction can naturally be extended to the twisted
setting. From
(\ref{xxxx}), it is natural to define the $\pp\cdot\pp$
normal ordering by
\beqa \label{xxxxh}
    && \lim_{x_1^{1/p} \to \om^s (x_2+x_0)^{1/p}}\pp
        X_M(u,x_1)X_M(v,x_2)\pp  \nn
    && \qquad\qquad =  \lim_{x_1^{1/p}
        \to \om^s (x_2+x_0)^{1/p}}
        \lt(\lt(\frc{x_1-x_2}{x_0}\rt)^{k(u,v)} X_M(u,x_1)X_M(v,x_2)\rt)
        \nn
    && \qquad \qquad = (x_2+x_0)^{{\rm wt}(u)} x_2^{{\rm wt}(v)}
        Y_M(Y(\nu^{-s}u,x_0)v,x_2)
\eeqa
for $s\in\Z$ and $u,\,v$ homogeneous in $V$. As in the previous
section (see (\ref{xxxxy}) and (\ref{xxxy})), we can generalize
this to a construction using an arbitrary change of variable
parametrized by a formal power series $f(y) \in \mathbb{C}[[y]]$
such that $f(y)=1+a_1y+a_2y^2+\cdots$, that is:
\beqa \label{xxxxhy}
    && \lim_{x_1^{1/p} \to \om^s f(y)^{1/p} x_2^{1/p}} \pp
    X_M(u,x_1)X_M(v,x_2)\pp \nn &&\qquad\qquad = \lim_{x_1^{1/p} \to
    \om^s f(y)^{1/p} x_2^{1/p}}
    \lt(\lt(\frc{x_1/x_2-1}{f(y)-1}\rt)^{k(u,v)}
    X_M(u,x_1)X_M(v,x_2)\rt) \nn && \qquad\qquad = (x_2f(y))^{{\rm
    wt}(u)} x_2^{{\rm wt}(v)} Y_M(Y(\nu^{-s}u,x_2(f(y)-1))v,x_2).
\eeqa
Again, taking the regular part in $y$ gives a family of
normal-ordering operations, which we denote $\aa\cdot\aa_f$, that
regularize products of homogeneous vertex operators and their
derivatives at the same point:
\beqa \label{xxxhy}
    && \lim_{x_1^{1/p} \to \om^s f(y)^{1/p} x_2^{1/p}} \aa
    X_M(u,x_1)X_M(v,x_2) \aa_f \\ && \qquad\qquad = \mbox{Reg}_y \
    \lt( (x_2f(y))^{{\rm wt}(u)}(x_2)^{{\rm wt}(v)}
    Y_M(Y(\nu^{-s}u,x_2(f(y)-1))v,x_2) \rt). \no
\eeqa
Specializing to the case $f(y)=e^y$, we get the twisted
generalization of the normal ordering $\aa\cdot\aa$ introduced in
(\ref{aaXXaa}):
\beqa
    && \lim_{x_1^{1/p} \to \om^s (e^y x_2)^{1/p}}
    \aa X_M(u,x_1)X_M(v,x_2) \aa \nn
    && \qquad\qquad = \mbox{Reg}_y \ \lt((x_2e^y)^{{\rm wt}(u)}
    (x_2)^{{\rm wt}(v)} X_M(Y(\nu^{-s}u,x_2 (e^y-1))v,x_2) \rt) \nn &&
    \qquad\qquad = X_M(Y^+[\nu^{-s}u,y]v,x_2),
\eeqa
for $s\in\Z$ and $u,\,v\,\in V$ (not necessarily homogeneous), where
$Y^+[u,y]$ is
the regular part in $y$ of $Y[u,y]$. Concerning the $\pp \cdot
\pp$ normal ordering, {}from (\ref{modwassht}) and
(\ref{xxxxhy}) and taking again $f(y)=e^y$ it easily follows that
\beqa\label{ppX}
    &&\lim_{x_1^{1/p} \to {\omega}_p^s (e^y x_2)^{1/p}} \pp
    X_M(u,x_1)X_M(v,x_2)\pp \nn &&\qquad\qquad = \lim_{x_1^{1/p} \to
    {\omega}_p^s (e^y
    x_2)^{1/p}}\lt(\lt(\frc{x_1/x_2-1}{e^y-1}\rt)^{k(u,v)}
    X_M(u,x_1)X_M(v,x_2)\rt)\nn && \qquad\qquad =
    X_M(Y[\nu^{-s}u,y]v,x_2)
\eeqa
which is the normal-ordering operation that will be the most
directly related to our construction of the algebra $\h{\D}^+$.
This is the generalization to the twisted setting of the formulas
(2.25), (2.26) of \cite{L5}.

\sect{Commutator formula for iterates on twisted modules}
\label{sec:commform}

Modified weak associativity for homogeneous twisted vertex
operators turns out to be a very useful calculational tool. Formal
limit operations respect products, under suitable conditions, and
using this principle, one can compute, for instance, commutators
of certain iterates in a natural way.  For our applications, an
important commutator is
$$
    [X_M(Y[u_1,y_1]v_1,x_1),X_M(Y[u_2,y_2]v_2,x_2)],
$$
which we would like to express in terms of similar iterates. Using
modified weak associativity (\ref{modwassht}) and the commutator
formula (\ref{comformt}), we find a generalization of the main
commutator formula of \cite{L5} (it is also related to similar
commutator formulas in \cite{M1}--\cite{M3}):

\begin{theo} \label{commiter}
For $X_M(\cdot,x)$ defined by (\ref{thomo}) and $u_1,v_1,u_2,v_2
\in V$,
\beqa   \label{gencomform}
    &&  [X_M(Y[u_1,y_1]v_1,x_1),X_M(Y[u_2,y_2]v_2,x_2)] = \\
\no&& \Res_y \frc1p \sum_{r=0}^{p-1} \lt\{
        \delta\lt(\om^{-r}\lt(e^{y_1-y}\frc{x_1}{x_2}\rt)^{1/p}\rt)
        X_M(Y[u_2,y_2]Y[\nu^{-r} v_1,-y_1+y]Y[\nu^{-r}u_1,y]v_2,x_2)
    \rt.
\\\no&& \quad +
        \delta\lt(\om^{-r}\lt(e^{-y}\frc{x_1}{x_2}\rt)^{1/p}\rt)
        X_M(Y[u_2,y_2]Y[\nu^{-r}u_1,y_1+y] Y[\nu^{-r} v_1,y]v_2,x_2)
\\\no&& \quad +
        \delta\lt(\om^{-r}\lt(e^{-y_2+y_1-y}\frc{x_1}{x_2}\rt)^{1/p}\rt)
        X_M(Y[Y[\nu^{-r}v_1,-y_1]Y[u_2,-y]\nu^{-r}u_1,y_2+y]v_2,x_2)
\\\no&& \lt.    \quad +
        \delta\lt(\om^{-r}\lt(e^{-y_2-y}\frc{x_1}{x_2}\rt)^{1/p}\rt)
        X_M(Y[Y[\nu^{-r}u_1,y_1]Y[u_2,-y]\nu^{-r}v_1,y_2+y]v_2,x_2)
    \rt\}.
\eeqa
\end{theo}

\begin{remark}
In order to obtain a form which more directly specializes to the
commutator formula of \cite{L5}, one has to modify the third term
on the right--hand side of (\ref{gencomform}). Using skew--symmetry
(\ref{skewsymmZhu}) to change the sign of $y_1$ in the
vertex operator $Y[\nu^{-r}v_1,-y_1]$ and using the
$L[-1]$--derivative property (\ref{Lm1derZhu}), one finds
\beqa   \label{gencomformLep}
&&  [X_M(Y[u_1,y_1]v_1,x_1),X_M(Y[u_2,y_2]v_2,x_2)] = \\
\no&& \Res_y \frc1p \sum_{r=0}^{p-1} \lt\{
        \delta\lt(\om^{-r}\lt(e^{y_1-y}\frc{x_1}{x_2}\rt)^{1/p}\rt)
        X_M(Y[u_2,y_2]Y[\nu^{-r} v_1,-y_1+y]Y[\nu^{-r}u_1,y]v_2,x_2)
    \rt.
\\\no&& \quad +
        \delta\lt(\om^{-r}\lt(e^{-y}\frc{x_1}{x_2}\rt)^{1/p}\rt)
        X_M(Y[u_2,y_2]Y[\nu^{-r}u_1,y_1+y] Y[\nu^{-r} v_1,y]v_2,x_2)
\\\no&& \quad +
        \delta\lt(\om^{-r}\lt(e^{-y_2+y_1-y}\frc{x_1}{x_2}\rt)^{1/p}\rt)
        X_M(Y[Y[Y[u_2,-y]\nu^{-r}u_1,y_1]\nu^{-r}v_1,y_2-y_1+y]v_2,x_2)
\\\no&& \lt.    \quad +
        \delta\lt(\om^{-r}\lt(e^{-y_2-y}\frc{x_1}{x_2}\rt)^{1/p}\rt)
        X_M(Y[Y[\nu^{-r}u_1,y_1]Y[u_2,-y]\nu^{-r}v_1,y_2+y]v_2,x_2)
    \rt\}
\eeqa
which specializes to the formula in Theorem 3.4 of \cite{L5} when
$p=1$.
\end{remark}

{\em Proof of Theorem \ref{commiter}:} Using (\ref{modwassht})
multiplied through by the factor $\lt(e^y-1\rt)^{-k(u,v)}$, we
rewrite the commutator as a commutator of quadratics:
\beqa &&
    [X_M(Y[u_1,y_1]v_1,x_1),X_M(Y[u_2,y_2]v_2,x_2)] =
    \nonumber \\
&&  \qquad\qquad \lim_{x_3 \to e^{y_1}x_1 \atop x_4 \to
    e^{y_2}x_2} \biggl\{
        \lt(\frc{x_3/x_1-1}{e^{y_1}-1}\rt)^{k(u_1,v_1)}
        \lt(\frc{x_4/x_2-1}{e^{y_2}-1}\rt)^{k_2} \ \cdot \nonumber \\
&& \qquad\qquad\qquad\qquad
    [X_M(u_1,x_3)X_M(v_1,x_1)\;,\;X_M(u_2,x_4)X_M(v_2,x_2)]
    \biggr\} \nonumber
\eeqa
where we take $k_2\ge k(u_2,v_2)$, otherwise unspecified for now.
Then we use (\ref{comformt}) to write the commutator as a sum of
four terms:
\beqa\label{comder1}
&&  [X_M(Y[u_1,y_1]v_1,x_1),X_M(Y[u_2,y_2]v_2,x_2)] = \\
\no&&   \Res_y \frc1p \sum_{r=0}^{p-1} \lim_{x_3 \to e^{y_1}x_1
    \atop x_4 \to e^{y_2}x_2}
    \lt\{
        \lt(\frc{x_3/x_1-1}{e^{y_1}-1}\rt)^{k(u_1,v_1)}
        \lt(\frc{x_4/x_2-1}{e^{y_2}-1}\rt)^{k_2}
    \rt.    \ \cdot
\\\no&& \qquad\qquad\lt(
        \delta\lt(\om^{-r}\lt(e^y\frc{x_4}{x_1}\rt)^{1/p}\rt)
        X_M(u_1,x_3)X_M(Y[\nu^rv_1,y]u_2,x_4) X_M(v_2,x_2)
    \rt.
\\\no&& \qquad\qquad+\ {}
        \delta\lt(\om^{-r}\lt(e^y\frc{x_2}{x_1}\rt)^{1/p}\rt)
        X_M(u_1,x_3)X_M(u_2,x_4) X_M(Y[\nu^rv_1,y]v_2,x_2)
\\\no&& \qquad\qquad+\ {}
        \delta\lt(\om^{-r}\lt(e^y\frc{x_4}{x_3}\rt)^{1/p}\rt)
        X_M(Y[\nu^ru_1,y]u_2,x_4) X_M(v_2,x_2)X_M(v_1,x_1)
\\\no&& \lt.\lt.   \qquad\qquad + \ {}
        \delta\lt(\om^{-r}\lt(e^y\frc{x_2}{x_3}\rt)^{1/p}\rt)
        X_M(u_2,x_4) X_M(Y[\nu^ru_1,y]v_2,x_2) X_M(v_1,x_1)
    \rt) \rt\}.
\eeqa
By the truncation property, formal series of the form $Y[u,y]v$ have
finitely many negative powers of $y$ for any $u,\,v\,\in V$. Also, in the
expression inside the braces on the right--hand side above, other
factors involving $y$ only contain nonnegative powers of $y$. Hence,
we conclude that when evaluating the residue in $y$, only a finite
number of terms are to be added, each term of course being a
doubly--infinite series in the variables
$x_1^{1/p},\,x_2^{1/p},\,x_3^{1/p},\,x_4^{1/p}$ with coefficients in
$\End M$. Each such term can be expressed as a product of three
homogeneous twisted vertex operators. In fact, each term contains a
factor of the form $X_M(u,x_4)X_M(v,x_2)$ for some $u\in V$ and $v\in
V$. Since there is a finite number of these terms, it is possible to
choose a value of $k_2$ which is greater than or equal to $k(u,v)$ for
all of the vectors $u,\,v$ involved in these factors (as well, of
course, as greater than or equal to $k(u_2,v_2)$). Doing so, we can
evaluate the limit $\lim_{x_4\to e^{y_2}x_2}$ of each term independently
using (\ref{modwassht}). This justifies the name ``resolving factor''
for the factor $\lt(\frc{x_4/x_2-1}{e^{y_2}-1}\rt)^{k_2}$ on the
right--hand side of (\ref{comder1}), for $k_2$ large enough (see
Remark \ref{remaresh}). This gives:
\beqa\label{comder2}
&&  [X_M(Y[u_1,y_1]v_1,x_1),X_M(Y[u_2,y_2]v_2,x_2)] = \\
\no&&   \Res_y \frc1p \sum_{r=0}^{p-1} \lim_{x_3\to e^{y_1}x_1}
    \lt\{
        \lt(\frc{x_3/x_1-1}{e^{y_1}-1}\rt)^{k(u_1,v_1)}
    \rt.    \ \cdot
\\\no&& \qquad\qquad\lt(
        \delta\lt(\om^{-r}\lt(e^{y_2+y}\frc{x_2}{x_1}\rt)^{1/p}\rt)
        X_M(u_1,x_3)X_M(Y[Y[\nu^rv_1,y]u_2,y_2]v_2,x_2)
    \rt.
\\\no&& \qquad\qquad+\ {}
        \delta\lt(\om^{-r}\lt(e^y\frc{x_2}{x_1}\rt)^{1/p}\rt)
        X_M(u_1,x_3) X_M(Y[u_2,y_2] Y[\nu^rv_1,y]v_2,x_2)
\\\no&& \qquad\qquad+\ {}
        \delta\lt(\om^{-r}\lt(e^{y_2+y}\frc{x_2}{x_3}\rt)^{1/p}\rt)
        X_M(Y[Y[\nu^ru_1,y]u_2,y_2]v_2,x_2) X_M(v_1,x_1)
\\\no&& \lt.\lt.  \qquad\qquad  + \ {}
        \delta\lt(\om^{-r}\lt(e^y\frc{x_2}{x_3}\rt)^{1/p}\rt)
        X_M(Y[u_2,y_2]Y[\nu^ru_1,y]v_2,x_2) X_M(v_1,x_1)
    \rt)^{\ } \rt\}.
\eeqa

In order to evaluate the limit $\lim_{x_3 \to e^{y_1}x_1}$ of each
term independently, one could replace $k(u_1,v_1)$ by an integer
$k_1\ge k(u_1,v_1)$ large enough. To make the procedure more
transparent, however, we will use extra resolving factors. First,
we apply (\ref{XYskewsymmt}) to the first and third terms on the
right--hand side of the previous equation and use the
delta-function properties:
\beqa\label{comder3}
&&  [X_M(Y[u_1,y_1]v_1,x_1),X_M(Y[u_2,y_2]v_2,x_2)] = \\
\no&&   \Res_y \frc1p \sum_{r=0}^{p-1} \lim_{x_3 \to e^{y_1}x_1}
    \lt\{
        \lt(\frc{x_3/x_1-1}{e^{y_1}-1}\rt)^{k(u_1,v_1)}
    \rt.    \ \cdot
\\\no&&  \lt(
        \delta\lt(\om^{-r}\lt(e^{y_2+y}\frc{x_2}{x_1}\rt)^{1/p}\rt)
        X_M(u_1,x_3)\lim_{x_1^{1/p} \to
        x_1^{1/p}\om^r} X_M(Y[v_2,-y_2]Y[\nu^rv_1,y]u_2,e^{-y}x_1)
    \rt.
\\\no&& + \ {}
        \delta\lt(\om^{-r}\lt(e^y\frc{x_2}{x_1}\rt)^{1/p}\rt)
        X_M(u_1,x_3) \lim_{x_1^{1/p} \to
        x_1^{1/p}\om^r} X_M(Y[u_2,y_2] Y[\nu^rv_1,y]v_2,e^{-y}x_1)
\\\no&&  + \ {}
        \delta\lt(\om^{-r}\lt(e^{y_2+y}\frc{x_2}{x_3}\rt)^{1/p}\rt)
        \lim_{x_3^{1/p} \to
        x_3^{1/p}\om^r} X_M(Y[v_2,-y_2]Y[\nu^ru_1,y]u_2,e^{-y}x_3)
        X_M(v_1,x_1)
\\\no&& \lt.\lt.   + \ {}
        \delta\lt(\om^{-r}\lt(e^y\frc{x_2}{x_3}\rt)^{1/p}\rt)
        \lim_{x_3^{1/p} \to
        x_3^{1/p}\om^r} X_M(Y[u_2,y_2]Y[\nu^ru_1,y]v_2,e^{-y} x_3)
        X_M(v_1,x_1)
    \rt)^{\ } \rt\}.
\eeqa
Then, inside the braces, we insert the resolving factor
\beq
    \lt(\frc{e^y x_3/x_1-1}{e^{y_1+y}-1}\rt)^{k_3} \, \lt(\frc{e^{-y}
    x_3/x_1-1}{e^{y_1-y}-1}\rt)^{k_3}
\eeq
where $k_3\ge 0$ is an integer large enough, yet unspecified. Note
that this factor gives 1 under $\lim_{x_3\to e^{y_1}x_1}$, so it
doesn't change the result, but it allows us to calculate it
easily. Indeed, looking at any {\em fixed} power of $y_2$, we can
use an argument similar to what we explained after (\ref{comder1})
and choose $k_3$ large enough (its minimum value depending on the
power of $y_2$) in order to evaluate the limit of each term inside
the braces independently using (\ref{modwassht}). Since $k_3$
clearly does not appear in the result for any power of $y_2$, this
procedure can be apply to the whole series in $y_2$.

It is convenient to write the resulting expression for the commutator
in a form which does not involve the factor $e^{-y}$ in the argument
of the homogeneous twisted vertex operators. This form can be obtained
by using the $L[-1]$--derivative property for homogeneous twisted
vertex operators (\ref{Lm1Xt}). Using further the identity
$e^{L[-1]y_2}Y[u,y_1] e^{-L[-1]y_2} = Y[u,y_1+y_2]$ (which comes from
the $L[-1]$--bracket property and the $L[-1]$--derivative property
(\ref{Lm1derZhu}) and (\ref{Lm1bracketZhu}) of Zhu's vertex
operators), we obtain:
\beqa   \label{comder4}
&&  [X_M(Y[u_1,y_1]v_1,x_1),X_M(Y[u_2,y_2]v_2,x_2)] = \\
\no&& \Res_y \frc1p \sum_{r=0}^{p-1} \lt\{
        \delta\lt(\om^{r}\lt(e^{-y_2-y}\frc{x_1}{x_2}\rt)^{1/p}\rt)
        X_M(Y[u_1,y_1]Y[\nu^{-r} v_2,-y_2-y]Y[\nu^{-r}u_2,-y]v_1,x_1)
    \rt.
\\\no&& \qquad +
        \delta\lt(\om^{r}\lt(e^{-y}\frc{x_1}{x_2}\rt)^{1/p}\rt)
        X_M(Y[u_1,y_1]Y[\nu^{-r}u_2,y_2-y] Y[\nu^{-r} v_2,-y]v_1,x_1)
\\\no&& \qquad +
        \delta\lt(\om^{r}\lt(e^{y_1-y_2-y}\frc{x_1}{x_2}\rt)^{1/p}\rt)
        X_M(Y[Y[\nu^{-r}v_2,-y_2]Y[u_1,y]\nu^{-r}u_2,y_1-y]v_1,x_1)
\\\no&& \lt.    \qquad +
        \delta\lt(\om^{r}\lt(e^{y_1-y}\frc{x_1}{x_2}\rt)^{1/p}\rt)
        X_M(Y[Y[\nu^{-r}u_2,y_2]Y[u_1,y]\nu^{-r}v_2,y_1-y]v_1,x_1)
    \rt\}.
\eeqa
This is completely equivalent to the commutator formula
(\ref{gencomform}). Indeed, the commutator is unchanged if we make
the transformation $x_1\leftrightarrow x_2$, $y_1\leftrightarrow
y_2$, $u_1\leftrightarrow u_2$ and $v_1\leftrightarrow v_2$, and
change the overall sign. Doing this operation on the right--hand
side and absorbing the overall sign in a change of sign of the
variable $y$ on which we take the residue, we obtain
(\ref{gencomform}). \eproof

\begin{remark}
It is important to note that although the derivation above does
not make it apparent, there are many subtleties related to
cancellation of ``singularities'' in evaluating the commutator by
taking the limits in (\ref{comder1}).  Recall that the formula
(\ref{comder1}) is valid for any $k_2 \ge k(u_2,v_2)$, but that
the limit on $x_4$ was taken by assuming $k_2$ large enough; for a
given vertex operator algebra, this can be (and typically is)
larger than $k(u_2,v_2)$. The result (\ref{comder2}) is valid no
matter the value of $k_2\ge k(u_2,v_2)$. What happens is that for
$k_2=k(u_2,v_2)$, for instance, the limit cannot typically be
taken independently on each of the four terms inside the braces on
the right--hand side of (\ref{comder1}), because they contain
``singular terms'' of the type $(x_4-x_2)^{-n}$ for $n$ positive
and large enough. However, the nontrivial statement is that all
these singular terms cancel out inside the braces, since the limit
must exist. In a particular vertex operator algebra, one could
(and often does, see \cite{FLM3} for instance) evaluate this limit
explicitly by taking $k_2=k(u_2,v_2)$ and observing the
cancellation of the remaining singular terms. A similar phenomenon
happens when evaluating the limit on $x_3$ in (\ref{comder3}).
\end{remark}

\sect{Main results} \label{sec:main}

We now obtain a representation of $\hat{\mathcal{D}}^+$ on a
certain natural module for a twisted affine Lie algebra based on a
finite-dimensional abelian Lie algebra (essentially a twisted
Heisenberg Lie algebra), generalizing Bloch's representation on
the module $S \simeq S(\hat{\a{h}}^-)$ constructed in Section
\ref{sec:untwistedconstruction}. This is also a generalization of
the twisted Virasoro algebra construction (see \cite{FLM2},
\cite{FLM3}, \cite{DL2}).

\ssect{Description of the twisted module}

Let $\a{h}$ be a finite-dimensional abelian Lie algebra (over
$\C$) of dimension $d$ on which there is a nondegenerate symmetric
bilinear form $\<\cdot,\cdot\>$.  Let $\nu$ be an isometry of
$\a{h}$ of period $p>0$:
$$
    \<\nu\alpha,\nu\beta\> = \<\alpha,\beta\> ,\quad
    \nu^p\alpha = \alpha
$$
for all $\alpha,\beta\in\a{h}$. Consider the affine Lie algebra
$\hat{\a{h}}$ and its abelian subalgebra $\hat{\a{h}}^-$ recalled
in Section \ref{sec:untwistedconstruction}. The induced
(level--one) $\hat{\a{h}}$-module $S \simeq S(\hat{\a{h}}^-)$
(linearly) carries a natural structure of vertex operator algebra.
This structure is constructed as follows (cf. \cite{FLM3}). First,
one identifies the vacuum vector as the element 1 in $S$: ${\bf
1}=1$. Recalling the notation $\alpha(n) \; (\alpha\in\a{h},\,
n\in\Z )$ for the action of $\alpha\otimes t^n \in \hat{\a{h}}$ on
$S$, one constructs the following formal series acting on $S$:
$$
    \alpha(x)=\sum_{n \in \Z} \alpha(n)x^{-n-1}\com{\alpha\in\a{h}}.
$$
Then, the vertex operator map $Y(\cdot,x)$ is given by
\beqa
&&    Y(\alpha_1(-n_1) \cdots \alpha_j(-n_j){\bf 1},x) \no\\
    && \qquad = \: \frc1{(n_1-1)!} \lt( \frc{d}{dx} \rt)^{n_1-1}
\alpha_1(x) \cdots
    \frc1{(n_j-1)!} \lt( \frc{d}{dx} \rt)^{n_j-1} \alpha_j(x) \:
\eeqa
for $\alpha_k\in\a{h},\, n_k\in\Z_+,\; k=1,2,\ldots,j$, for all $j\in\N$,
where
$\:\cdot\:$ is the usual normal ordering, which brings $\alpha(n)$
with $n>0$ to the right. Choosing an orthonormal basis
$\{\alpha_q|q=1,\ldots,d\}$ of $\a{h}$, the conformal vector is
$\omega= \frc12 \sum_{q=1}^d \alpha_q(-1) \alpha_q(-1){\bf 1}$.
This implies in particular that the weight of $\alpha(-n){\bf 1}$
is $n$:
\[
    L(0) \alpha(-n){\bf 1} = n\alpha(-n){\bf 1} \com{\alpha\in\a{h},\;
    n\in\Z_+}.
\]
The isometry $\nu$ on $\a{h}$ lifts naturally to an automorphism
of the vertex operator algebra $S$, which we continue to call
$\nu$, of period $p$.

We now proceed as in \cite{L1}, \cite{FLM2}, \cite{FLM3} and
\cite{DL2} to construct a space $S[\nu]$ that carries a natural
structure of $\nu$--twisted module for the vertex operator algebra
$S$. In these papers, the twisted module structure was observed
assuming that $\nu$ preserves a rational lattice in $\a{h}$. Since
the properties of this twisted module will be essential in our
argument below, we make the same assumption here.

Recalling our primitive $p$--th root of unity $\om$, for $r \in
\Z$ set
$$
    \a{h}_{(r)} = \{\alpha\in\a{h} \;|\; \nu\alpha = \om^r\alpha\}
    \subset \a{h}.
$$
For $\alpha\in\a{h}$, denote by $\alpha_{(r)},\; r\in\Z$, its
projection on $\a{h}_{(r)}$. Define the $\nu$-twisted affine Lie
algebra $\hat{\a{h}}[\nu]$ associated with the abelian Lie algebra
$\a{h}$ by
\beq
    \hat{\a{h}}[\nu] = \coprod_{n\in\frc1p\Z} \a{h}_{(pn)} \otimes t^n
    \oplus \C C
\eeq
with
\beqa
    [\alpha\otimes t^m,\beta\otimes t^n] &=& \<\alpha,\beta\>
    m\delta_{m+n,0}\, C \com{\alpha\in\a{h}_{(pn)},\;
    \beta\in\a{h}_{(pm)},\; m,n\in\frc1p\Z} \no\\ {}
    [C,\hat{\a{h}}[\nu]] &=& 0 .
\eeqa
Set
\beq
    \hat{\a{h}}[\nu]^+ = \coprod_{n>0} \a{h}_{(pn)}\otimes t^n ,
    \qquad \hat{\a{h}}[\nu]^- = \coprod_{n<0} \a{h}_{(pn)} \otimes t^n.
\eeq
The subalgebra
\beq
    \hat{\a{h}}[\nu]^+ \oplus \hat{\a{h}}[\nu]^- \oplus \C C
\eeq
is a Heisenberg Lie algebra. Form the induced (level-one)
$\hat{\a{h}}[\nu]$-module
\beq
    S[\nu] = \mathcal{U}(\hat{\a{h}}[\nu])
    \otimes_{\mathcal{U}\lt( \hat{\a{h}}[\nu]^+
        \oplus \a{h}_{(0)} \oplus \C C\rt)} \C
        \simeq S(\hat{\a{h}}[\nu]^-) \for{(linearly),}
\eeq
where $\hat{\a{h}}[\nu]^+\oplus \a{h}_{(0)}$ acts trivially on
$\C$ and $C$ acts as 1; $\mathcal{U}(\cdot)$ denotes universal
enveloping algebra.  Then $S[\nu]$ is irreducible under the
Heisenberg algebra $\hat{\a{h}}[\nu]^+ \oplus \hat{\a{h}}[\nu]^-
\oplus \C C$. We will use the notation
$\alpha^\nu(n)\;(\alpha\in\a{h}_{(pn)},\, n\in\frc1p\Z)$ for the
action of $\alpha\otimes t^n \in \hat{\a{h}}[\nu]$ on $S[\nu]$.

\begin{remark}
The special case where $p=1$ ($\nu=1_{\a{h}}$) corresponds to the
$\hat{\a{h}}$-module $S$ discussed in Section
\ref{sec:untwistedconstruction}.
\end{remark}

As we mentioned above, the $\hat{\a{h}}[\nu]$-module $S[\nu]$ is
naturally a $\nu$--twisted module for the vertex operator algebra
$S$. One first constructs the following formal series acting on
$S[\nu]$:
\beq\label{apar}
    \alpha^\nu(x) = \sum_{n\in\frc1p\Z} \alpha^\nu(n) x^{-n-1},
\eeq
as well as the formal series $W(v,x)$ for all $v\in S$:
\beqa
&&    W(\alpha_1(-n_1) \cdots \alpha_j(-n_j) {\bf 1},x) \no\\
&& \qquad    = \: \frc1{(n_1-1)!} \lt( \frc{d}{dx} \rt)^{n_1-1}
        \alpha_1^\nu(x) \cdots
        \frc1{(n_j-1)!} \lt( \frc{d}{dx} \rt)^{n_j-1} \alpha_j^\nu(x) \:
\eeqa
where $\alpha_k\in\a{h},\, n_k\in\Z_+,\; k=1,2,\ldots,j$, for all
$j\in\N$. The
twisted vertex operator map $Y_{S[\nu]}(\cdot,x)$ acting on
$S[\nu]$ is then given by
\beq\label{tvertex}
    Y_{S[\nu]}(v,x) = W(e^{\Delta_x} v,x) \com{v\in S}
\eeq
where $\Delta_x$ is a certain formal operator involving the formal
variable $x$ \cite{FLM2}, \cite{FLM3}, \cite{DL2}. This operator
is trivial on $\alpha(-n){\bf 1}\in S\;(n\in\Z_+)$, so that one
has in particular
\beq\label{tvertexalphax}
    Y_{S[\nu]}(\alpha(-n){\bf 1},x) = \frc1{(n-1)!}
        \lt( \frc{d}{dx} \rt)^{n-1}
    \alpha^\nu(x).
\eeq

One crucial role of the formal operator $\Delta_x$ is to make the
fixed--point subalgebra $\{u\,|\, \nu u=u\}$ act according to a
true module action.  This property will be essential below in
constructing a representation of the algebra $\h{D}^+$ on the
twisted space $S[\nu]$. For instance, the conformal vector
$\omega$ is in the fixed point subalgebra, so that the vertex
operator $Y_{S[\nu]}(\omega,x)$ generates a representation of the
Virasoro algebra on the space $S[\nu]$. This representation of the
Virasoro algebra was explicitly constructed in \cite{DL2}. As one
can see in the results of \cite{DL2} and as will become clear
below, the resulting representation of the Virasoro generator
$L_0$ is not an (infinite) sum of normal-ordered products the type
$\sum_{n\in\frc1p\Z}\,\:\alpha(n)\beta(-n)\:\;$; rather, there is
an extra term proportional to the identity on $S[\nu]$, the
so-called correction term, which appears because of the operator
$\Delta_x$. The correction term was calculated in \cite{DL2} using
the explicit action of $e^{\Delta_x}$ on $\omega$. In the case of
the period--2, $\nu=-1$ automorphism, this action is given by
\cite{FLM2}, \cite{FLM3}:
$$
    e^{\Delta_x}\omega = \omega + \frc1{16} (\dim h) x^{-2},
$$
and for general automorphism, the calculation was carried out in
\cite{DL2} (see also \cite{FFR} and \cite{FLM3}). This is
relevant, for instance, in the construction of the moonshine
module \cite{FLM3}.

In order to have the correction term for the representation of the
algebra $\h{D}^+$ on the twisted space $S[\nu]$, one can calculate
the action of $e^{\Delta_x}$ for a general automorphism on the
vectors generating the representation of the whole algebra
$\h\D^+$. This is a complicated problem, mainly because generators
of $\hat{\mathcal{D}}^+$ have arbitrary large weights. Below we
will calculate the correction terms using the general theory of
twisted modules for vertex operator algebras, in particular using
the modified weak associativity relation for twisted operators, as
well as the simple result (\ref{tvertexalphax}). Hence in our
argument, the explicit action of $\Delta_x$ on vectors generating
the representation of the algebra $\h\D^+$ is not of importance;
all we need to know is that {\em there exists} such an operator
$\Delta_x$ giving to the space $S[\nu]$ the properties of a
twisted module for the vertex operator algebra $S$.

\ssect{Construction of the corresponding representation of $\h\D^+$}

In order to construct a representation of $\h\D^+$ on $S[\nu]$ we
need to consider certain homogeneous twisted vertex operators
associated to $Y_{S[\nu]}(\cdot,x)$. For $\alpha \in \a{h}$ we
define the following series acting on $S[\nu]$:
\beq\label{alphanuX}
    \alpha^{\nu}\< x\>=X_{S[\nu]}(\alpha(-1){\bf 1},x)=
    \sum_{n \in \frac{1}{p}\mathbb{Z}} \alpha^\nu(n)x^{-n}
\eeq
where we used (\ref{tvertexalphax}) and the fact that
$\alpha(-1){\bf 1}$ has weight 1. Recalling the orthonormal basis
$\{\alpha_q|q=1,\ldots,d\}$ of $\goth{h}$, we define the following
two formal series acting on $S[\nu]$:
\beqa \label{Lnuy1y2x}
    L^{\nu;y_1,y_2}\<x\>&=& \frc12 \sum_{q=1}^{d}
    \: \alpha^\nu_q\langle e^{y_1}x\rangle \alpha^\nu_q\langle
        e^{y_2}x\rangle  \: \no \\
    && - \frc12 \frc{\d}{\d y_1} \lt(\sum_{k=0}^{p-1}
    \frac{(e^{\frac{k(-y_1+y_2)}{p}}
     - 1)\; {\rm dim} \ \goth{h}_{(k)}}
    {1-e^{-y_1+y_2}} \rt)
\eeqa
and
\beqa \label{bLnuy1y2x}
    {\bar{L}}^{\nu;y_1,y_2}\<x\>&=& \frc12 \sum_{q=1}^{d}
    \: \alpha^\nu_q\langle e^{y_1}x\rangle
    \alpha^\nu_q\langle e^{y_2}x\rangle  \: \no \\
    && - \frc12 \frc{\d}{\d y_1} \lt(\sum_{k=0}^{p-1}
        \frac{e^{\frac{k(-y_1+y_2)}{p}}
    {\rm dim} \ \goth{h}_{(k)}}
    {1-e^{-y_1+y_2}} \rt).
\eeqa
\begin{remark} \label{remaLuntwisted}
In the special case $p=1$ and $d=1$, the operators
$L^{\nu;y_1,y_2}\<x\>$ and ${\bar L}^{\nu;y_1,y_2}\<x\>$,
respectively, specialize to the operators $L^{(y_1,y_2)}(x)$ and
${\bar L}^{(y_1,y_2)}(x)$ of \cite{L4}, \cite{L5}.
\end{remark}

We now prove an important statement for our construction of the
representation of $\h\D^+$ on the twisted space $S[\nu]$:

\begin{propo} \label{mathcliter}
The series ${\bar{L}}^{\nu;y_1,y_2}\<x\>$ as defined in
(\ref{bLnuy1y2x}) can be identified with the following iterate of
vertex operators:
\beq \label{iterbar}
    {\bar{L}}^{\nu;y_1,y_2}\<x\>= X_{S[\nu]}
    \lt(\frac{1}{2} \sum_{q=1}^d Y[\alpha_q(-1) {\bf 1} ,y_1-y_2]
    \alpha_q(-1) {\bf 1} ,e^{y_2}x \rt).
\eeq
\end{propo}

\proof One first rewrites the formal series (\ref{bLnuy1y2x}) in
the form:
\beq\label{Lbarlimit}
    \b{L}^{\nu;y_1,y_2}\<x_2\>
    =
    \frc12 \lim_{x_1 \to  x_2}
    \sum_{q=1}^d \lt(\lt( \frac{ \frac{x_1}{x_2}e^{y_1-y_2} -1 }
    {e^{y_1-y_2} -1} \rt)^k
        \alpha^\nu_q \< e^{y_1}x_1\> \alpha^\nu_q\< e^{y_2}x_2\> \rt)
\eeq
for any fixed $k\in\N,\,k\ge2$, which comes {}from
\[
    \sum_{q=1}^{d}
    \alpha^\nu_q \< e^{y_1}x_1\> \alpha^\nu_q\<e^{y_2}x_2\>
    = \sum_{q=1}^{d} \:\alpha^\nu_q \< e^{y_1}x_1\>
    \alpha^\nu_q\<e^{y_2}x_2\>\: -
    \frc{\d}{\d y_1} \lt(
    \sum_{k=0}^{p-1} \frc{e^{\frc{k(-y_1+y_2)}p} \; {\rm dim} \
    \goth{h}_{(k)}}{
    1-\frc{x_2}{x_1}e^{-y_1+y_2}} \rt).
\]
Then one uses modified weak associativity (\ref{modwassht}) with
$s=0,\; y = y_1-y_2$ and the replacements $x_2\mapsto e^{y_2}x_2,\;
x_1 \mapsto e^{y_1}x_1$, along with the definition (\ref{alphanuX}).
\eproof

\begin{remark}
{}From the expression (\ref{Lbarlimit}) of the operator
$\b{L}^{\nu;y_1,y_2}(x)$, we see that it can be written in terms
of the normal ordering introduced in (\ref{ppX}):
\beq
    \b{L}^{\nu;y_1,y_2}\<x\>
    = \frc12 \pp \sum_{q=1}^d
    \alpha^\nu_q \< e^{y_1}x\> \alpha^\nu_q\<e^{y_2}x\>
    \pp\;.
\eeq
This generalizes formula (3.17) of \cite{L4} (equivalently,
formula (1.42) of \cite{L5}), and the Proposition above in
particular proves formulas (2.25), (2.26) of \cite{L5}.
\end{remark}

Once Proposition \ref{mathcliter} is established, we can use the
general theory of twisted modules for vertex operator algebras in
order to easily extend some results in the untwisted setting to
the twisted setting. This is the way in which we chose to prove,
next, the commutator formula for the formal series
$\b{L}^{\nu;y_1,y_2}\<x\>$ acting on $S[\nu]$, Proposition
\ref{propLbarbracket}. The corresponding commutator formula for
the untwisted operators ${\bar L}^{(y_1,y_2)}(x)$ was announced in
\cite{L4} (a proof was given in \cite{M2}). Below, we give an
alternative proof of this untwisted commutator formula (the
untwisted case of Proposition \ref{propLbarbracket}) by
specializing to the untwisted setting our general commutator
formula (\ref{gencomform}). Then we will use Proposition
\ref{mathcliter} to extend the proof to the twisted setting.

\begin{propo} \label{propLbarbracket}
The series ${\bar{L}}^{\nu;y_1,y_2}\<x\>$ as defined by
(\ref{bLnuy1y2x}) satisfies the following bracket relation:
\begin{eqnarray}\label{Lbarbracketsalg}
    \lefteqn{[{\bar L}^{\nu;y_1,y_2}\<x_1\>,{\bar
    L}^{\nu;y_3,y_4}\<x_2\>]} \\
    &&= - {\frac{1}{2}} \frac{\partial}{\partial y_1} \biggl({\bar
    L}^{\nu;-y_1+y_2+y_3,y_4}\<x_2\> \delta
    \left({\frac{e^{y_1}x_1}{e^{y_3}x_2}}\right)+ {\bar
    L}^{\nu;-y_1+y_2+y_4,y_3}\<x_2\>
    \delta \left({\frac{e^{y_1}x_1}{e^{y_4}x_2}}\right)\biggr)\nonumber\\
    &&\quad - {\frac{1}{2}} \frac{\partial}{\partial y_2} \biggl({\bar
    L}^{\nu;y_1-y_2+y_3,y_4}\<x_2\> \delta
    \left({\frac{e^{y_2}x_1}{e^{y_3}x_2}}\right)+ {\bar
    L}^{\nu;y_1-y_2+y_4,y_3}\<x_2\> \delta
    \left({\frac{e^{y_2}x_1}{e^{y_4}x_2}}\right)\biggr). \no
\end{eqnarray}
\end{propo}

\proof Let us first prove this formula in the untwisted case $p=1$.
We specialize our general commutator formula
(\ref{gencomform}) to the case $p=1$, with $u_2 = v_2
=\alpha_{q}(-1){\bf 1},\; u_1 = v_1 = \alpha_{q'}(-1){\bf 1}$ and
the replacements $y_2 \mapsto y_3-y_4,\; y_1 \mapsto y_1-y_2,\; x_2 \mapsto
e^{y_4}x_2,\; x_1 \mapsto e^{y_2} x_1$. We sum independently over $q$
and $q'$, {}from 1 to $d$, and multiply through by a factor of
$1/4$. We recall that here $\alpha_q,\; q=1,\ldots,d$, form an
orthonormal basis for $\a{h}$. In fact, in order to directly
obtain the form of the right--hand side as written in
(\ref{Lbarbracketsalg}), it is preferable to use Equation
(\ref{comder4}) instead of Equation (\ref{gencomform}), although
they are equivalent. Also, an important formula for our purposes is
\[
    Y[\alpha(-1){\bf 1},y]\beta(-1){\bf 1} = \<\alpha,\beta\>
    y^{-2}{\bf 1} + \mbox{ series in nonnegative powers of $y$}
\]
for $\alpha,\beta\,\in \a{h}$. Applying this to the first term on
the right--hand side of the commutator formula (\ref{comder4}), we
find the term
\beqa
    && \frc14 \sum_{q=1}^d \Res_y \bigg\{ y^{-2}\; \delta\lt(e^{y_2-y_3-y}
    \frc{x_1}{x_2} \rt)\; \cdot \nn
    && \qquad\qquad \cdot\; X_{S}( Y[\alpha_q(-1){\bf 1} ,y_1-y_2]
    Y[\alpha_q(-1){\bf
    1},-y_3+y_4-y] {\bf 1}, e^{y_2} x_1) \bigg\}. \no
\eeqa
Recall that skew--symmetry of Zhu's vertex operators (\ref{skewsymmZhu})
gives
in general $Y[u,y] {\bf 1} = e^{L[-1]y} u$, and that the
$L[-1]$--derivative property and the $L[-1]$--bracket property of
Zhu's vertex operators (\ref{Lm1derZhu}) and (\ref{Lm1bracketZhu})
give $e^{L[-1]y_2} Y[u,y_1] e^{-L[-1]y_2} = Y[u,y_1+y_2]$. We use
these two properties along with the $L[-1]$--derivative property
of homogeneous vertex operators (\ref{Lm1X}) to write the
expression above in the form
\beqa
    && \frc14 \sum_{q=1}^d \Res_y \bigg\{ y^{-2}\;
    \delta\lt(e^{y_2-y_3-y} \frc{x_1}{x_2}
    \rt)\; \cdot \nn
    && \qquad\qquad \cdot\;
    X_{S}( Y[\alpha_q(-1){\bf 1} ,y_1-y_2+y_3-y_4+y] \alpha_q(-1){\bf
    1}, e^{y_2-y_3+y_4-y} x_1) \bigg\}. \no
\eeqa
Using the main delta--function property and the fact that $\Res_y
(y^{-2} f(y-y_2)) = -\frc{\d}{\d y_2} f(-y_2)$ for $f(y)$ a formal
series in nonnegative powers of $y$, the expression above gives
the first term in the second parentheses on the right--hand side
of (\ref{Lbarbracketsalg}). Similar arguments applied to the
second term on the right--hand side of the commutator formula
(\ref{comder4}) lead to the second term in the second parentheses
on the right--hand side of (\ref{Lbarbracketsalg}).

For the third term on the right--hand side of the commutator
formula (\ref{comder4}), one first uses (\ref{XYskewsymmt}) to
bring it to the form
\[
    \delta\lt(\om^{r}\lt(e^{y_1-y_2-y}\frc{x_1}{x_2}\rt)^{1/p}\rt)
    X_M(Y[v_1,-y_1+y]Y[\nu^{-r}v_2,-y_2]Y[u_1,y]\nu^{-r}u_2,e^{y_1-y}
    x_1).
\]
Then similar arguments as above applied to this expression lead to
the first term in the first parentheses on the right--hand side of
(\ref{Lbarbracketsalg}). In a similar way, the last term on the
right--hand side of the commutator formula (\ref{comder4}) gives
the second term in the first parentheses on the right--hand side
of (\ref{Lbarbracketsalg}).

The same arguments could be used to prove the formula
(\ref{Lbarbracketsalg}) in the twisted case. Instead, we will use
a general and simple argument {}from the theory of vertex operator
algebras. First, it is easy to see that the formal series
$\sum_{q=1}^d Y[\alpha_q(-1) {\bf 1} ,y_1-y_2] \alpha_q(-1) {\bf
1}$, with coefficients in $S$, appearing in the vector argument of
the homogeneous twisted vertex operator $X_{S[\nu]}(\cdot,x)$ in
(\ref{iterbar}), is invariant under the automorphism $\nu$. Hence
by the properties of twisted modules for vertex operator algebras,
the space $S[\nu]$ is a true module for the algebra satisfied by
the particular elements of $S$ generated in $y_1$ and $y_2$ by the formal
series $\sum_{q=1}^d Y[\alpha_q(-1) {\bf 1} ,y_1-y_2] \alpha_q(-1)
{\bf 1}$. The homogeneous twisted vertex operator
$X_{S[\nu]}(\cdot,x)$ gives actions of these elements on $S[\nu]$,
which are then in agreement with their actions on $S$ itself
for all twisting automorphisms. Hence the form of the bracket
relations for the formal series $\b{L}^{\nu;y_1,y_2}\<x\>$ acting
on $S[\nu]$ is independent of the twisting automorphism $\nu$,
which, combined with the proof of (\ref{Lbarbracketsalg}) in the
untwisted case above, proves (\ref{Lbarbracketsalg}) in the
twisted case. \eproof

Our twisted construction of $\hat{\mathcal{D}}^+$ is then a simple
consequence of the fact that the form of the previous commutator
formula is independent of the twisting automorphism $\nu$, which,
as emphasized in the proof above, is due to Proposition
\ref{mathcliter} and to aspects of the general theory of twisted
modules for vertex operator algebras. Moreover, as we will see,
Proposition \ref{propLbarbracket} immediately leads to the main
properties of the generators $\b{L}_n^{(r)}$ (see Section
\ref{sec:untwistedconstruction}) of the algebra $\h\D^+$, in
particular to the monomial central term in the bracket relations.

The formal series $\b{L}^{\nu;y_1,y_2}\<x\>$ (\ref{bLnuy1y2x})
generates a representation of the Lie algebra
$\hat{\mathcal{D}}^+$ on the $\hat{\a{h}}[\nu]$-module $S[\nu]$,
generalizing the untwisted case studied in \cite{Bl}. In fact,
combined with what we said in Remark \ref{remacylcoord},
Proposition \ref{mathcliter} suggests that
$\b{L}^{\nu;y_1,y_2}\<x\>$ should be a generating function for
modes that generalize the Virasoro modes in ``cylindrical
coordinates'' as defined after equation (\ref{Tcylcoord}). More
precisely, let
\beqa\label{Lnuy1y2exp}
    {L}^{\nu;y_1,y_2}\<x\>&=&
        \sum_{n\in\Z,\;r_1,r_2\,\in\N} {L}^{\nu;r_1,r_2}(n) x^{-n}
        \frc{y_1^{r_1}y_2^{r_2}}{r_1!r_2!}, \\
    \b{L}^{\nu;y_1,y_2}\<x\>&=& \frc{1}2 \frc{d}{(y_1-y_2)^2} +
        \sum_{n\in\Z,\;r_1,r_2\,\in\N} \b{L}^{\nu;r_1,r_2}(n) x^{-n}
        \frc{y_1^{r_1}y_2^{r_2}}{r_1!r_2!}.
    \label{bLnuy1y2exp}
\eeqa
Then the following holds (recall the generators (\ref{Lnrdef}),
(\ref{bLnrbloch}) of $\hat{\mathcal{D}}^+$):
\begin{theo} \label{main1}
Let
\beqa
    {L}^{\nu;r}(n) &=& {L}^{\nu; r,r}(n) \com{n\in\Z,\,r\in\N}, \no \\
    \b{L}^{\nu;r}(n) &=& \b{L}^{\nu;r,r}(n) \com{n\in\Z,\,r\in\N}. \no
\eeqa
\begin{itemize}
\item[(a)] The assignment
$$ L_n^{(r)} \mapsto L^{\nu;r}(n), \ \ c \mapsto d,$$
defines a representation of the Lie algebra
$\hat{\mathcal{D}}^+$ on $S[\nu]$.
\item[(b)]
The assignment
$$\bar{L}_n^{(r)} \mapsto \bar{L}^{\nu;r}(n), \ \ c \mapsto d $$
also defines a representation of the Lie algebra
$\hat{\mathcal{D}}^+$, with the central term being a pure
monomial, as in (\ref{bl2coc}).
\end{itemize}
\end{theo}

\proof We first prove assertion (a). The formal series
$L^{\nu;y_1,y_2}\<x\>$ defined by (\ref{Lnuy1y2x}) specializes, in
the untwisted case $p=1$, to the following
series acting on $S$ (see also Remark \ref{remaLuntwisted}):
\[
   L^{(y_1,y_2)}\<x\> = \frc12 \sum_{q=1}^{d}
    \: \alpha_q\langle e^{y_1}x\rangle \alpha_q\langle e^{y_2}x\rangle\:
\]
where $\alpha\<x\> = x\alpha(x)$ for $\alpha\in\a{h}$. It is easy
to see that an expansion in $y_1$ and $y_2$ as above:
\[
    {L}^{(y_1,y_2)}\<x\>=
    \sum_{n\in\Z,\;r_1,r_2\,\in\N} {L}^{(r_1,r_2)}(n) x^{-n}
    \frc{y_1^{r_1}y_2^{r_2}}{r_1!r_2!},
\]
leads to operators $L^{(r)}(n) = L^{(r,r)}(n)$ that represent on
$S$ generators $L_n^{(r)}$ of the algebra $\h\D^+$, as in
(\ref{Lnr}).

The other operators $L^{(r_1,r_2)}(n),\; r_1\neq r_2$ are linear
combinations of the operators $L^{(r)}(n)$. Indeed, they can be
written
\beq\label{Lr1r2n}
    L^{(r_1,r_2)}(n) = \frc12\sum_{q=1}^d \sum_{j\in\Z}
    j^{r_1}(n-j)^{r_2}\,\: \alpha_q(j)\alpha_q(n-j) \:
    \com{n\in\Z,\,r_1,r_2\,\in\N}.
\eeq
To see that these are linear combinations of $L^{(r)}_n$, note
that because the normal ordering is symmetric, the coefficients of
the operators $\: \alpha_q(j)\alpha_q(n-j) \:$ can be symmetrized:
\[
    L^{(r_1,r_2)}(n) = \frc12\sum_{q=1}^d \sum_{j\in\Z}
    ( j^{r_1}(n-j)^{r_2} + (n-j)^{r_1}(j)^{r_2} )\,\:
\alpha_q(j)\alpha_q(n-j)
    \:\;.
\]
In general, any symmetric polynomial in two variables $z_1,\,z_2$
can be written in a unique way as a polynomial in $z_1+z_2$ and
$z_1z_2$. Hence, the symmetric polynomial
\[
    \frc12\lt(j^{r_1}(n-j)^{r_2} + (n-j)^{r_1}(j)^{r_2}\rt)
\]
in $j,\,n-j$ can be written in a unique way as a polynomial in $n$
and $j(n-j)$. The coefficient (which is itself a polynomial in
$n$) of $(j(n-j))^r$ in this polynomial is the coefficient of
$L^{(r)}(n)$ in the linear combination representing
$L^{(r_1,r_2)}(n)$. That is, denoting this coefficient by
$C^{(r_1,r_2|r)}_n$, we have
\beq\label{lincombL}
    L^{(r_1,r_2)}(n) = \sum_{0\le r \le (r_1+r_2)/2}
    C^{(r_1,r_2|r)}_n L^{(r)}(n).
\eeq
In fact, the operators $L^{(r,r)}(n)$ for all $r\in\N,\;n\in\Z$
form a basis in the linear space spanned by $L^{(r_1,r_2)}(n)$ for
all $r_1,r_2\,\in\N,\;n\in\Z$. We will not need the explicit form
of these coefficients, except for $C^{(r_1,r_2|r)}_0$, which can
be found easily:
\beq\label{coeff0}
    C^{(r_1,r_2|r)}_0 = (-1)^{r_1} \delta_{r_1+r_2,2r}.
\eeq

In the twisted case, the operators $L^{\nu;r_1,r_2}(n)$ are also
linear combinations of operators $L^{\nu;r}(n)$, with the same
coefficients $C^{(r_1,r_2|r)}_n$. This is trivial for $n\neq0$.
For $n=0$, we only need to verify that the relation
(\ref{lincombL}) is still satisfied if one replaces the operators
$L^{(r_1,r_2)}(0)$ and $L^{(r)}(0) = L^{(r,r)}(0)$ by the
coefficients of $y_1^{r_1}y_2^{r_2}/(r_1!r_2!)$ and of
$(y_1y_2)^{r}/(r!)^2$, respectively, in the following formal power
series in $y_1-y_2$:
\[
    f(y_1-y_2) =- \frc12 \frc{\d}{\d y_1} \lt(\sum_{k=0}^{p-1}
    \frac{(e^{\frac{k(-y_1+y_2)}{p}}
     - 1)\; {\rm dim} \ \goth{h}_{(k)}}
    {1-e^{-y_1+y_2}} \rt),
\]
which appears in the definition of $L^{\nu;y_1,y_2}\<x\>$
(\ref{Lnuy1y2x}). This formal power series is even under change of
sign of its argument: $f(y) = f(-y)$. It is easy to see that this
indeed gives
\beq
    f(y_1-y_2) = \sum_{r_1,r_2\ge0} f^{(r_1,r_2)}
    \frc{y_1^{r_1}y_2^{r_2}}{r_1!r_2!}
\eeq
where the coefficients $f^{(r_1,r_2)}$ are of the form
\beq
    f^{(r_1,r_2)} = \sum_{r\ge0} C^{(r_1,r_2|r)}_0 f^{(r,r)}
\eeq
with (\ref{coeff0}).

Hence, the operators $L^{\nu;r_1,r_2}(n)$ are linear combinations of
operators $L^{\nu;r}(n)$ with the same coefficients
$C^{(r_1,r_2|r)}_n$ as those appearing in (\ref{lincombL}), the
form of the bracket relations (\ref{Lbarbracketsalg}) is independent
of the twisting automorphism, and the difference $\frc12
\frc{\d}{\d y_1} \frc{d}{1-e^{-y_1+y_2}}$ between the formal series
$L^{\nu;y_1,y_2}\<x\>$ and $\b{L}^{\nu;y_1,y_2}\<x\>$ is also
independent of the twisting automorphism. This implies that the operators
$L^{\nu;r}(n)$
satisfy, for all twisting automorphisms $\nu$, the same bracket
relations as do the operators $L^{(r)}(n)$ that generate the algebra
$\h\D^+$. This proves assertion (a).

In order to prove assertion (b), we only need to prove that the
central term in the commutator of generators $\b{L}^{\nu;r}(n)$ is
a monomial. This can be seen {}from the bracket formula
(\ref{Lbarbracketsalg}). The source of the central term in this
bracket formula is the term $\frc12 \frc{d}{(y_1-y_2)^2}$ in
(\ref{bLnuy1y2exp}). As expected, it appears on the right--hand
side of (\ref{Lbarbracketsalg}) only when the power of $x_1$ is
the same as that of $x_2^{-1}$. On both sides of this formula,
take then the term in $(x_1/x_2)^m$ for some fixed $m\in\Z$, and
fix the powers of $y_1,\, y_2,\, y_3$ and $y_4$. This selects, on
the left--hand side of (\ref{Lbarbracketsalg}), two elements of
the form $\b{L}^{\nu;r_1,r_2}(0)$ and $\b{L}^{\nu;r_3,r_4}(0)$ of
which we take the bracket. On the right--hand side, in the first
term, for instance, the part of $\b{L}^{\nu;-y_1+y_2+y_3,y_4}(0)$
relevant to the calculation of the central term, from the formula
(\ref{bLnuy1y2exp}), is the formal series $\frc12
\frc{d}{(-y_1+y_2+y_3-y_4)^2}$. This is a sum of terms in which
the sum of the powers of $y_1,\,y_2,\,y_3$ and $y_4$ is fixed to
-2. In (\ref{Lbarbracketsalg}), this formal series is multiplied
by the formal delta--function $\delta(e^{y_1-y_3}x_1/x_2)$, which
gives a sum of contributions of the form $m^k (y_1-y_3)^k$ for all
nonnegative values of $k$. In a contribution $(y_1-y_3)^k$, the
sum of the powers of $y_1$ and of $y_3$ is $k$. But since the
powers of $y_1,\,y_2,\,y_3$ and $y_4$ are fixed in
(\ref{Lbarbracketsalg}), the value of $k$ is also fixed.
Specifically, in the bracket of $\b{L}^{\nu;r_1,r_2}(0)$ with
$\b{L}^{\nu;r_3,r_4}(0)$, it is given by $k=r_1+r_2+r_3+r_4+3$.
The other terms on the right--hand side of (\ref{Lbarbracketsalg})
lead to the same value of $k$. Hence, only one power of $m$
appears; this proves assertion (b). \eproof

\begin{remark}
Equivalently, it is possible to represent the algebra $\h{\D}^+$
using the operators $L^{(r,0)}(n)$ instead of the operators
$L^{(r,r)}(n)$. Indeed, expanding $j^{r_1}(n-j)^{r_2}$ in a
polynomial in $j$ and $n$ in (\ref{Lr1r2n}), one can write any
operator $L^{(r_1,r_2)}(n)$ as a linear combination of
$L^{(r,0)}(n)$. {}From these linear combinations, one can see that
the operators $L^{(r,0)}(n)$ for all $r\in\N,\;n\in\Z$ form a
basis in the linear space spanned by $L^{(r_1,r_2)}(n)$ for all
$r_1,r_2\,\in\N,\;n\in\Z$.
\end{remark}

Explicit expressions for the operators $L^{\nu;r}(n)$ and
$\b{L}^{\nu;r}(n)$, involving Bernoulli polynomials, are easy to
obtain {}from (\ref{Lnuy1y2x}) and (\ref{bLnuy1y2x}):
\beqa
    L^{\nu;r}(n)&=& \frc12 \sum_{q=1}^{d} \sum_{j\in \frac{1}{p}\Z}
    j^{r}(n-j)^{r}\,\: \alpha_q^\nu (j)\alpha_q^\nu(n-j) \: \no\\   &&
    - \delta_{n,0}\, \frc{(-1)^r}{4(r+1)} \sum_{k=0}^{p-1}
    \dim\a{h}_{(k)}
        \lt( B_{2(r+1)}(k/p) - B_{2(r+1)} \rt)
\eeqa
and
\beqa
    \b{L}^{\nu;r}(n) &=& \frc12 \sum_{q=1}^d \sum_{j\in \frac{1}{p}
    \Z} j^{r}(n-j)^{r}\,
    \: \alpha_q^\nu (j)\alpha_q^\nu (n-j) \: \no\\   &&  -
    \delta_{n,0} \, \frc{(-1)^r}{4(r+1)} \sum_{k=0}^{p-1}
    \dim\a{h}_{(k)} B_{2(r+1)}(k/p).
\eeqa
{}From our construction, the appearance of Bernoulli polynomials
is seen to be directly related to general properties of
homogeneous twisted vertex operators.

The next result is a simple consequence of Theorem \ref{main1}. It
describes the action of the ``Cartan subalgebra'' of
$\hat{\mathcal{D}}^+$ on a highest weight vector of a canonical
quasi-finite $\hat{\mathcal{D}}^+$--module; here we are using the
terminology of \cite{KR}. This corollary gives the ``correction''
terms referred to in the introduction.

\begin{corol} \label{kacradul}
Given a highest weight $\hat{\mathcal{D}}^+$--module $W$, let
$\delta$ be the linear functional on the ``Cartan subalgebra'' of
$\hat{\mathcal{D}}^+$ (spanned by $L^{(k)}_{0}$ for $k\in \N$)
defined by
$$
    {L}^{(k)}_{0} \cdot w =(-1)^k\delta\lt( {L}^{(k)}_{0} \rt)w,
$$
where $w$ is a generating highest weight vector of $W$, and let
$\Delta(x)$ be the generating function
$$
    \Delta(x)=\sum_{k \ge 1} \frac{\delta({L}^{(k)}_0)x^{2k}}{(2k)!}
$$
(cf. \cite{KR}). Then for every automorphism $\nu$ of period $p$
as above,
$$
    \mathcal{U}(\hat{\mathcal{D}}^+) \cdot 1 \subset S[\nu]
$$
is a quasi--finite highest weight $\hat{\mathcal{D}}^+$--module
satisfying
\beq\label{Delta}
    \Delta(x) = \frac{1}{2} \frac{d}{dx} \sum_{k=0}^{p-1} \frac{
    (e^{\frac{kx}{p}} -1) {\rm dim} \ \goth{h}_{(k)}} {1-e^x}.
\eeq
\end{corol}

\proof Clearly $\mathcal{U}(\hat{\mathcal{D}}^+) \cdot 1 \subset
S[\nu]$ is a module for $\h\D^+$, submodule of $S[\nu]$. Its
highest weight vector is $1$, and the action of the algebra
element $L_0^{(k)}$ on this vector is given by the operator
$L^{\nu;k,k}(0)$ in the formal series (\ref{Lnuy1y2exp}) defined
by (\ref{Lnuy1y2x}). This action comes entirely {}from the term
proportional to the identity operator in $(\ref{Lnuy1y2x})$, which
immediately gives (\ref{Delta}). \eproof

Finally, we have an additional result (for the untwisted bosonic case
an equivalent result was obtained in \cite{Bl} and for the spinor
constructions in \cite{M2}):

\begin{corol}
The generating function
\beq \label{generators}
    X_{S[\nu]}(\sum_{q=1}^d \alpha_q(-m-1)\alpha_q(-m-1){\bf 1},x),
\eeq
$m \in \mathbb{N}$, defines the same $\hat{\mathcal{D}}^+$--module
as in Theorem \ref{main1}.  That is, every operator $L^{(r)}(n)$
(or equivalently $\bar{L}^{(r)}(n)$) can be expressed as a linear
combination of the expansion coefficients of the operator
(\ref{generators}), and vice versa.
\end{corol}

\proof Clearly, we can replace $X$--operators (\ref{generators})
by
$$
    Y_{S[\nu]}(\sum_{q=1}^d \alpha_q(-m-1)\alpha_q(-m-1){\bf 1},x).
$$
{}From \cite{DL2} we have
$$
    Y_{S[\nu]}(\alpha_q(-m-1)\alpha_q(-m-1){\bf 1},x)=
    W(e^{\Delta_x}\alpha_q(-m-1)\alpha_q(-m-1){\bf 1},x),
$$
where $\Delta_x$ is given by formula (4.42) in \cite{DL2} and
$$
    W(\alpha_q(-m-1)\alpha_q(-m-1){\bf 1},x)=\:
    \left(\frac{1}{m!} \left( \frac{d}{dx} \right)^m \alpha^\nu_q(x)
    \right)\left(\frac{1}{m!} \left( \frac{d}{dx} \right)^m
    \alpha_q^\nu (x) \right) \:\; .
$$
In addition it is not hard to see that
\beqa \label{wcorr}
    && W(e^{\Delta_x}\alpha_q(-m-1)\alpha_q(-m-1){\bf 1},x)=\nn &&
    W(\alpha_q(-m-1)\alpha_q(-m-1){\bf 1},x)+x^{-2m-2}f(q,m),
\eeqa
where $f(q,m) \in \mathbb{C}$. Consider
\beq \label{wqad}
    \sum_{q=1}^d W(\alpha_q(-m-1)\alpha_q(-m-1){\bf 1},x).
\eeq
As in the untwisted case \cite{M2}, it follows that the space
spanned by the expansion coefficients of (\ref{wqad}) defines the
same space of operators as the space spanned by the expansion
coefficients of
\beq \label{wqad1}
    \sum_{q=1}^{d} \: \alpha^\nu_q\langle e^{y_1}x\rangle
    \alpha^\nu_q\langle e^{y_2}x\rangle  \:\; .
\eeq
In other words every generator of $\hat{\mathcal{D}}^+$ of degree
$\neq 0$  (described in Theorem \ref{main1}) is a linear
combination of the Fourier coefficients in (\ref{generators}) and
vice versa. In addition, because of the twisted Jacobi identity
(\ref{jacobiXt}), the operators of the form (\ref{generators}) are
closed with respect to the commutator; therefore they generate a
Lie algebra. But the generators of nonzero degree uniquely
determine the action of the whole Lie algebra
$\hat{\mathcal{D}}^+$ (this fact follows by induction). The result
follows. \eproof

{\ }\\

\small{

{\em Acknowledgments.} B.D. gratefully acknowledges partial
support {}from an NSERC Postgraduate Scholarship. J.L. and A.M.
gratefully acknowledge partial support {}from NSF grant
DMS-0070800.

}

\noindent {\small \sc Department of Physics,  Rutgers University,
Piscataway,
NJ 08854} \\
\noindent {\small Current address: \sc Rudolf Peierls Centre for
Theoretical Physics, University of Oxford,
Oxford, OX1 3NP, UK} \\
{\em E-mail address}: b.doyon1@physics.ox.ac.uk \\
\\
\noindent {\small \sc Department of Mathematics,  Rutgers
University, Piscataway,
NJ 08854} \\
{\em E--mail address}: lepowsky@math.rutgers.edu \\
\\
\noindent {\small \sc Department of Mathematics,  University of
Arizona, Tucson, AZ 85721} \\
{\small Current address: \sc Department of Mathematics, University at
Albany, SUNY, Albany, NY 12222} \\
{\em E--mail address}: amilas@math.albany.edu


\begin{thebibliography}{DHVW2}

\bibitem[AFOQ]{AFOQ}
H. Awata, M. Fukuma, S. Odake and Y--H. Quano, Eigensystem and
full character formula of the $W_{1+\infty}$ algebra with $c=1$,
{\em  Lett. Math. Phys.} {\bf 31} (1994), 289--298.

\bibitem[AFMO]{AFMO}
H. Awata, M. Fukuma, Y. Matsuo and S. Odake, Representation theory
of the ${W}_{1+\infty}$ algebra, Quantum field theory, integrable
models and beyond (Kyoto, 1994), {\em Progr. Theoret. Phys.
Suppl.} {\bf 118} (1995),  343--373.

\bibitem[BDM]{BDM}
K. Barron, C. Dong and G. Mason, Twisted sectors for tensor
product vertex operator algebras associated to permutation groups,
{\em Comm. Math. Phys.} {\bf 227} (2002), 349--384.

\bibitem[BPZ]{BPZ}
A.~A.~Belavin, A.~M.~Polyakov and A.~B.~Zamolodchikov, Infinite
conformal symmetries in two-dimensional quantum field theory, {\em
Nucl. Phys.} {\bf B241} (1984), 333--380.

\bibitem[Bl]{Bl}
S.~Bloch,
Zeta values and differential operators on the circle,
{\em J. Algebra} {\bf 182} (1996), 476--500.

\bibitem[BO]{BO}
S. Bloch and A. Okounkov, The character of the infinite wedge
representation, {\em  Adv. Math.} {\bf 149} (2000), 1--60.

\bibitem[Bo]{Bo}
R.~E.~Borcherds,
Vertex algebras, Kac-Moody algebras, and the Monster,
{\em Proc. Natl. Acad. Sci. USA} {\bf 83} (1986), 3068--3071.

\bibitem[BHS]{BHS}
L. Borisov, M. B. Halpern and C. Schweigert, Systematic approach
to cyclic orbifolds, {\em Inter. J. of Mod. Phys.} {\bf A13}
(1998), 125--168.

\bibitem[dBHO]{dBHO}
J. de Boer, M. B. Halpern and N. A. Obers, The operator algebra and
twisted KZ equations of WZW orbifolds, {\em J. High Energy Phys.}
{\bf 10} (2001), 1.

\bibitem[DVVV]{DVVV}
R. Dijkgraaf, C. Vafa, E. Verlinde and H. Verlinde, Operator algebras
of orbifold models, {\em Commun. Math. Phys.} {\bf 123} (1989),
485--526.

\bibitem[DFMS]{DFMS}
L. Dixon, D. Friedan, E. Martinec and S. Shenker, The conformal field
theory of orbifolds {\em Nucl. Phys.} {\bf B282} (1987), 13--73.

\bibitem[DHVW1]{DHVW1}
L. Dixon, J. A. Harvey, C. Vafa and E. Witten, Strings on orbifolds, I
{\em Nucl. Phys.} {\bf B261} (1985), 678--686.

\bibitem[DHVW2]{DHVW2}
L. Dixon, J. A. Harvey, C. Vafa and E. Witten, Strings on orbifolds,
II {\em Nucl. Phys.} {\bf B274} (1986), 285--314.

\bibitem[D]{D}
C. Dong, Twisted modules for vertex algebras associated with even
lattices, {\em J. Algebra} {\bf 165} (1994), 91--112.

\bibitem[DL1]{DL1}
C.~Dong and J.~Lepowsky,
\newblock {\em Generalized Vertex Algebras and Relative Vertex
Operators},
\newblock Progress in Mathematics, Vol. 112, Birkh\"{a}user,
Boston, 1993.

\bibitem[DL2]{DL2}
C.~Dong and J.~Lepowsky,
\newblock The algebraic structure of relative twisted vertex
operators, \newblock {\em Journal of Pure and Applied Algebra}
{\bf 110} (1996), 259--295.

\bibitem[DLM]{DLM}
C. Dong, H. Li and G. Mason, Modular invariance of trace functions
in orbifold theory and generalized Moonshine, {\em Comm. Math.
Phys.} {\bf 214} (2000), 1--56.

\bibitem[DLMi]{DLMi}
B. Doyon, J. Lepowsky and A. Milas, Twisted modules for vertex
operator algebras and Bernoulli polynomials, {\em I.M.R.N.} {\bf
44} (2003), 2391-2408.

\bibitem[FFR]{FFR}
A. J. Feingold, I. B. Frenkel and J. F. X. Ries, {\em Spinor
construction of vertex operator algebras, triality and
$E_8^{(1)}$}, Contemporary Math., Vol. 121, Amer. Math. Soc.,
Providence, 1991.

\bibitem[FrB]{FrB}
E. Frenkel and D. Ben--Zvi, {\em Vertex Algebras and Algebraic
Curves}, Mathematical Surveys and Monographs, Vol. 88, Amer. Math.
Soc., Providence, 2001.

\bibitem[FrS]{FrS}
E. Frenkel and M. Szczesny, Twisted modules over vertex algebras on
algebraic curves, {\em Adv. Math.} {\bf 187} (2004), 195--227.

\bibitem[FHL]{FHL}
I.~B.~Frenkel, Y.-Z.~Huang and J.~Lepowsky,
On axiomatic approaches to vertex operator algebras
and modules, preprint (1989);
{\em Memoirs Amer. Math. Soc.} {\bf 104} (1993).

\bibitem[FKRW]{FKRW}
E. Frenkel, V. Kac, A. Radul and W. Wang, $W_{1+\infty}$ and
$W(gl_N)$ with central charge $N$, {\em Comm. Math. Phys.} {\bf
170} (1995), 337-357.

\bibitem[FLM1]{FLM1}
I. B. Frenkel, J. Lepowsky and A. Meurman, A natural
representation of the Fischer-Griess Monster with the modular
function $J$ as character, {\em Proc. Natl. Acad. Sci. USA} {\bf
81} (1984), 3256-3260.

\bibitem[FLM2]{FLM2}
I.~B. Frenkel, J.~Lepowsky and A.~Meurman, Vertex operator
calculus, in: {\em Mathematical Aspects of String Theory}, Proc.
1986 Conf., San Diego, ed. by S.-T. Yau,  World Scientific,
Singapore, 1987, 150-188.

\bibitem[FLM3]{FLM3}
I.~B. Frenkel, J.~Lepowsky and A.~Meurman,
{\em Vertex Operator Algebras and the Monster},
Pure and Appl. Math., Vol. 134, Academic Press, Boston, 1988.

\bibitem[GHHO]{GHHO}
O. Ganor, M. B. Halpern, C. Helfgott and N. A. Obers, The
outer-automorphic WZW orbifolds on ${\frak{so}}(2n)$, including five
triality orbifolds on ${\frak{so}}(8)$, {\em J. High Energy Phys.}
{\bf 212} (2002), 19.

\bibitem[G]{G}
P.~ Goddard,
Meromorphic conformal field theory, in:
{\em Infinite Dimensional Lie Algebras and Groups, Advanced
Series in Math. Physics,} Vol. 7, ed. by V.~Kac, World Scientific,
Singapore, 1989, 556--587.

\bibitem[HH]{HH}
M. B. Halpern and C. Helfgott, The general twisted open WZW string,
{\em Int. J. Mod. Phys.} {\bf A20} (2005), 923.

\bibitem[HO]{HO}
M. B. Halpern and N. A. Obers, Two large examples in orbifold theory:
abelian orbifolds and the charge conjugation orbifold on
${\frak{su}}(n)$, {\em Int. J. of Mod. Phys.} {\bf A17} (2002), 3897.

\bibitem[H1]{H1}
Y.-Z.~Huang, Applications of the geometric interpretation of
vertex operator algebras, in: {\em Proc. 20th International
Conference on Differential Geometric Methods in Theoretical
Physics, New York, 1991}, ed. by S.~Catto and A.~Rocha, World
Scientific, Singapore, 1992, 333--343.

\bibitem[H2]{H2}
Y.-Z. Huang,
{\em Two-dimensional Conformal Geometry and Vertex
Operator Algebras},
Progress in Math., Vol. 148, Birkh\"{a}user,
Boston, 1997.

\bibitem[KR]{KR}
V. Kac and A. Radul, Quasifinite highest weight modules over the
Lie algebra of differential operators on the circle, {\em Comm.
Math. Phys.} {\bf 157} (1993), 429--457.

\bibitem[KWY]{KWY}
V. G. Kac, W. Wang and C. H. Yan, Quasifinite representations of
classical Lie subalgebras of $\mathcal{W}_{1+\infty}$, {\em Adv.
Math.} {\bf 139} (1998), 56--140.

\bibitem[L1]{L1}
J.~Lepowsky, Calculus of twisted vertex operators, {\em Proc.
Natl. Acad. Sci. USA} {\bf 82} (1985), 8295--8299.

\bibitem[L2]{L2}
J.~Lepowsky, Perspectives on vertex operators and the Monster, in:
{\em Proceedings of the Symposium on the Mathematical Heritage of
Hermann Weyl}, Duke University, May, 1987, Proc. Symp. Pure Math., Vol.
48, Amer. Math. Soc., Providence, 1988, 181--197.

\bibitem[L3]{L3}
J.~Lepowsky,
Remarks on vertex operator algebras and moonshine, in:
{\em Proc. 20th International Conference on Differential
Geometric Methods in Theoretical Physics, New York, 1991}, ed. by
S.~Catto and A.~Rocha, World Scientific, Singapore, 1992,
362--370.

\bibitem[L4]{L4}
J.~Lepowsky, Vertex operator algebras and the zeta function, in: {\em
Recent Developments in Quantum Affine Algebras and Related
Topics}, ed. by N.~Jing and K.~C.~Misra, Contemporary Math., Vol.
248, Amer. Math. Soc., Providence, 1999, 327-340.

\bibitem[L5]{L5}
J.~Lepowsky, Application of a `Jacobi Identity' for vertex
operator algebras to zeta values and differential operators, {\em
Lett. Math. Phys.} {\bf 53} (2000), 87-103.

\bibitem[LL]{LL}
J. Lepowsky and H. Li, {\em Introduction to Vertex Operator
Algebras and Their Representations}, Progress in Mathematics, Vol. 227,
Birkh\"auser, Boston, 2003.

\bibitem[LW]{LW}
J. Lepowsky and R.L. Wilson, Construction of the affine Lie
algebra $A_{1}^{(1)}$, {\em Comm. Math. Phys.} {\bf 62}
(1978), 43--53.

\bibitem[Li1]{Li1}
H. Li, Local systems of vertex operators, vertex superalgebras and
modules, {\em Journal of Pure and Applied Algebra} {\bf 109} (1996),
143--195.

\bibitem[Li2]{Li2}
H. Li, Local systems of twisted vertex operators, vertex operator
superalgebras and twisted modules, in: {\em Moonshine, the
Monster, and Related Topics (South Hadley, MA, 1994)},
Contemporary Math., Vol. 193, Amer. Math. Soc., Providence, 1996,
203--236.

\bibitem[M1]{M1}
A. Milas, Correlation functions, vertex operator algebras and
$\zeta$--functions, Ph.D. thesis, Rutgers University, 2001.

\bibitem[M2]{M2}
A. Milas, Formal differential operators, vertex operator algebras
and zeta--values, I, {\em Journal of Pure and Applied Algebra}
{\bf 183} (2003), 129-190.

\bibitem[M3]{M3}
A. Milas, Formal differential operators, vertex operator algebras
and zeta--values, II, {\em Journal of Pure and Applied Algebra}
{\bf 183} (2003), 191-244.

\bibitem[N]{N}
S.--H. Ng, The Lie bialgebra structures on the Witt and Virasoro algebras,
Ph.D. thesis, Rutgers University, 1997.

\bibitem[OP]{OP}
A. Okounkov and R. Pandharipande, Gromov--Witten theory, Hurwitz
theory, and completed cycles, arXiv:math.AG/0204305.

\bibitem[Z1]{Z1}
Y.~Zhu,
Vertex operators, elliptic functions and modular forms,
Ph.D. thesis, Yale University, 1990.

\bibitem[Z2]{Z2}
Y.~Zhu,
Modular invariance of characters of vertex operator algebras,
{\em J. Amer. Math. Soc.} {\bf 9} (1996), 237--307.

\end{thebibliography}
\end{document}